\pdfoutput=1
\RequirePackage{ifpdf}
\ifpdf 
\documentclass[pdftex]{sigma}
\else
\documentclass{sigma}
\fi

\usepackage{bbm,tikz,stmaryrd}
\usepackage[OT2,T1]{fontenc}
\DeclareSymbolFont{cyrletters}{OT2}{wncyr}{m}{n}
\DeclareMathSymbol{\Sha}{\mathalpha}{cyrletters}{"58}

\numberwithin{equation}{section}

\newtheorem{Theorem}{Theorem}[section]
\newtheorem*{Theorem*}{Theorem}

\newtheorem{Conjecture}[Theorem]{Conjecture}
\newtheorem*{Conjecture*}{Conjecture}
\newtheorem{Lemma}[Theorem]{Lemma}
\newtheorem{Proposition}[Theorem]{Proposition}
 { \theoremstyle{definition}
\newtheorem{Definition}[Theorem]{Definition}

\newtheorem{Example}[Theorem]{Example}
 }

\newcommand{\red}[1]{{\color{red} #1}}
\newcommand{\green}[1]{{\color{green} #1}}

\def\BZ{\mathbb{Z}}
\def\BQ{\mathbb{Q}}
\def\BR{\mathbb{R}}
\def\BC{\mathbb{C}}

\def\SL{\mathrm{SL}}

\def\tq{\tilde{q}}

\def\Res{\mathrm{Res}}

\def\z{\zeta}

\def\ae{\sim}

\def\be{\begin{equation}}
\def\ee{\end{equation}}
\def\ba{\begin{aligned}}
\def\ea{\end{aligned}}

\def\Li{\mathrm{Li}}

\def\e{\mathbf e}

\def\Phih{\widehat{\Phi}}
\def\VC{\mathrm{VC}}

\def\PSL{\mathrm{PSL}}
\def\SL{\mathrm{SL}}

\def\GL{\mathrm{GL}}

\def\UHP{\mathfrak{h}}

\def\Bor{\mathcal{B}}
\def\Lap{\mathcal{L}}

\def\denom{\mathrm{denom}}
\def\numer{\mathrm{numer}}
\def\CS{\mathrm{CS}}

\def\SU{\mathrm{SU}}

\def\WRT{\Sha}

\def\CYCDS{D}
\def\etaroots{\varepsilon}

\def\ms{\mathbf{j}}
\def\Om{\Omega}

\def\BZ{\mathbb Z}
\def\BQ{\mathbb Q}
\def\BR{\mathbb R}
\def\BC{\mathbb C}

\def\calS{\mathcal S}

\def\SL{\mathrm{SL}}

\def\tq{\tilde{q}}
\def\Res{\mathrm{Res}}

\def\ve{\varepsilon}

\def\Li{\operatorname{Li}}
\def\PSL{\mathrm{PSL}}

\def\Re{\operatorname{Re}}
\def\Im{\operatorname{Im}}

\def\be{\begin{equation}}
\def\ee{\end{equation}}

\def\Bor{\mathcal{B}}
\def\Lap{\mathcal{L}}

\def\z{\zeta}

\def\GL{\mathrm{GL}}

\def\VC{\mathrm{V}}
\def\triv{\rho_0}
\def\Om{\Omega}

\usepackage{accents}

\def\SU{\mathrm{SU}}

\usepackage{listings}

\usepackage[textsize=footnotesize]{todonotes}
\definecolor{codegreen}{rgb}{0,0.6,0}
\definecolor{codegray}{rgb}{0.5,0.5,0.5}
\definecolor{codepurple}{rgb}{0.58,0,0.82}
\definecolor{backcolour}{rgb}{0.95,0.95,0.92}

\lstdefinelanguage{PARIGP}{
 keywords={typeof, new, true, false, catch, function, return, null, catch, switch, var, if, in, while, do, else, case, break},
 keywordstyle=\color{blue}\bfseries,
 ndkeywords={class, export, boolean, throw, implements, import, this},
 ndkeywordstyle=\color{darkgray}\bfseries,
 identifierstyle=\color{black},
 sensitive=false,
 comment=[l]{//},
 morecomment=[s]{/*}{*/},
 commentstyle=\color{purple}\ttfamily,
 stringstyle=\color{red}\ttfamily,
 morestring=[b]',
 morestring=[b]"
}

\lstdefinestyle{code}{
 backgroundcolor=\color{backcolour},
 commentstyle=\color{codegreen},
 keywordstyle=\color{magenta},
 numberstyle=\tiny\color{codegray},
 stringstyle=\color{codepurple},
 basicstyle=\ttfamily\footnotesize,
 breakatwhitespace=false,
 breaklines=true,
 captionpos=b,
 keepspaces=true,
 numbers=left,
 numbersep=5pt,
 showspaces=false,
 showstringspaces=false,
 showtabs=false,
 tabsize=2
}

\begin{document}
\allowdisplaybreaks

\newcommand{\arXivNumber}{2308.03265}

\renewcommand{\PaperNumber}{004}

\FirstPageHeading

\ShortArticleName{Quantum Modularity for a Closed Hyperbolic 3-Manifold}

\ArticleName{Quantum Modularity\\ for a Closed Hyperbolic 3-Manifold}

\Author{Campbell WHEELER}

\AuthorNameForHeading{C.~Wheeler}

\Address{Institut des Hautes \'Etudes Scientifiques, Le Bois-Marie, Bures-sur-Yvette, France}
\Email{\href{mailto:wheeler@ihes.fr}{wheeler@ihes.fr}}
\URLaddress{\url{https://www.ihes.fr/~wheeler}}

\ArticleDates{Received January 11, 2024, in final form December 23, 2024; Published online January 08, 2025}

\Abstract{This paper proves quantum modularity of both functions from $\mathbb{Q}$ and $q$-series associated to the closed manifold obtained by $-\frac{1}{2}$ surgery on the figure-eight knot, $4_1(-1,2)$. In a sense, this is a companion to work of Garoufalidis--Zagier, where similar statements were studied in detail for some simple knots. It is shown that quantum modularity for closed manifolds provides a unification of Chen--Yang's volume conjecture with Witten's asymptotic expansion conjecture. Additionally we show that $4_1(-1,2)$ is a counterexample to previous conjectures of Gukov--Manolescu relating the Witten--Reshetikhin--Turaev invariant and the $\widehat{Z}(q)$ series. This could be reformulated in terms of a ``strange identity'', which gives a volume conjecture for the $\widehat{Z}$ invariant. Using factorisation of state integrals, we give conjectural but precise $q$-hypergeometric formulae for generating series of Stokes constants of this manifold. We find that the generating series of Stokes constants is related to the 3d index of $4_1(-1,2)$ proposed by Gang--Yonekura. This extends the equivalent conjecture of Garoufalidis--Gu--Mari\~no for knots to closed manifolds. This work appeared in a similar form in the author's Ph.D.~Thesis.}

\Keywords{3d index; asymptotic expansions; Borel resummation; character varieties; Chern--Simons invariants; circle method; closed three-manifolds; cocycles; dilogarithm; duality; Faddeev quantum dilogarithm; factorisation; flat connections; hyperbolic manifolds; modularity; perturbative invariants; $q$-difference equations; $q$-hypergeometric functions; quadratic relations; quantum invariants; quantum modular forms; resurgence; surgery; state integrals; stationary phase; Stokes constants; Stokes phenomenon; strange identity; three-manifolds; volume conjecture; Witten--Reshetikhin--Turaev invariants; $\widehat{Z}$ invariants}

\Classification{57N10; 57K16; 57K14; 57K10}

\tableofcontents

\section{Introduction}

This paper is an in-depth case study of the $\mathfrak{sl}_{2}$-quantum invariants~\cite{Jones,RT:3man,Turaev:book,Witten} of the closed hyperbolic manifold obtained from $-1/2$ surgery on the figure-eight knot, $4_{1}(-1,2)$. This example serves to unify a variety of phenomenon under the umbrella of quantum modularity~\cite{Zagier:Qmod}. The~main importance of this example is the fact it is a hyperbolic manifold. For non-hyperbolic manifolds -- closed or not -- it has long been understood that their quantum invariants are related to Ramanujan's mock modular forms~\cite{AM:res,Cheng:3dmod,Hikami,LawrenceZag,Nahm:Blochgrp,Zagier:Vas}. The story for hyperbolic manifolds starts with Kashaev's volume conjecture~\cite{Kashaev:VC,Murakami2}. This required the introduction of a new concept called a quantum modular form~\cite{Zagier:Qmod}. The ideas surrounding this quantum modularity have been clarified in recent years~\cite{BetDrap,Cheng:3dmod,GK:rat,GK:qser,GZ:qser,GZ:RQMOD,Zagier:qmod.hol.lec} resulting in a detailed understanding for simple hyperbolic knots. This work provides the same depth for a closed hyperbolic manifold. For knots, the invariant is the Kashaev invariant or the coloured Jones polynomial. While here, the invariant is the Witten--Reshetikhin--Turaev (WRT) invariant. These invariants give the story at $\BQ$, however, recently this has been extended into $\BC\smallsetminus \BR$. To extend away from $\BQ$, we need the new $q$-series invariants of three-manifolds~\cite{BeemDimPas,Cheng:3dmod,GZ:qser,GukovMan,GPPV,GPV} that have been conjectured to exist and constructed in certain examples. These $q$-series invariants are also expected to provide an important step towards categorifying~\cite{Khovanov} the WRT invariant. The manifold $4_1(-1,2)$ illustrates the relation between the $q$-series invariants $\widehat{Z}$ and the WRT invariant, which behave in new ways when a closed hyperbolic manifold is considered. Regardless of any conjectured invariance, these $q$-series allow for the conjectural computation of Stokes data and Borel resummations~\cite{BeemDimPas,GGM:II,GGM:I,GGMW,GZ:qser,GMP} of asymptotic series~\cite{DimGar:QC,DGLZ,Gang,GSW,Hikami:stateint} associated to knots and closed three-manifolds. These collections of Stokes data have been related to the 3d index of Dimofte--Gaiotto--Gukov~\cite{DGG,GGM:I} for knots. For closed manifolds, an extension of this 3d index has been proposed~\cite{Gang3d} and an analogous conjecture seems to hold. The main results of this paper can be described as follows.
\begin{Theorem}
 The WRT invariant and the $\widehat{Z}$ invariant of the manifold $4_1(-1,2)$ are part of a~matrix valued quantum modular form.
\end{Theorem}

\subsection{Asymptotics of quantum invariants}

Quantum invariants of three-manifolds were discovered around the end of the 1980s in a union of algebra with topology. This all stemmed from Witten's interpretation~\cite{Witten} of Jones's polynomial invariants~\cite{Jones,JonesII}. Witten's physical ideas were powerful enough to suggest axioms that would completely define a theory of three-manifold invariants. This was all mathematically described by Reshetikhin--Turaev~\cite{RT:rib,RT:3man} leading to a function $\WRT(M,q)$ from fourth roots of unity to $\BC$ for closed three-manifolds $M$. While the set of axioms was enough to define a theory, Witten's original ideas have stronger links to the geometry of the three-manifolds. In particular, these new invariants should be related to connections on the three-manifold via a semi-classical limit~\cite{Andersen:WAEC,Witten}. Letting $\e(x)=\exp(2\pi {\rm i}x)$ this conjecture is formulated precisely as follows.
\begin{Conjecture}[Witten's asymptotic expansion conjecture~\cite{Witten}]
If $M$ is closed $3$-manifold, then there exists a finite set $S$ of flat connections $\rho\in\mathrm{Hom}(\pi_{1}(M),\SU(2))/\SU(2)$, $d_{\rho}\in\frac{1}{2}\BZ$ and $a^{(\rho)}_k\in\BC$ such that for $1/\hbar\in\BZ$ and $\hbar\rightarrow0$
\[
 \WRT(M,\e(\hbar))
 \sim
 \sum_{\rho\in S}\exp\left(-\frac{\CS(\rho)}{2\pi {\rm i}\hbar}\right)\hbar^{d_\rho}\bigl(a_0^{(\rho)}+a_{1/2}^{(\rho)}\hbar^{1/2}+a_{1}^{(\rho)}\hbar^{1}+\cdots\bigr) ,
\]
where $\CS$ is the Chern--Simons invariant of the connection.
\end{Conjecture}
While originally Witten only defined his invariants when $1/\hbar\in\BZ$, Reshetikhin--Turaev had no such restriction. This has led to the attempt to analytically continue Witten's original theories~\cite{Gukov,Witten:ancon}, however, no such attempt has led to enough axioms to pin down a theory at the physical level. Such a theory is expected -- at least in part -- to be described by state integrals originally introduced by Hikami~\cite{Dimofte,DGLZ,Hikami:stateint} and described precisely in work of Andersen--Kashaev~\cite{AK:II,AK:I,AK:III}. Regardless, from the mathematical point of view, we can have $\hbar\in\BQ$. Using the full theory of Reshetikhin--Turaev, a conjecture of a different form has arisen for the asymptotics of the WRT invariant of closed hyperbolic manifolds when we take $1/\hbar\notin\BZ$. This states that the invariants of Reshetikhin--Turaev should contain information about the geometric connection, which is an $\PSL_{2}(\BC)$-connection encoding the hyperbolic structure of the three-manifold. The~Chern--Simons invariant of the geometric connection stores the information of the volume of the hyperbolic three-manifold~\cite{Neu:Ext,NeuZag,Thurston:8geo,Yoshida:eta}. Therefore, this conjecture is similar to Kashaev's volume conjecture~\cite{Kashaev:VC,Murakami2} for hyperbolic knots. This conjecture was originally introduced by Chen--Yang~\cite{ChenYang} and is given as follows.
\begin{Conjecture}[Chen--Yang's volume conjecture~\cite{ChenYang}]
Let $M$ be a closed hyperbolic $3$-manifold. Then there is some $d_M\in\frac{1}{2}\BZ$ and $a_M\in\BC$ such that for $N\in\BZ$ and $N\rightarrow\infty$,
\[
 \WRT\biggl(M;\e\biggl(-\frac{1}{N+1/2}\biggr)\biggr)
 \sim
 \exp\biggl(\frac{\VC(M)}{2\pi {\rm i}}\biggl(N+\frac{1}{2}\biggr)\biggr)\biggl(N+\frac{1}{2}\biggr)^{d_M}\bigl(a_M+O\bigl(N^{-1}\bigr)\bigr) ,
\]
where $\VC(M)$ is the complexified volume of $M$, i.e., the $\SL_2(\BC)$ Chern--Simons invariant of the geometric connection.
\end{Conjecture}
These two conjectures describe very different behaviour while both describing the asymptotics of the WRT invariant as we send its argument (a root of unity) to $1$. The first conjecture states that the WRT invariant of a closed hyperbolic manifold has polynomial growth, while the second states that it has exponential growth determined by the volume. These conjectures will be seen to be unified through quantum modularity.

Before discussing quantum modularity, we should compare these asymptotic statements to the asymptotics of the $q$-series invariants $\widehat{Z}(q)$. Based on ideas in physics and explicit examples, $q$-series invariants of three-manifolds have been conjectured to exist~\cite{BeemDimPas,Cheng:3dmod,GZ:qser,GukovMan,GPPV,GPV,Hikami,LawrenceZag}. An approach using $q$-difference equations and surgery was given by Gukov--Manolescu~\cite{GukovMan} and has led to many computational examples. These $q$-series are expected to be related to the WRT invariant of the three-manifold. It was expected that the radial limits as $q$ approaches roots of unity of certain combinations of these $q$-series should be given by the WRT invariant. This works for many non-hyperbolic manifolds~\cite{AM:res,Cheng:3dmod}, however, does not seem to hold for hyperbolic manifolds. Indeed, $4_1(-1,2)$ is an integer homology sphere and there was only expected to be one series $\widehat{Z}$ so that the limit of this series as $q$ tends to a root of unity should have reproduced the WRT invariant. The series $\widehat{Z}(q)$ for this manifold was first computed by Gukov--Manolescu~\cite{GukovMan} and has expansion
\[
 \widehat{Z}(q)
  =
 -q^{-1/2}\bigl(1 - q + 2q^3 - 2q^6 + q^9 + 3q^{10} + q^{11} - q^{14} - 3q^{15}+\cdots\bigr) .
\]
For this example, the invariant does not have radial limits and instead it has exponential growth like the Chern--Simons invariant of an $\SL_{2}(\BR)$-connection. As $q$ tends to $1$ this is given explicitly in Theorem~\ref{thm:zhatasym}. Therefore, we see that while the asymptotics of $\widehat{Z}(q)$ sees certain connections, it does not always see the trivial connection, which should correspond to the WRT invariant at leading order. Indeed, if we take the asymptotics as $q$ tend to a root of unity on a small angle as opposed to radially, then the asymptotics of $\widehat{Z}$ should in fact see the geometric connections for any hyperbolic three-manifold. We can make the following conjecture similar to Chen--Yang's as the argument tends to $1$.
\begin{Conjecture}[$\widehat{Z}$ volume conjecture]
Let $M$ be a closed hyperbolic integer homology sphere. Then there is some $\epsilon\in\BR_{>0}$, $d_M\in\frac{1}{2}\BZ$ and $a_M\in\BC$ such that for $\tau\in\UHP=\{z\in\BC\mid \Im(z)>0\}$ with $0<\arg(\tau)<\epsilon$ and as $\tau\rightarrow \infty$ $($with unbounded imaginary part$)$
\[
 \widehat{Z}(M;\e(-1/\tau))
 \sim
 \exp\biggl(\frac{\VC(M)}{2\pi {\rm i}}\tau\biggr)\tau^{d_M}\bigl(a_M+O\bigl(\tau^{-1}\bigr)\bigr) .
\]
\end{Conjecture}
To try and make a precise conjecture on the asymptotics as $\tau\rightarrow {\rm i}\infty$ for $\tau\in {\rm i}\BR$ seems interesting but heavily dependent on the shape of the $\widehat{Z}$ series along with the Chern--Simons values. Without a canonical connection one could expect would dominate the asymptotics -- such as the previous expectation of the trivial connection -- a formulation is not currently known.

\subsection{Modularity conjectures}

The quantum modularity of the invariants of a hyperbolic knot was first noticed by Zagier~\cite{Zagier:Qmod} for the figure-eight knot $4_1$. The observation -- which is now a theorem~\cite{BetDrap,GZ:RQMOD} -- was that the Kashaev invariant of the figure-eight knot $J_{4_1}(q)=\sum_{k=0}^{\infty}(q;q)_{k}\bigl(q^{-1};q^{-1}\bigr)_{k}$, where $(x;q)_k=\smash{\prod_{j=0}^{k-1}\bigl(1-q^jx\bigr)}$ satisfied the following asymptotic statement for $k\in\BQ$ with bounded denominator and $k\rightarrow\infty$
\[
 J_{4_1}(\e(-1/k))
 \ae
 k^{3/2}\exp\biggl(\frac{\VC_{4_1}}{2\pi {\rm i}}k\biggr)\Phi\biggl(\frac{2\pi {\rm i}}{k}\biggr)
 J_{4_1}(\e(k)) ,
\]
for some $\Phi(\hbar)\in(-3)^{-1/4}\BQ\bigl(\sqrt{-3}\bigr)\llbracket\hbar\rrbracket$. Of course, if $k\in\BZ$, then $J_{4_1}(\e(k))=1$ and we find the all order volume conjecture~\cite{GarLe:asymp,Gukov}, however, if $k\notin\BZ$, then to obtain the same asymptotics we must divide by the Kashaev invariant at an argument obtained by a M\"obius transformation. This led to the quantum modularity conjecture~\cite{Zagier:Qmod} generalising this statement to all hyperbolic knots. While modularity had been observed for non-hyperbolic manifolds previously~\cite{Cheng:3dmod,Hikami,LawrenceZag}, an~equivalent conjecture has not been tested for a closed hyperbolic manifold. In Theorem~\ref{thm:qmodold}, we prove that the WRT invariant of $4_1(-1,2)$ satisfies a quantum modularity conjecture, which can be generalised to all closed hyperbolic manifolds as follows.
\begin{Conjecture}[quantum modularity conjecture for the WRT invariant]
Let $M$ be a closed hyperbolic three-manifold. Then there exist $d\in\frac{1}{2}\BZ$ and $\Phi(\hbar)\in\BC\llbracket\hbar\rrbracket$ such that for $k\in\BQ$ with bounded denominator and $k\rightarrow\infty$
\begin{gather}
 \WRT(M;\e(-1/k))
\nonumber\\
\qquad{} \ae
 k^{d}\exp\biggl(\frac{\VC(M)}{2\pi {\rm i}}k\biggr)
 \Phi(2\pi {\rm i}/k)(1-\e(k))\WRT(M;\e(k))
 +o\biggl(\exp\biggl(\frac{\VC(M)}{2\pi {\rm i}}k\biggr)\biggr).\label{eq:qmod.conj}
\end{gather}
\end{Conjecture}
The factor $(1-\e(k))$ appeared originally in work of Lawrence--Zagier~\cite{LawrenceZag} and appears again for the manifold $4_1(-1,2)$. This factor is responsible for the union of Chen--Yang's volume conjecture with Witten's. Indeed, when $k\in\BZ$ we see that the right-hand side of equation~\eqref{eq:qmod.conj} becomes simply \smash{$o\bigl(\!\exp\bigl(\frac{\VC(M)}{2\pi {\rm i}}k\bigr)\bigr)$}. This states that the WRT invariant grows smaller than the exponential of the volume, which is consistent with Witten's conjecture stating it should grow polynomially. To completely link these conjectures, we need to use the refined modularity Conjecture~\ref{conj:theconj}. Before discussing this refinement, we can make an analogous conjecture for the $\widehat{Z}$ series.

In the years following the modularity conjecture for the Kashaev invariant, many $q$-series associated to knots were discovered and defined~\cite{BeemDimPas,DGG,GZ:qser,GukovMan}. It was noticed by Garoufalidis--Zagier~\cite{GZ:qser} that these $q$-series have the same asymptotics as the Kashaev invariant as $q$ approaches a root of unity on a small angle. Moreover, these series start to see exponentially small corrections as we decrease the angle at which $\tau$ tends to ${\rm i}\infty$. Remarkably, these exponentially small corrections again arise from modularity. If one bounds the imaginary part of $\tau$, these corrections become $O(1)$ and therefore must be taken into account if we make an asymptotic statement. We can formulate this into a precise conjecture for $\widehat{Z}$. For the manifold $4_{1}(-1,2)$, this statement is given in Theorem~\ref{thm:qmodold.zhat}.
\begin{Conjecture}[quantum modularity conjecture for $\widehat{Z}$]
Let $M$ be a closed hyperbolic three-manifold. Then there exist $d\in\frac{1}{2}\BZ$ and $\Phi(\hbar)\in\BC\llbracket\hbar\rrbracket$ such that for $\tau\in\UHP$ with imaginary part bounded away from $0$ and $\infty$ and $\tau\rightarrow\infty$
\[
 \frac{\widehat{Z}(M;\e(-1/\tau))}{\widehat{Z}(M;\e(\tau))}
 \ae
 \tau^{d}\exp\biggl(\frac{\VC(M)}{2\pi {\rm i}}\tau\biggr)\,
 \Phi(2\pi {\rm i}/\tau) .
\]
\end{Conjecture}
The modularity conjectures along with Witten's conjecture have hinted that underneath the leading order of these asymptotics there are well behaved exponentially small corrections. Of~course, to make sense of these corrections, we can not have asymptotic series and we need some kind of analytic functions. To deal with this, there are two approaches. One uses Borel resummation and one uses state integrals.

\subsection{Stokes constants and a unification}

In work of Garoufalidis--Gu--Mari\~no~\cite{GGM:I}, it is conjectured that state integrals and the Borel resummations of the asymptotic series are one in the same. Along with conjectures on the behaviour of Stokes constants, this leads to powerful computational techniques. At the heart of the computations is a refined modularity conjecture paired with the resurgent structure. Resurgence for quantum invariants of three-manifolds was first conjectured by Garoufalidis in~\cite{Gar:resCS}. The form of these conjectures we will present here will be closer to the work of~\cite{GGM:I}.

In the language of~\cite{GGM:I}, given an asymptotic series $\Phi(\hbar)=\sum_{k=0}^{\infty}a_{k}\hbar^{k}\in\BC\llbracket\hbar\rrbracket$ whose coefficients grow factorially, we take the Borel transform $(\Bor\Phi)(\xi)=\sum_{k=0}^{\infty}a_{k}\xi^{k}/k!$. If this function has endless analytic continuation, certain growth conditions and if $\arg(\hbar)$ does not agree with the argument of a singularity, then we can take the Laplace transform $(\Lap f)(\hbar)=\int_{0}^{\infty}\exp(-\xi)f(\hbar\xi){\rm d}\xi$. Combining the Borel transform followed by the Laplace transform gives rise to an analytic function $s(\Phi)=\Lap\Bor\Phi$ called the Borel resummation that has the same asymptotics as the original function.

As part of the structure of resurgence, we obtain a set of asymptotic series $\Phi^{(\rho_{i})}(\hbar)$. The~Borel transforms of these series have singularities with potential branch cuts. The jumping across these branch cuts is determined completely by the collection of series we start with and a~collection of constants called Stokes constants. Garoufalidis--Gu--Mari\~no~\cite{Gar:resCS,GGM:I} conjecture that the singularities in Borel plane are given by the difference between two Chern--Simons values divided by $2\pi {\rm i}$ and that the series associated to two singularities differing by an integer times $2\pi {\rm i}$ are the same up to the exponential singularity, which shifts by the same multiple of~$2\pi {\rm i}$. Similar statements were considered for some non-hyperbolic manifolds in~\cite{AM:res,CG:KZser}. Therefore, the singular rays in Borel plane form a peacock pattern. The towers of singularities then all have the same asymptotic series but each contains the data of an integral Stokes constant. This is depicted in Figure~\ref{fig:bplane.stokes} for the asymptotic series associated to the trivial connection $\Phi^{(\rho_{0})}$ of the manifold~$4_1(-1,2)$.

\begin{figure}[t]\centering
\begin{tikzpicture}[scale=2.2]

\draw[<->] (-1.5,0) -- (1.5,0) node[right] {$\Re(z)$};
\draw[<->] (0,-4) -- (0,4) node[above] {$\Im(z)$};


\foreach \y in {
2-0.51413,3-0.51413,4-0.51413}
{
 \filldraw[blue] (0,\y) circle (0.5pt);
}


\draw[blue] (0,2-0.51413) node[right] {1};
\draw[blue] (0,3-0.51413) node[right] {3};
\draw[blue] (0,4-0.51413) node[right] {6};


\foreach \y in {0.17188,1+0.17188,2+0.17188,3+0.17188}
{
 \filldraw[green] (0,\y) circle (0.5pt);
}


\draw[green] (0,0.17188) node[left] {-1};
\draw[green] (0,1+0.17188) node[left] {-1};
\draw[green] (0,2+0.17188) node[left] {-3};
\draw[green] (0,3+0.17188) node[left] {-5};


\foreach \y in {0.0029434,1+0.0029434,2+0.0029434,3+0.0029434}
{
 \filldraw[purple] (0,\y) circle (0.5pt);
}


\draw[purple] (0,0.0029434) node[right] {-1};
\draw[purple] (0,1+0.0029434) node[right] {1};
\draw[purple] (0,2+0.0029434) node[right] {1};
\draw[purple] (0,3+0.0029434) node[right] {2};


\foreach \y in {1-0.23516,2-0.23516,3-0.23516,4-0.23516}
{
 \filldraw[orange] (0,\y) circle (0.5pt);
}


\draw[orange] (0,1-0.23516) node[right] {-1};
\draw[orange] (0,2-0.23516) node[right] {-2};
\draw[orange] (0,3-0.23516) node[right] {-3};
\draw[orange] (0,4-0.23516) node[right] {-8};


\foreach \y in {1-0.053934,2-0.053934,3-0.053934,4-0.053934}
{
 \filldraw[yellow] (0,\y) circle (0.5pt);
}


\draw[yellow] (0,1-0.053934) node[left] {1};
\draw[yellow] (0,2-0.053934) node[left] {1};
\draw[yellow] (0,3-0.053934) node[left] {3};
\draw[yellow] (0,4-0.053934) node[left] {9};


\foreach \y in {2.87670,3.87670}
{
 \filldraw[pink] (-25*0.035425,\y) circle (0.5pt);
 \draw[pink,->] (-25*0.035425,\y) -- (-100*0.035425/\y,4);
}

\foreach \y in {0.87670,1.87670}
{
 \filldraw[pink] (-25*0.035425,\y) circle (0.5pt);
 \draw[pink,->] (-25*0.035425,\y) -- (-1.5,\y*1.6937);
}

\foreach \y in {-3.12330}
{
 \filldraw[pink] (-25*0.035425,\y) circle (0.5pt);
 \draw[pink,->] (-25*0.035425,\y) -- (100*0.035425/\y,-4);
}

\foreach \y in {-0.12330,-1.12330,-2.12330}
{
 \filldraw[pink] (-25*0.035425,\y) circle (0.5pt);
 \draw[pink,->] (-25*0.035425,\y) -- (-1.5,\y*1.6937);
}

\draw[pink] (-25*0.035425,1-0.12330) node[right] {-1};
\draw[pink] (-25*0.035425,2-0.12330) node[right] {-4};
\draw[pink] (-25*0.035425,3-0.12330) node[right] {-12};
\draw[pink] (-25*0.035425,4-0.12330) node[right] {-38};

\draw[pink] (-25*0.035425,-0.12330) node[right] {1};
\draw[pink] (-25*0.035425,-1-0.12330) node[right] {-3};
\draw[pink] (-25*0.035425,-2-0.12330) node[right] {-1};
\draw[pink] (-25*0.035425,-3-0.12330) node[right] {2};

\foreach \y in {2.87670,3.87670}
{
 \filldraw[magenta] (25*0.035425,\y) circle (0.5pt);
 \draw[magenta,->] (25*0.035425,\y) -- (100*0.035425/\y,4);
}

\foreach \y in {0.87670,1.87670}
{
 \filldraw[magenta] (25*0.035425,\y) circle (0.5pt);
 \draw[magenta,->] (25*0.035425,\y) -- (1.5,\y*1.6937);
}


\foreach \y in {-0.12330,-1.12330}
{
 \filldraw[magenta] (25*0.035425,\y) circle (0.5pt);
 \draw[magenta,->] (25*0.035425,\y) -- (1.5,\y*1.6937);
}

\draw[magenta] (25*0.035425,1+0.02-0.12330) node[left] {1};
\draw[magenta] (25*0.035425,2+0.02-0.12330) node[left] {2};
\draw[magenta] (25*0.035425,3+0.02-0.12330) node[left] {1};
\draw[magenta] (25*0.035425,4+0.02-0.12330) node[left] {1};

\draw[magenta] (25*0.035425,0.02-0.12330) node[left] {1};
\draw[magenta] (25*0.035425,0.02-1-0.12330) node[left] {-1};

\end{tikzpicture}

\caption{Picture of the Borel plane labelled by the Stokes constants for the asymptotic series $\Phi^{(\rho_0)}$ associated to the trivial connection. The points are located at \color{blue}{$\frac {\VC_{\rho_1}}{2\pi {\rm i}}+(2\pi {\rm i})\BZ$}\color{black}{,} \color{green}{$\frac {\VC_{\rho_2}}{2\pi {\rm i}}+(2\pi {\rm i})\BZ$}\color{black}{,} \color{purple}{$\frac {\VC_{\rho_3}}{2\pi {\rm i}}+(2\pi {\rm i})\BZ$}\color{black}{,} \color{orange}{$\frac {\VC_{\rho_4}}{2\pi {\rm i}}+(2\pi {\rm i})\BZ$}\color{black}{,} \color{yellow}{$\frac {\VC_{\rho_5}}{2\pi {\rm i}}+(2\pi {\rm i})\BZ$}\color{black}{,} \color{pink}{$\frac {\VC_{\rho_6}}{2\pi {\rm i}}+(2\pi {\rm i})\BZ$}\color{black}{,} \color{magenta}{$\frac {\VC_{\rho_7}}{2\pi {\rm i}}+(2\pi {\rm i})\BZ$}\color{black}{.}\color{black}{}
See Section~\ref{stomat} for more details.}\label{fig:bplane.stokes}\vspace{-2mm}
\end{figure}

Therefore, using the conjecture that the Borel resummation is a combination of state integrals along with the factorisation of state integrals~\cite{GK:rat,GK:qser}, we are led to the following conjecture.\looseness=1
\begin{Conjecture}\label{conj:theconj}
Let $M$ by a closed $3$-manifold. Then there exists finite set of size $r$ of flat $\SL_{2}(\BC)$-connections $\rho$ with Chern--Simons values $\VC_{\rho}$ such that
\begin{itemize}\itemsep=0pt
 \item the WRT invariant $(1-q)\WRT(M;q)$ and the $\widehat{Z}(M;q)$ invariants are part a vector of functions $\mathbf{Z}\colon \UHP\cup\BQ\cup\overline{\UHP}\rightarrow\BC^{r}$ indexed by $\rho$.
 \item there exists a matrix of Borel resummable asymptotic series $\widehat{\Phi}^{(\rho)}(M;\hbar)$ indexed by $\rho$ with singularities in Borel plane given by $\VC_{\rho'}-\VC_{\rho'}$ and exponential singularities $\exp(\VC_{\rho}/\hbar)$, which resum to a matrix of functions $s\bigl(\widehat{\Phi}\bigr)(M;\hbar)$ analytic on some cut plane,
 \item for $\tau\in\UHP\cup\BQ\cup\overline{\UHP}$ we have
 \[
 \mathbf{Z}(-1/\tau)
  =
 s\bigl(\widehat{\Phi}\bigr)(M;2\pi {\rm i}/\tau)\,P_{R}(\e(\tau))\,\mathbf{Z}(\tau) ,
 \]
 where $P_{R}(q)$ is an invertible matrix of Laurent polynomials in $q$ and $R$ denotes the cone in which $\arg(\tau)$ is between the argument of two adjacent singularities.
\end{itemize}
\end{Conjecture}
This conjecture contains three parts: extension, resurgence, and modularity. The extension refers to the existence of an extension of the quantum invariants to a vector of functions from the complex numbers excluding the irrational reals~\cite{GZ:qser,GZ:RQMOD}. The resurgence refers to the conjectures of~\cite{Gar:resCS,GGM:I} on the resuregence structure of the quantum invariants. The modularity then describes how this theory is related when we have a small and large quantum parameter analogous to~$\hbar$~\mbox{\cite{GZ:qser,GZ:RQMOD,Zagier:Qmod}}. These three aspects could all be phrased independently however their interplay leads to the best depth of understanding.

Conjecture~\ref{conj:theconj} recovers all of the previous conjectures stated throughout the introduction with some additional properties of the function $\mathbf{Z}$. The leading order asymptotics of this conjecture when $\tau\in\BZ$ recover Witten's conjecture assuming that $\mathbf{Z}^{(\rho)}(\BZ)$ vanish for all connections $\rho$ with $\Im(V_{\rho})\geq0$ that are not valued in $\SU(2)$. Similarly, Chen--Yang's conjecture follows from the conjecture assuming that for the geometric connection $\mathbf{Z}^{(\rho_{\mathrm{geo}})}\bigl(\tfrac{1}{2}+\BZ\bigr)\neq0$. The quantum modularity conjectures follow noting that the geometric connection gives the largest possible contribution when $\tau\in\BR$.

Even this conjecture can be extended. Indeed, in the factorisation of state integrals~\cite{GK:rat,GK:qser} these kind of statements can be lifted to matrix equations as was emphasised in~\cite{GZ:qser,GZ:RQMOD}. This lift to a matrix equation was utilised in~\cite{GGM:II,GGM:I,GGMW} to calculate Stokes constants associated to the asymptotic series of Conjecture~\ref{conj:theconj}. In the case of $4_{1}(-1,2)$, we utilise the same methods to conjecturally compute the generating series of Stokes constants. If $s\bigl(\widehat{\Phi}^{(\rho_{i})}\bigr)(2\pi {\rm i}/\tau)$ is the Borel resummation of \smash{$\widehat{\Phi}^{(\rho_{i})}$} then as the argument crosses \smash{$\frac{\VC_{\rho_{j}}-\VC_{\rho_{i}}}{2\pi {\rm i}}+2\pi {\rm i}k$} it is conjectured that~the~series jumps by the addition of an integer constant times $q^{k}s\bigl(\Phi^{(\rho_{j})}\bigr)(2\pi {\rm i}/\tau)$. These constants can be~stored\footnote{To explicitly extract what are usually called Stokes constants, one must factorise this matrix into unipotent matrices associated to each argument in the upper or lower half plane and take their logarithm. See Section~\ref{stomat}.} in a matrix of $q$-series with entries given by
\begin{gather*}
 \mathsf{S}_{+}^{(\rho_i,\rho_j)}(q)
 =
 \sum_{k\colon \Re\bigl(\frac{\VC_{\rho_{j}}-\VC_{\rho_{i}}}{4\pi^2}\bigr)+k>0}\mathsf{S}^{(\rho_i,\rho_j)}_{k}q^{k},
 \\
 \mathsf{S}_{-}^{(\rho_i,\rho_j)}(q)
 =
 \sum_{k\colon \Re\bigl(\frac{\VC_{\rho_{j}}-\VC_{\rho_{i}}}{4\pi^2}\bigr)+k<0}\mathsf{S}^{(\rho_i,\rho_j)}_{k}q^{-k}.
\end{gather*}
These generating series do not collect Stokes constants on the real line but as a pair see all others.
Then using the techniques introduced in~\cite{GGM:I}, we can find combinations of entries of~$\mathbf{Z}(q)$ and~$\mathbf{Z}(q)^{-1}$ conjecturally giving the matrices of Stokes constants. The second is closely related to~$\mathbf{Z}\bigl(q^{-1}\bigr)^{t}$. Therefore, we can express the matrices $\mathsf{S}_{\pm}(q)$ as explicit combinations of the $q$-series of equation~\eqref{eq:theqser}.
In~\cite[Conjecture 2]{GGM:I},
an explicit conjecture for the generating series of the Stokes constants of the geometric connection with itself was given in terms of the 3d index of Dimofte--Gaiotto--Gukov~\cite{DGG}. A surgery formula for the 3d index was conjectured by Gang--Yonekura in~\cite{Gang3d} that gave a conjectural extension of the 3d index to closed manifolds. These conjectures seems to hold for this closed manifold and indeed one finds
\begin{gather*}
 1+\mathsf{S}_{+}^{(\rho_6,\rho_6)}(q)+\mathrm{O}\bigl(q^{10}\bigr)
  =
 \mathcal{I}_{4_{1}(-1,2)}(q)+\mathrm{O}\bigl(q^{10}\bigr)\\
 \qquad{} =
 1 - 2q - 3q^2 + 3q^4 + 10q^5 + 14q^6 + 22q^7 + 20q^8 + 14q^9+\mathrm{O}\bigl(q^{10}\bigr) ,
\end{gather*}
where $\rho_6$ is the geometric connection. For many more Stokes constants see Appendix~\ref{app:stokes}.

\section[The manifold $4_1(-1,2)$]{The manifold $\boldsymbol{4_1(-1,2)}$}

In this section, we will discuss some of the classical and quantum aspects of the manifold~$4_1(\!-1,2)$. In particular, we will study its $\SL_{2}(\BC)$ character variety. We will then consider its associated~$\mathfrak{sl}_2$ quantum invariant, that is, its WRT invariant. This will take the form of an element in the Habiro ring, which -- in the weakest sense -- is a function from roots of unity to $\mathbb{C}$. Finally, we will discuss a construction of $q$-series from the WRT invariant that matches $q$-series~$\widehat{Z}$ found by Gukov--Manolescu in~\cite{GukovMan}.

\subsection[The $\SL_{2}(\BC)$-character variety]{The $\boldsymbol{\SL_{2}(\BC)}$-character variety}

The manifold $4_1(\!-1,2)$ is a closed, hyperbolic, integer homology sphere of volume $1.398508884\dots$ To compute the various quantities associated with $4_1(-1,2)$, we take the standard (for example~\cite{Thurston:GT3M}) two tetrahedral triangulation of $4_1$ with shapes $w_1$, $w_2$ and add an equation to fill in the cusp. The character variety of the figure eight knot consists of the component with trivial longitude and the geometric component (compare the appendix of~\cite{CCGLS} to the example~\cite[equation~(4.20)]{DimGar:QC}). All connections on $4_1(-1,2)$ therefore come from either a point on the geometric component or the trivial connection. The geometric component is contained in the solutions to the gluing equations and therefore all connections beside the trivial connection come from solutions to the gluing equations.

This is implemented in SnapPy~\cite{SnapPy} and these equations can be computed using the code:
\begin{lstlisting}
M=Manifold('4_1')
M.dehn_fill((-1,2))
M.gluing_equations(form='rect')
/*[([2, 2] ,[-1, -1], 1), ([-2, -2] ,[1, 1], 1), ([-1, 7] ,[0, -3], -1)]*/
\end{lstlisting}
Therefore, we find that
\begin{align*}
 w_1^{-2}w_2^{-2}(1-w_1)(1-w_2) =1,\qquad
 w_1^{-1}w_2^7(1-w_2)^{-3} =-1 .
\end{align*}
Take a change of variables
\[
 z_{1} =w_{1}
 \qquad\text{and}\qquad
 z_{2} =\bigl(1-w_2^{-1}\bigr)^{-1} .
\]

The trace field $F$ of $4_1(-1,2)$ is the number field of type $[5,1]$ with discriminant $-7215127$ (a prime number),
generated by a root of $p(x)$, where
\be
\label{tracefield}
p(x) =x^7 - x^6 - 2 x^5 + 6 x^4 - 11 x^3 + 6 x^2 + 3 x - 1  .
\ee
The shape field is a quadratic extension of $F$ generated by a root of $pp(z)$, where
\begin{gather}
pp(z) =z^{14} + 2 z^{13} + z^{12} - 4 z^{10} - 8 z^9 - 10 z^8 - 13 z^7 -
10 z^6 - 8 z^5 \nonumber\\ \hphantom{pp(z) =}{}
 - 4 z^4 + z^2 + 2 z + 1 .\label{shapefield}
\end{gather}
If $p(\xi)=0$, then the quadratic extension is given by the polynomial
\[
x^2-\bigl(1+4\xi+4\xi^2-3\xi^3+\xi^4-\xi^5-\xi^6\bigr)x+1 =0  .
\]
We find that for $pp(\psi)=0$ the shapes are given by
\begin{gather}
z_1  = 2 \psi^{13} + 3 \psi^{12} + \psi^{11} - 9 \psi^9
- 12 \psi^8 - 15 \psi^7 - 20 \psi^6 - 10 \psi^5 \nonumber\\
\hphantom{z_1  =}{}
- 13 \psi^4 - \psi^3 + \psi^2 + 2 \psi + 4 ,\nonumber \\
z_2  = -2 \psi^{13} - 3 \psi^{12} - \psi^{11} + 9 \psi^9
+ 12 \psi^8 + 15 \psi^7 + 20 \psi^6 + 10 \psi^5 \nonumber\\
\hphantom{z_2  =}{}
+ 12 \psi^4 + \psi^3 - 2 \psi - 3.\label{eq.z}
\end{gather}
These points come in seven pairs and correspond to points on the character variety and therefore give flat $\SL_{2}(\BC)$ connections or representations in $\mathrm{Hom}(\pi_{1}(4_1(-1,2)),\SL_{2}(\BC))/\SL_{2}(\BC)$. Most of the points actually correspond to representations landing in a subgroup such as $\SL_{2}(\BR)$ or~$\mathrm{SU}(2)$.\looseness=1

Each point on the character variety has a complex volume associated to it~\cite{Neu:Ext}. The geometric point on the character variety gives rise to the actual hyperbolic volume. To compute these volumes let $R$ be the enhanced Rogers dilogarithm,
\[
R(z, p, q)  = \Li_2(z)
+ \frac{1}{2}(\log(z) + 2 \pi {\rm i} p) (\log(1 - z) - 2 \pi {\rm i} q)  .
\]
The set of complex volumes of $\SL(2,\BC)$ representations was given to me by S. Garoufalidis who received them from C. Zickert.
There are $7$ complex volumes and each has the form
\begin{gather*}
\VC_{\rho_j}  = - R(z_1(\alpha_j),p_{1,j},q_{1,j})
+ R(z_2(\alpha_j),p_{2,j},q_{2,j}),
\end{gather*}
where $\alpha_j$ for $j=1,\dots,14$ are the roots of the polynomial of
equation~\eqref{shapefield},
$z_1$ and $z_2$ are given by~\eqref{eq.z}, and
the pairs $r_{i,j}:=(p_{i,j},q_{i,j})$ for $i=1,2$ and $j=1,\dots,7$
are given in the following table.

\begin{center}
\def\arraystretch{1.2}
\begin{tabular}{|l|l|l|l|l|l|l|l|}\hline
$j$ & $1$ & $2$ & $3$ & $4$ & $5$ & $6$ & $7$
\\ \hline
$\Re(\alpha_j)$ & $0.477$ & $0.156$ & $0.616$ & $-0.194$ & $-0.424$ & $-0.693$ & $-0.693$
\\ \hline
$\Im(\alpha_j)$ & $0.879$ & $0.988$ & $0$ & $0.981$ & $0.906$ & $0.0194$ & $-0.0194$
\\ \hline
$r_{1,j}$ & $(0,0)$ & $(0,1)$ & $(0,0)$ & $(0,0)$ & $(0,-1)$ & $(0,0)$ & $(0,0)$ \\ \hline
$r_{2,j}$ & $(-2,-3)$ & $(0,0)$ & $(0,0)$ & $(1,1)$ & $(-1,-2)$ & $(-1,-2)$ & $(1,2)$ \\ \hline
\end{tabular}
\end{center}

Besides the complex volumes, there is an important element of the trace field associated to each point on the character variety called the $1$-loop invariant. This is related to the Reidermeister torsion and given in this example by
\[
 \delta =-12\xi^6+15\xi^5+31\xi^4-74\xi^3+133\xi^2-66\xi-74 ,
\]
where $p(\xi)=0$ with $p$ of equation~\eqref{tracefield}. The seven complex volumes \big(modulo $(2\pi {\rm i})^2$\big) and $1$-loop invariants with corresponding generator of an embedding of the trace field are given numerically~by
\begin{gather*}
\begin{aligned}
&\VC_{\rho_1} = 20.297\dots, && \delta_{\rho_1}  = -11.578\dots, \\
&\xi_{\rho_1} = -2.2411\dots,\\
&\VC_{\rho_2} = -6.7857\dots, && \delta_{\rho_2}  = -12.636\dots \\
& \xi_{\rho_2}  = -0.43760\dots,\\
&\VC_{\rho_3} = 39.362\dots, && \delta_{\rho_3}  = -83.275\dots,\\
 & \xi_{\rho_3}  = 0.25599\dots,\\
&\VC_{\rho_4} = 9.2837\dots, && \delta_{\rho_4}  = -7.0205\dots,\\
 & \xi_{\rho_4}  = 1.3348\dots,\\
&\VC_{\rho_5} = 2.1292\dots, && \delta_{\rho_5}  = -5.3937\dots,\\
 & \xi_{\rho_5}  = 1.3483\dots,\\
&\VC_{\rho_6} = 4.8678\ldots - {\rm i}1.3985\dots, \qquad && \delta_{\rho_6}  = 3.9517\ldots - {\rm i}0.15252\dots,\\
 & \xi_{\rho_6}  = 0.36981\ldots-{\rm i}1.4410\dots,\\
&\VC_{\rho_7} = 4.8678\ldots + {\rm i}1.3985\dots, \qquad && \delta_{\rho_7}  = 3.9517\ldots + {\rm i}0.15252\dots,\\
& \xi_{\rho_7}  = 0.36981\ldots+{\rm i}1.4410\dots .
\end{aligned}
\end{gather*}
In terms of representations, $\rho_{3}$ corresponds to an $\SL_{2}(\BR)$ representation as the corresponding roots of $p$ and $pp$ are real, while $\rho_{6}$ and $\rho_{7}$ are conjugate $\SL_{2}(\BC)$ representations as both the corresponding roots of $p$ and $pp$ are imaginary giving the geometric and anti-geometric connections respectively, which have $\Re({\rm i}\VC_{\rho_{6}})=-\Re({\rm i}\VC_{\rho_{7}})=\mathrm{Vol}(4_{1}(-1,2))$. The rest are all $\mathrm{SU}(2)$ representations as the corresponding roots of $p$ are real but those of $pp$ are imaginary. We will often denote $\xi$ as the generator of the field defined by $p(\xi)=0$ and use $\rho_i$ when we specify an embedding into~$\BC$.

\subsection{The WRT invariant}

The WRT invariant~\cite{RT:3man} is given as a sum over the coloured Jones polynomial~\cite{Jones,JonesII,RT:rib} of a~framed link representing the three-manifold. In the case of $4_{1}(-1,2)$, we have a nice formula for the coloured Jones polynomial of $4_1$ and therefore it is best to use a surgery formula even if it is not integral surgery. This formula was given by Beliakova--Blanchet--L\^e in~\cite{BelBlaLe} and they show that the WRT invariant -- normalised by a factor of $(1-q)$ observed in~\cite{LawrenceZag} -- of $4_1(-1,2)$ is given as an element of Habiro ring by
\be\label{WRTformula}
 \WRT(q)  =
 \sum_{0\leq\ell\leq k}\alpha_{k,\ell}(q)
  =
 \sum_{0\leq\ell\leq k}(-1)^{k}q^{-\frac{1}{2}k(k+1)+\ell(\ell+1)}\frac{(q;q)_{2k+1}}{(q;q)_{\ell}(q;q)_{k-\ell}} .
\ee
This formula will be our starting point. As this is an element of the Habiro ring, we can formally evaluate at $q={\rm e}^\hbar$ so that
\[
 \WRT\bigl(e^\hbar\bigr) =-\hbar -\frac{25}{2}\hbar^2 - \frac{1621}{6}\hbar^3 - \frac{195601}{24}\hbar^4 - \frac{37907101}{120}\hbar^5+\cdots .
\]
Using this normalisation of the WRT invariant, which matches some of the conventions used in~\cite{LawrenceZag}, we note that
\[
 \WRT(1) =0 .
\]
This is of absolutely fundamental importance for relating various conjectures~\cite{ChenYang,Witten} on the asymptotics of WRT invariants as described in Section~\ref{volconj} and we will come back to it later.
The~$q$-hypergeometric expression~\eqref{WRTformula} has a natural $q$-holonomic module associated to it generated by
\[
 \WRT_{m,n}(q)
  =
 \sum_{0\leq\ell\leq k}\alpha_{k,\ell}(q)q^{mk+n\ell}
  =
 \sum_{0\leq\ell\leq k}(-1)^{k}q^{-\frac{1}{2}k(k+1)+\ell(\ell+1)+mk+n\ell}\frac{(q;q)_{2k+1}}{(q;q)_{\ell}(q;q)_{k-\ell}} .
\]
Noting that
\begin{gather}
 \bigl(1-q^{k+1-\ell}\bigr)\alpha_{k+1,\ell}(q)
  =-q^{-k-1}\bigl(1-q^{2k+2}\bigr)\bigl(1-q^{2k+3}\bigr)\alpha_{k,\ell}(q),\nonumber\\
 \bigl(1-q^{\ell+1}\bigr)\alpha_{k,\ell+1}(q)
  =q^{2\ell+2}\bigl(1-q^{k-\ell}\bigr)\alpha_{k,\ell}(q) ,\label{alpharec}
\end{gather}
we can sum both sides to determine relations between $\WRT_{m,n}(q)$ for various $m,n$.
This $q$-holonomic module is rank eight. The dependence on $m$ is of most interest and using standard $q$-holonomic techniques it is not hard to show that
\begin{gather}
 q^{2m+2}\WRT_{m,n}(q)+\bigl(q^{m+1}+q^{m+2}\bigr)\WRT_{m+1,n}(q) \nonumber\\
 \qquad{}+\bigl(-q^{2m+4}-q^{2m+5}-q^{2m+6}-q^{2m+7}-q^{m+2}+1\bigr)\WRT_{m+2,n}(q) \nonumber\\
 \qquad{}+\bigl(q^{n+m+3}-q^{m+3}-2q^{m+4}-q^{m+5}-1\bigr)\WRT_{m+3,n}(q) \nonumber\\
 \qquad{}+\bigl(q^{2m+7}+q^{2m+8}+2q^{2m+9}+q^{2m+10}+q^{2m+11}+q^{m+4}+q^{m+5}\bigr)\WRT_{m+4,n}(q) \nonumber\\
 \qquad{}+\bigl(-q^{n+m+5}-q^{n+m+6}q^{m+6}+q^{m+7}\bigr)\WRT_{m+5,n}(q)\nonumber\\
 \qquad{}+\bigl(-q^{2m+11}-q^{2m+12}-q^{2m+13}-q^{2m+14}-q^{m+7}\bigr)\WRT_{m+6,n}(q) \nonumber\\
 \qquad{}+q^{n+m+8}\WRT_{m+7,n}(q)+q^{2m+16}\WRT_{m+8,n}(q) =0 .\label{8thorderqdiff}
\end{gather}
Taking the classical limit \big($q\to 1$, $q^m\to 1$, $\WRT_{m+j,n}(q)\to z^j$\big) of the above equation gives a~polynomial equation
\begin{gather*}
 z^8+z^7-5z^6+8z^4-4z^3-4z^2+2z+1\nonumber\\
 \qquad{}=(z-1)\bigl(z^7+2z^6-3z^5-3z^4+5z^3+z^2-3z-1\bigr)=0 .\label{eq:z.pol}
\end{gather*}
The field defined by the degree $7$ polynomial is of course the trace field of $4_1(-1,2)$ indeed $7\xi^6 + 2\xi^5 - 11\xi^4 + 28\xi^3 - 42\xi^2 - 10\xi + 6$ with $p(\xi)=0$ is a solution. This factorisation could suggest a refinement of this equation. One can in fact show that
\begin{gather*}
 q^{2m+2}\WRT_{m,n}(q)+\bigl(q^{2m+4}+q^{m+1}+q^{m+2}\bigr)\WRT_{m+1,n}(q)\nonumber\\
 \qquad{}+\bigl(-q^{2m+7}-q^{2m+5}-q^{2m+4}+q^{m+3}+1\bigr)\WRT_{m+2,n}(q)\nonumber\\
 \qquad{}+\bigl(-q^{2m+9}-q^{2m+7}-q^{2m+6}+q^{m+n+3}-q^{m+5}-q^{m+4}-q^{m+3}\bigr)\WRT_{m+3,n}(q)\nonumber\\
 \qquad{}+\bigl(q^{2m+10}+q^{2m+9}+q^{2m+7}+q^{m+n+4}-q^{m+6}\bigr)\WRT_{m+4,n}(q)\nonumber\\
 \qquad{}+\bigl(q^{2m+12}+q^{2m+11}+q^{2m+9}-q^{m+n+6}+q^{m+6}\bigr)\WRT_{m+5,n}(q)\nonumber\\
 \qquad{}+\bigl(-q^{2m+12}-q^{m+n+7}\bigr)\WRT_{m+6,n}(q)-q^{2m+14}\WRT_{m+7,n}(q) =(1-q) .\label{7thorderqdiffinhom}
\end{gather*}
These $q$-difference equations are related to those constructed in~\cite[Theorem 1]{Wh:thesis} and come from this particular surgery presentation of the manifold. The vector space associated to closed manifolds -- which was shown to be an invariant -- is then expected to agree with the module coming from skein theory of~\cite{GarLe:module}. There it has been shown that the module for $4_1(-1,2)$ has rank at most eight and using Theorem~\ref{thm:qmodold} we can show that it has rank exactly eight as mentioned in~\mbox{\cite[Example 3.3]{GarLe:module}}.
\begin{Proposition}
The $\BQ(q)$--vector space $\mathrm{Span}_{\BQ(q)}\{\WRT_{m,n}\mid m,n\in\BZ\}$ has dimension eight.
\end{Proposition}
\begin{proof}
Firstly, using the $q$-holonomic relations, we can show that $\mathrm{Span}_{\BQ(q)}\{\WRT_{m,n}\mid m,n\in\BZ\}$ is equal to $\mathrm{Span}_{\BQ(q)}\{\WRT_{m,0}\mid m\in\BZ\}$ and that this has dimension at most eight. Then using a~simple extension of the asymptotics of Theorem~\ref{thm:qmodold} to the case $\WRT_{m,0}(q)$ one can show that $\WRT_{m,0}(q)/\WRT_{0,0}(q)\rightarrow z^m$ as $q\rightarrow 1$ where $q=\e(1/(n+1/2))$ for $n\in\BZ$ say and where $z$ makes the order seven polynomial factor of equation~\eqref{eq:z.pol} vanish. Therefore, the leading asymptotics span at least a seven dimensional space and give rise to a $q$-holonomic module, which is just the homogeneous part of equation~\eqref{7thorderqdiffinhom}. If $\WRT_{m,0}(q)$ for $m=0,\dots,7$ satisfy a relation then the leading asymptotics must. This relation must therefore be a multiple of the homogeneous part of equation~\eqref{7thorderqdiffinhom}. However, the WRT invariant satisfies the inhomogeneous version and therefore $\WRT_{m,0}(q)$ for $m=0,\dots,7$ cannot satisfy a homogeneous equation. Therefore, the dimension is at least eight.
\end{proof}

\subsection[The $\widehat{Z}$ series]{The $\boldsymbol{\widehat{Z}}$ series}

In~\cite{GukovMan}, Gukov and Manolescu give a construction of a two variable series associated to knot complements. This two variable series has variables $x$ and $q$ where $|q|<1$ is the variable replacing roots of unity and $x$ behaves like $x=q^{N}$, where $N$ represents the colour of the coloured Jones polynomial~\cite{Jones,RT:rib}. This series is computed by solving the $\hat{A}$-polynomial associated to the knot as a $q$-series and choosing initial conditions coming from the expansion of the coloured Jones polynomial~\cite{BarNatanGar,MelvinMorton,Rozansky}. Using this two variable series, Gukov and Manolescu propose a~surgery formula similar to~\cite{BelBlaLe}, which conjecturally leads to invariants of closed three-manifolds. For~$4_1(-1,2)$, they compute~\cite[Table 10]{GukovMan},
\be\label{zhatGM}
 \widehat{Z}(q)
  =
 -q^{1/2}\bigl(1 - q + 2q^3 - 2q^6 + q^9 + 3q^{10} + q^{11} - q^{14} - 3q^{15}+\cdots\bigr) .
\ee
This $q$-series will be seen to come from the specialisation of a solution of the $q$-difference equation~\eqref{8thorderqdiff} associated to the WRT invariant of $4_1(-1,2)$. Indeed, it is a natural $q$-series that pairs with $\WRT(q)$. In particular, consider $\alpha$ from equation~\eqref{WRTformula} as a function
\[
 \alpha\colon\ \BZ^2\rightarrow\BQ(q)
 \qquad \text{s.t.}\qquad
 \alpha_{k,\ell}(q)
  =\begin{cases}
 \dfrac{(-1)^{k}q^{-\frac{1}{2}k(k+1)+\ell(\ell+1)}(q;q)_{2k+1}}{(q;q)_{\ell}(q;q)_{k-\ell}} & \text{if }0\leq\ell\leq k,\\
 0 & \text{otherwise}.
 \end{cases}
\]
Equations~\eqref{alpharec} still hold for $\alpha$ with this extended domain. Moreover, equations~\eqref{alpharec} completely determine the module. Therefore, it is natural to ask whether there is another solution to equations~\eqref{alpharec}. Indeed, there is another solution given by
\begin{gather*}
 \beta\colon \ \BZ^2\rightarrow\BQ(q)\\
\text{s.t.}\qquad
 \beta_{k,\ell}(q)
  =\begin{cases}
 \dfrac{(-1)^{k+\ell}q^{\frac{1}{2}3k(k+1)+\frac{1}{2}\ell(\ell+1)+k+1}(q;q)_{-\ell-1}}{(q;q)_{-2k-2}(q;q)_{k-\ell}} & \text{if }\ell\leq k\leq-1,\\
 0 & \text{otherwise}.
 \end{cases}
\end{gather*}
For positive $0\leq k\leq\ell$, we have
\[
 \beta_{-k-1,-\ell-1}(q)
  =
 (-1)^{k+\ell}q^{\frac{1}{2}3k(k+1)+\frac{1}{2}\ell(\ell+1)-k}\frac{(q;q)_{\ell}}{(q;q)_{2k}(q;q)_{\ell-k}} .
\]
We see that the sum of these terms is convergent when $|q|<1$. Indeed, we find that
\begin{align}
 Z_{0,0}(q)
  =
 \sum_{k,\ell\in\BZ}\beta_{k,\ell}(q)
 &{}=
 \sum_{0\leq k\leq\ell}(-1)^{k+\ell}q^{\frac{1}{2}3k(k+1)+\frac{1}{2}\ell(\ell+1)-k}\frac{(q;q)_{\ell}}{(q;q)_{2k}(q;q)_{\ell-k}}\nonumber\\
 &{}=
 1 - q + 2q^3 - 2q^6 + q^9 + 3q^{10} + q^{11} - q^{14} - 3q^{15}+\cdots .\label{eq:Z00}
\end{align}
We will define more generally
\begin{align*}
 Z_{m,n}(q)
 & =
 \sum_{k,\ell\in\BZ}\beta_{k,\ell}(q)q^{mk+n\ell}\\
 & =
 \sum_{0\leq k\leq\ell}(-1)^{k+\ell}q^{\frac{1}{2}3k(k+1)+\frac{1}{2}\ell(\ell+1)-(m+1)k-m-n\ell-n}\frac{(q;q)_{\ell}}{(q;q)_{2k}(q;q)_{\ell-k}} .
\end{align*}
We see that up to a factor of $-q^{1/2}$, this is the $\widehat{Z}$-series associated to $4_1(-1,2)$~\cite{GukovMan} given in equation~\eqref{zhatGM}. We can use a formula for the two variable series of Gukov--Manolescu found by Park~\cite{GGMW,Park} to prove the equality. Indeed, following the notation of~\cite{GukovMan,Park}, we have
\begin{align*}
 \Xi_{4_1}(x;q)
 &{}=
 \bigl(x^{1/2}-x^{-1/2}\bigr)\sum_{k=0}^{\infty}(-1)^k\frac{q^{k(k+1)/2}}{(x;q)_{k+1}(x^{-1};q)_{k+1}}\\
 &{}=
 \sum_{k,j,\ell=0}^{\infty}\bigl(x^{k+j+\ell+1/2}-x^{k+j+\ell+3/2}\bigr)\binom{k+j}{j}_q\binom{k+\ell}{\ell}_{q^{-1}},
\end{align*}
and so
\begin{gather*}
 F_{4_1}(x;q)
  =\frac{1}{2}\bigl(\Xi(x,q)-\Xi\bigl(x^{-1},q\bigr)\bigr) ,\\
\hphantom{F_{4_1}(x;q)  =}{} =
 \frac{1}{2}\sum_{k,j,\ell=0}^{\infty}\bigl(x^{k+j+\ell+1/2}-x^{k+j+\ell+3/2}-x^{-k-j-\ell-1/2}+x^{-k-j-\ell-3/2}\bigr)\\
\hphantom{F_{4_1}(x;q)=\frac{1}{2}\sum_{k,j,\ell=0}^{\infty}}{}\quad \times \binom{k+j}{j}_q \binom{k+\ell}{\ell}_{q^{-1}}.
\end{gather*}
Then, using the surgery formula in~\cite{GukovMan}, where for $\mathcal{L}_{-1/2}^{(0)}(x^uq^v)=q^{2u^2+v}$
\begin{align*}
 \widehat{Z}(q)
 &{}=
 q^{3/8}\mathcal{L}_{-1/2}^{(0)}F_{4_1}(x;q)
 \\
 &{}=
 q^{3/8}\sum_{k,j,\ell=0}^{\infty}\bigl(q^{2(k+j+\ell+3/4)^2}-q^{2(k+j+\ell+7/4)^2}-q^{2(k+j+\ell+1/4)^2}+q^{2(k+j+\ell+5/4)^2}\bigr)\\
 &\hphantom{=q^{3/8}\sum_{k,j,\ell=0}^{\infty}}{}\ \times\binom{k+j}{j}_q\binom{k+\ell}{\ell}_{q^{-1}}\\
 &{}=
 -q^{1/2}\bigl(1 - q + 2q^3 - 2q^6 + q^9 + 3q^{10} + q^{11} - q^{14} - 3q^{15}+\cdots\bigr) .
\end{align*}

\begin{Proposition}\label{prop:zhat}
Let
\begin{align*}
 \widehat{Z}_{m}(q)&{}=q^{3/8}\sum_{k,j,\ell=0}^{\infty}\bigl(q^{2(k+j+\ell+3/4)^2}-q^{2(k+j+\ell+7/4)^2}-q^{2(k+j+\ell+1/4)^2}+q^{2(k+j+\ell+5/4)^2}\bigr)\\
 &\hphantom{=q^{3/8}\sum_{k,j,\ell=0}^{\infty}}{}\ \times\binom{k+j}{j}_q\binom{k+\ell}{\ell}_{q^{-1}}q^{-mk-m} .
\end{align*}
We have the following identity:
\[
 \widehat{Z}_{m}(q) =-q^{1/2}Z_{m,0}(q) ,\qquad \text{in particular,} \qquad
 \widehat{Z}(q) =-q^{1/2}Z_{0,0}(q) .
\]
\end{Proposition}
\begin{proof}
One can show using holonomic function techniques~\cite{WilfZeil} that $\widehat{Z}_{m}$ satisfies the $q$-difference equation~\eqref{8thorderqdiff} with $n=0$. Therefore, as both are power series in $q^{-m}$ and the $q$-difference equations are second order in $q^{m}$, if the coefficients of $q^{-m}$ and $q^{-2m}$ agree this proves the result. The first of these equalities is given by
\begin{gather*}
 -q^{-1/8}\!\sum_{j,\ell=0}^{\infty}\!\bigl(q^{2(j+\ell+3/4)^2}\!-q^{2(j+\ell+7/4)^2}\!-q^{2(-j-\ell-1/4)^2}\!+q^{2(-j-\ell-5/4)^2}\bigr)
 \! =\!\sum_{j=0}^{\infty}(-1)^jq^{j(j+1)/2} ,
\end{gather*}
which can be proved by direct computation and similarly for the second.
\end{proof}

\section{Asymptotics}

Recently, volume conjectures have been introduced for quantum invariants of closed hyperbolic three-manifolds~\cite{ChenYang}. Moreover, Witten's asymptotic expansion conjecture~\cite{Witten} states that the WRT invariant should be expressed as a sum of oscillating asymptotic series of polynomial growth. As in the case of knots, these conjectures have been proved for some families of hyperbolic closed $3$-manifolds~\cite{CharlesMarche,Ohtsuki:VCYVC} including our example of interest $4_1(-1,2)$ but remain open in general.~On the other hand, the $\widehat{Z}$ invariants predicted in~\cite{GukovMan,GPV} had conjectural radial limits to WRT invariants. We will study the asymptotics of both of these invariants and see that they are related by quantum modularity. This will extend the volume conjecture for WRT invariants. For the $q$-series, it will in fact give a counterexample to the conjectured radial limits and provide a~reformulation consistent with the previous examples.

\subsection{Volume conjectures and asymptotic series}
\label{volconj}

To get a feeling for the behaviour of the WRT invariant, it is best to simply compute. Consider the following values of the WRT invariant:
\begin{align*}
& \WRT\biggl(\e\biggl(\frac{1}{1000}\biggr)\biggr) =-6.3258\ldots + 14.804\dots {\rm i} ,\\
& \WRT\biggl(\e\biggl(\frac{1}{1000+1/2}\biggr)\biggr) =(1.6433 \ldots + 0.067579\dots {\rm i})\times 10^{98} ,\\
& \WRT\biggl(\e\biggl(\frac{1}{1000+1/3}\biggr)\biggr) =(1.0551\ldots + 0.87759\dots {\rm i})\times 10^{98} .
\end{align*}
The immediate observation is the order of magnitude shift. The behaviour of the first two values is expected from Witten's asymptotic expansion conjecture~\cite{Witten} and Chen--Yang's volume conjecture for WRT invariants~\cite{ChenYang}. The third value is predicted from a generalisation of the second. Indeed, notice that
\begin{align*}
& \WRT\biggl(\e\biggl(\frac{1}{1000+1/2}\biggr)\biggr)\e\biggl(-\frac{\VC_{\rho_7}}{(2\pi {\rm i})^2}(1000+1/2)\biggr) =-22.044\ldots + 22.943\dots {\rm i} ,\\
& \WRT\biggl(\e\biggl(\frac{1}{1000+1/3}\biggr)\biggr)\e\biggl(-\frac{\VC_{\rho_7}}{(2\pi {\rm i})^2}(1000+1/3)\biggr) =-26.465\ldots + 7.6613\dots {\rm i} .
\end{align*}
More generally, we have the following observation which is in fact a theorem.\footnote{See Theorem~\ref{thm:qmodold} and its proof in Section~\ref{statphase} for a proof of a stronger statement.}
\begin{Theorem}
For $r\in\BQ\smallsetminus \BZ$, we have
\begin{align*}
 \lim_{k\rightarrow\infty}\frac{\log\bigl(\WRT\bigl(\e\bigl(\frac{1}{k+r}\bigr)\bigr)\e\bigl(-\frac{\VC_{\rho_7}}{(2\pi {\rm i})^2}(k+r)\bigr)\bigr)}{k+r}
  =0 .
\end{align*}
\end{Theorem}
This theorem is a slight generalisation of the Chen--Yang volume conjecture~\cite{ChenYang}. However, as in~\cite{Ohtsuki:VCYVC}, we can go further and find a full asymptotic expansion. Indeed, numerically (using extrapolation methods described for example in~\cite{GrMoree,Wh:thesis,Zagier:Vas}) one finds that there is some \hbox{$X(r)\in\BZ[\e(-r)]$} such that\footnote{We take $x_n\sim C^nn^d\bigl(A_0+A_1n^{-1}+A_2n^{-2}+\cdots\bigr)$ when $C^{-n}n^{-d}x_n-\sum_{k=0}^{K}A_kn^{-k}=\mathrm{o}\bigl(n^{K}\bigr)$ for all $K$.}
\be\label{qmodnum}
 \WRT\biggl(\e\biggl(\frac{1}{k+r}\biggr)\biggr)
 \sim
 X(r)\,\e(3/8)\,\sqrt{k+r}\,\e\biggl(\frac{\VC_{\rho_7}}{(2\pi {\rm i})^2}(k+r)\biggr)\,\Phi^{(\rho_7)}\biggl(\frac{2\pi {\rm i}}{k+r}\biggr),
\ee
where
\[
 \Phi^{(\rho_7)}(\hbar)
  =
 a^{\rho_7}_0+a^{\rho_7}_1\hbar+a^{\rho_7}_2\hbar^2+\dots
  =
 \frac{1}{\sqrt{\delta_{\rho_7}}}\bigl(1+A^{\rho_7}_1\hbar+A^{\rho_7}_2\hbar^2+\cdots\bigr)
\]
is a factorially divergent series with $A_{i}\in F$. The first two values $A_i$ (for more see Appendix~\ref{app:numerical}) are given by
\begin{gather}
 \begin{pmatrix}
 24\delta_{\rho_7}^3A^{\rho_7}_1\\
 1152\delta_{\rho_7}^6A^{\rho_7}_2
 \end{pmatrix}
  =
 \begin{pmatrix}
 1497746 & 3014838521575\\
 1345119 & 2732414541176\\
 -3675733 & -7414786842283\\
 2082815 & 4197826806919\\
 -839488 & -1690529009777\\
 -283405 & -574198051621\\
 383432 & 771765277669
 \end{pmatrix}^{\mathsf{T}}
 \begin{pmatrix}
 1\\
 \xi_7\\
 \xi_7^2\\
 \xi_7^3\\
 \xi_7^4\\
 \xi_7^5\\
 \xi_7^6
 \end{pmatrix} \nonumber\\
 \hphantom{\begin{pmatrix}
 24\delta_{\rho_7}^3A^{\rho_7}_1\\
 1152\delta_{\rho_7}^6A^{\rho_7}_2
 \end{pmatrix}}{}
  =
 \begin{pmatrix}
 -158.75\ldots - {\rm i}57.225\dots\\
 84862.\ldots + {\rm i}924.76\dots
 \end{pmatrix}.\label{eq:first.two.coeffs}
\end{gather}
Indeed, $X(1/2)=2$ and for $x=\frac{2\pi {\rm i}}{1000+1/2}$
\begin{align*}
 &\left|\WRT(\exp(x))\right|
 =1.6447\ldots\times 10^{98} ,\\
 &\left|\WRT(\exp(x))-2\e(3/8)\e\left(\frac{\VC_{\rho_7}}{2\pi {\rm i}x}\right)\sqrt{\frac{2\pi {\rm i}}{\delta_{\rho_7}x}}\right|
 =1.1735\ldots\times10^{95} ,\\
 &\left|\WRT(\exp(x))-2\e(3/8)\e\left(\frac{\VC_{\rho_7}}{2\pi {\rm i}x}\right)\sqrt{\frac{2\pi {\rm i}}{\delta_{\rho_7}x}}\bigl(1+A^{\rho_7}_1x\bigr)\right|
 =1.2530\ldots\times10^{92} ,\\
 &\left|\WRT(\exp(x))-2 \e(3/8) \e\left(\frac{\VC_{\rho_7}}{2\pi {\rm i}x}\right)\sqrt{\frac{2\pi {\rm i}}{\delta_{\rho_7}x}}\bigl(1+A^{\rho_7}_1x+A^{\rho_7}_2x^2\bigr)\right|
 =3.9674\ldots\times10^{89} .
\end{align*}
We will explore the function $X(r)$ in Section~\ref{qmodWRT}. However, those familiar with~\cite{Zagier:Qmod} will surely guess what the function is. The numbers $A_{i}$ can be computed from a formal Gaussian integration~\cite{Gang}, which is obtained from a conjectural surgery formula of the series discussed in~\cite{DimGar:QC}. Although these agree for this example, their topological invariance has not been established. The invariance of the series in~\cite{DimGar:QC} for knots was only recently proved in~\cite{GSW}. Recall the notation of~\cite{GSW},
\begin{gather*}
 \psi_{\hbar}(x,z)
 :=
 \exp\Biggl( -  \sum_{\underset{k+\frac{\ell}{2}>1}{k,\ell\in\BZ_{\geq0}}}
 \frac{B_{k} x^{\ell} \hbar^{k+\frac{\ell}{2}-1}}{\ell!k!}\Li_{2-k-\ell}(z)\Biggr)
 \in1+\hbar^{\frac{1}{2}}\BQ\bigl[z,(1-z)^{-1}\bigr][x]\bigl\llbracket\hbar^{\frac{1}{2}}
 \bigr\rrbracket ,
\end{gather*}
where $B_{k}=B_{k}(0)$ are the Bernoulli numbers and for $f_\hbar(x,z)\in\BQ(z)[x]\bigl\llbracket \hbar^{\frac12} \bigr\rrbracket$
\begin{align*}
\langle
f_\hbar(x,z)
\rangle_{x,\Lambda}
:={}&  \exp\Biggl(
{\frac 1 2 \sum_{i,j=1}^N \bigl(\Lambda^{-1}\bigr)_{i,j}
\frac{\partial}{\partial x_i}
\frac{\partial}{\partial x_j}}\Biggr)
f_\hbar(x,z) \Biggl|_{x=0}\\
  ={}&
\frac{\int e^{-\frac12x^t \Lambda \,x}f_{\hbar,z}(x) \, {\rm d}x}{
 \int e^{-\frac12x^t \Lambda \,x} \, {\rm d}x}
\in \BQ(z)\llbracket \hbar \rrbracket .
\end{align*}
Using this, we can define\footnote{There is an ambiguity of a sign due to the $\sqrt{\delta}$ but this can be chosen and then used consistently.} for each embedding of the trace field $\rho_1,\dots,\rho_7$ and corresponding solutions~\eqref{eq:Xsols} to equations~\eqref{eq:Xequs}
\begin{gather*}
  \Phi^{(\rho_{j})}_{m}(\hbar)
:=
 z_{j}^m\frac{\exp(-m\hbar)}{\sqrt{\delta_{\rho_j}}}
 \bigl\langle
 \exp\bigl(-mx_1\hbar^{\frac{1}{2}}+x_2\hbar^{\frac{1}{2}}\bigr)
 \psi_{\hbar}\bigl(2x_1,X_{1,j}^2\bigr)
 \psi_{\hbar}(x_2,X_{2,j})\\[1mm]
 \hphantom{\Phi^{(\rho_{j})}_{m}(\hbar):= z_{j}^m\frac{\exp(-m\hbar)}{\sqrt{\delta_{\rho_j}}} \bigl\langle}
 {}\times
 \psi_{\hbar}\bigl(-x_1+x_2,X_{1,j}^{-1}X_{2,j}\bigr)
 \bigr\rangle_{(x_1,x_2),\Lambda} ,\\[1mm]
 \Lambda
 =
 \begin{pmatrix}
 -2+2\Li_{0}\bigl(X_{1,j}^2\bigr)+\Li_{0}\bigl(X_{1,j}^{-1}X_{2,j}\bigr) & -1-\Li_{0}\bigl(X_{1,j}^{-1}X_{2,j}\bigr)\\[1mm]
 -1-\Li_{0}\bigl(X_{1,j}^{-1}X_{2,j}\bigr) & \Li_{0}(X_{2,j})+\Li_{0}\bigl(X_{1,j}^{-1}X_{2,j}\bigr)
 \end{pmatrix} .
\end{gather*}
This corresponds to formally applying the saddle point method to the integral of equation~\eqref{eq:state.int}. We can also define
\begin{align*}
 \widehat{\Phi}^{(\rho_0)}_{m}(\hbar) :={}& \WRT_{m,0}(\exp(\hbar))\\
  ={}&
 -\hbar - \frac{25}{2}\hbar^2 + \biggl(-12m - \frac{1621}{6}\biggr)\hbar^3 + \biggl(-6m^2 - 510m - \frac{195601}{24}\biggr)\hbar^4+\cdots ,
\end{align*}
and $V_{\rho_0}=0$, the asymptotic series associated to the trivial connection $\rho_0$. Using similar notation to~\cite{GZ:RQMOD}, we can take extended formal power series
\[
 \widehat{\Phi}^{(\rho_{j})}_{m}(\hbar)=\exp(\VC_{\rho_{j}}/\hbar)\Phi^{(\rho_{j})}_{m}(\hbar) .
\]
These extended series can be put into an $8\times 8$ matrix indexed by $m$ and $\rho$. This matrix gives a Wronskian solving the $q$-difference equation~\eqref{8thorderqdiff} when $q=e^{\hbar}$. To simplify notation, let
\[
\widehat{\Phi}^{(\rho_{j})}(\hbar)=\widehat{\Phi}^{(\rho_{j})}_{0}(\hbar).
\]

This collection of series is consistent with the conjectured arithmetic nature of perturbative quantum invariants for discrete character varieties~\cite{DGLZ,Gar:resCS}, which states that $\Phi^{(\rho_j)}$ for $j>0$ are all Galois conjugates of $\Phi^{(\rho_7)}$. These calculations give us a way to find the asymptotics of the $\WRT(\e(1/k))$ using Galois conjugation. Indeed, as predicted by Witten's asymptotic expansion conjecture, we have the following theorem, which to leading order appeared in~\cite{CharlesMarche} while to all orders expected in upcoming work of~\cite{AM:asymp}.
\begin{Theorem}[\cite{AM:asymp,CharlesMarche}]\label{thm:WAEC}
For $k\in\BZ$ as $k\rightarrow\infty$,
\[
 \WRT(\e(1/k))
 \sim
 \sum_{\rho\in\{\rho_0,\rho_1,\rho_2,\rho_4,\rho_5\}}\!\!\mu_{\rho}k^{d_\rho}\widehat{\Phi}^{(\rho)}(2\pi {\rm i}/k) ,
\]
where $d_{\rho_j}=1/2,\,\mu_{\rho_j}=(-1)^j\e(-1/8)$ and $\Phi^{(\rho_j)}$ are the Galois conjugates of $\Phi^{(\rho_7)}$, and finally a~term corresponding to the trivial connection on $M$, which has $d_{\triv}=0,\,\mu_{\rho_0}=1$.
\end{Theorem}
The comparison between the leading order and the series is given in Figure~\ref{fig:WAECpic}.
\begin{figure}[t]\centering
\vspace{-1mm}

\caption{Plots of the WRT invariant of $4_1(-1,2)$ against the first order approximation in Witten's asymptotic expansion conjecture where $N\in\BZ$ and $\e(x)=\exp(2\pi {\rm i}x)$.}\label{fig:WAECpic}\vspace{-2mm}
\end{figure}
Figure~\ref{fig:WAECpic} seems to be a good match however the relative error is actually only a decimal or two. However, a~more detailed numerical analysis can be applied. Indeed, we find that for $x=\frac{2\pi {\rm i}}{1000}$,
\begin{gather*}
 \left|\WRT(\exp(x))\right|
  = 16.099\dots,\\
 \left|\WRT(\exp(x))-      \sum_{\rho\in\{\rho_0,\rho_1,\rho_2,\rho_4,\rho_5\}}\mu_{\rho}\left(\frac{2\pi {\rm i}}{x}\right)^{d_\rho}\exp\left(\frac{\VC_{\rho}}{2\pi {\rm i}x}\right)a_0^{\rho}\right|
  = 0.019710\dots,\\
 \left|\WRT(\exp(x))-      \sum_{\rho\in\{\rho_0,\rho_1,\rho_2,\rho_4,\rho_5\}}\mu_{\rho}\left(\frac{2\pi {\rm i}}{x}\right)^{d_\rho}\exp\left(\frac{\VC_{\rho}}{2\pi {\rm i}x}\right)\bigl(a_0^{\rho}+a_1^{\rho}x\bigr)\right|
  = 0.00048434,\dots\\
 \left|\WRT(\exp(x))-     \!\!\! \sum_{\rho\in\{\rho_0,\rho_1,\rho_2,\rho_4,\rho_5\}}\!\!\!\!\mu_{\rho}\!\left(\frac{2\pi {\rm i}}{x}\right)^{d_\rho}\!\exp\!\left(\frac{\VC_{\rho}}{2\pi {\rm i}x}\right)\!\bigl(a_0^{\rho}+a_1^{\rho}x+a_2^{\rho}x^2\bigr)\right|\!
  = 6.5311\ldots\times10^{-5}.
\end{gather*}
The form of these volume conjectures make them seem closely related. We will find that this is indeed the case and the best description is through quantum modularity discussed in Section~\ref{sec:qmod}.

\subsection{Behaviour of the coefficients}
\label{asympser}

Using various extrapolation methods, I could compute around $110$ coefficients of the asymptotic series as exact numbers. With these values, we can compute the asymptotics of the coefficients. These asymptotic series are expected to be resurgent~\cite{Gar:resCS} and are therefore expected to have very precise behaviour. This behaviour allows the use of optimal truncation and using this we can find the first few sub-leading terms. Numerically, we find that for $a_{k}^{\rho}=A_{k}^{\rho}/\sqrt{\delta_{\rho}}$
\[
 a_{k}^{\rho_j}
 \ae
 \sum_{i\neq j}\frac{M_{j,i}}{2\pi {\rm i}}\sum_{\ell\geq0}\frac{\Gamma(k-\ell)}{(\VC_{\rho_j}-\VC_{\rho_i})^{k-\ell}}a_{\ell}^{\rho_i},\\
\]
where
\be\label{eq:closeststokes}
M
 =
\begin{pmatrix}
 ? & ? & ? & 0 & ? & 1 & 1\\
 ? & ? & ? & ? & 0 & -1 & -1\\
 0 & ? & ? & ? & ? & 1 & 1\\
 ? & ? & ? & ? & 0 & -1 & -1\\
 ? & ? & ? & ? & ? & 1 & 1\\
 ? & ? & ? & 1 & -1 & 0 & 0\\
 ? & ? & ? & 1 & -1 & 0 & 0\\
\end{pmatrix} .
\ee
The question marks appear as these contributions are too small to see numerically with something like optimal truncation. These will be computed via a different approach in Section~\ref{res}. We can perform a similar analysis to find that for the trivial connection
\[
 \WRT\bigl(e^{\hbar}\bigr)
  =
 \sum_{k=0}^{\infty}a_{k}^{\rho_0}\hbar^{k}
  =
 -\hbar - \frac{25}{2}\hbar^2 - \frac{1621}{6}\hbar^3 - \frac{195601}{24}\hbar^4+\cdots ,
\]
and numerically we observe that
\[
 a_{k}^{\rho_0}
 \ae
 \sum_{i=1}^{7}\frac{M_{0,i}}{\sqrt{-2\pi}}\sum_{\ell\geq0}\frac{\Gamma(k-\ell+1/2)}{(0-\VC_{i})^{k-\ell+1/2}}a_{\ell}^{\rho_i} ,
\]
where
\[
M_{0,\cdot}
 =
\begin{pmatrix}
 ? & ? & -1 & ? & 0 & 1 & 1
\end{pmatrix} .
\]
The leading order of this expansion was noticed originally by S. Garoufalidis in a note he gave me when I started work on this subject. These asymptotics illustrate the well-known lack of symmetry where the trivial connection sees the others connections while the others do not see the trivial.\footnote{This is not surprising from the perspective of the $q$-holonomic modules where the non-trivial series actually give rise to a rank seven $q$-holonomic module and the trivial is an inhomogeneous extension.}

As predicted by resurgence~\cite{Gar:resCS}, we see that the Borel transform of these series appears to be convergent with radius of convergence given by the difference between two volumes. Taking Pad\'e approximants along with these asymptotics then seem to give rise to functions with endless analytic continuation.

\subsection{The constant terms at all roots of unity}

As advocated in~\cite{DimGar:II,GZ:qser,GZ:RQMOD,Zagier:Qmod}, we can in fact extract more information. In particular, we get asymptotic series at each root of unity. These asymptotic series can be analogously computed using Gaussian integration.
However, we will not deal with detailed analysis of these coefficients to large order and they are expected to behave much in the same way as the coefficients for the series associated to $q=1$. However, we will be interested in the constants of these asymptotics series. These are constructed in a similar way to those in~\cite{GZ:RQMOD}. We use the solutions to
\begin{align}
 &0 =(1-X_2)\bigl(1-X_1^{-1}X_2\bigr)-X_1 ,\nonumber\\
 &0 =\bigl(1-X_1^2\bigr)^2-X_1^2X_2\bigl(1-X_1^{-1}X_2\bigr) .\label{eq:Xequs}
\end{align}
These have solutions with respect to our generator of the trace field
\begin{gather}
 X_{1,j} =-3+11\xi_{\rho_j}+20\xi_{\rho_j}^2-15\xi_{\rho_j}^3+6\xi_{\rho_j}^4-2\xi_{\rho_j}^5-4\xi_{\rho_j}^6 ,\nonumber\\
 X_{2,j} =-9+19\xi_{\rho_j}+76\xi_{\rho_j}^2-52\xi_{\rho_j}^3+20\xi_{\rho_j}^4-4\xi_{\rho_j}^5-13\xi_{\rho_j}^6 .\label{eq:Xsols}
\end{gather}
Then letting
\be\label{eq:Delta}
 \Delta
  =
 -257+806\xi+947\xi^2-749\xi^3+331\xi^4-133\xi^5-213\xi^6 ,
\ee
the constants of the asymptotic series can be computed (or here defined) for $\z=\e(N/M)$ and $j=1,\dots,7$ as follows:
\begin{gather}
 \Phi^{(\rho_j)}_{N/M,m}(0):=\frac{-{\rm i}}{M\etaroots(\z)}\frac{\bigl(1-X_{1,j}^2\bigr)^{3/2}(1-X_{2,j})\bigl(1-X_{1,j}^{-1}X_{2,j}\bigr)} {\sqrt{\Delta}\,D_{\z}(X_{2,j})D_{\z}\bigl(X_{1,j}^2\bigr)D_{\z}\bigl(X_{1,j}^{-1}X_{2,j}\bigr)}\nonumber\\
 \hphantom{\Phi^{(\rho_j)}_{N/M,m}(0):=}{}
 \times\sum_{k,\ell\in\BZ/M\BZ}
 \frac{\z^{k^2+k\ell-mk+\ell}X_{1}^{\frac{2k+\ell-m}{M}}X_{2,j}^{\frac{k+1}{M}}}{\bigl(\z X_{2,j}^{1/M};\z\bigr)_{\ell}\bigl(\z X_{1,j}^{2/M};\z\bigr)_{2k}\bigl(\z X_{1,j}^{-1/M}X_{2,j}^{1/M};\z\bigr)_{\ell-k}} ,\label{eq:phicont}
\end{gather}
where
\[
 \etaroots(\z)
  =
 \sqrt{-{\rm i}}\prod_{\ell=1}^{\mathrm{ord}(\z)-1}\bigl(1-\z^{\ell}\bigr)^{\frac{1}{2}-\frac{\ell}{\mathrm{ord}(\z)}}
 \qquad\text{and}\qquad
 D_{\z}(x)
  =
 \prod_{j=1}^{\mathrm{ord}\bigl(\z)}(1-\z^{j}x\bigr)^{\frac{j}{\mathrm{ord}(\z)}-\frac{1}{2}}
\]
is the multiplier system of the Dedekind $\eta$-function and the cyclic dilogarithm respectively. Then we can check the asymptotics of these coefficients to find they behave in the same way as those of the WRT invariant shown in equation~\eqref{qmodnum},
\be\label{eq:phiqmod}
 \Phi^{(\rho_j)}_{-1/x,m}(0)
 \ae\e\biggl({-}\frac{\VC_{\rho_j}}{(2\pi {\rm i})^2\denom(x)\numer(x)}\biggr)X^{(\rho_j)}(x)
 \Phih^{(\rho_6)}_{0,m}\biggl({-}\frac{2\pi {\rm i}}{x}\biggr) ,
\ee
for some periodic functions $X^{(\rho_j)}(x)=X^{(\rho_j)}(x+1)$. In this case, those familiar with~\cite{GZ:RQMOD} will already guess what these functions are while others will have to wait for Section~\ref{sec:RQMOD}.

\subsection[Radial asymptotics of $\widehat{Z}$]{Radial asymptotics of $\boldsymbol{\widehat{Z}}$}
\label{radzhat}

Based on extensive examples for non-hyperbolic manifolds~\cite{Cheng:3dmod,GukovMan,GPV,Hikami,LawrenceZag}, it has been conjectured that there should exist $q$-series invariants of closed manifolds referred to as $\widehat{Z}$ and that their radial asymptotics at roots of unity should recover the WRT invariant of the same manifold~\cite{GukovMan}. This asymptotic equivalence property does no uniquely define the $q$-series as there can of course be lower order corrections. An important point to note is that for non-hyperbolic manifolds all of the Chern--Simons values or complixified volumes at flat connections are real numbers. Therefore, radial asymptotics also determine the asymptotics on angles. This is not the case for manifolds with hyperbolic pieces and therefore we could expect new phenomenon. This appears already for the example $4_1(-1,2)$. For $4_1(-1,2)$, we saw in Proposition~\ref{prop:zhat} that up to a factor of $-q^{1/2}$ the series $\widehat{Z}(q)$ is given by
\[
 Z_{0,0}(q)
  =
 \sum_{0\leq k\leq\ell}(-1)^{k+\ell}q^{\frac{1}{2}3k(k+1)+\frac{1}{2}\ell(\ell+1)-k}\frac{(q;q)_{\ell}}{(q;q)_{2k}(q;q)_{\ell-k}} .
\]
We can apply similar methods used later in Section~\ref{statphase} and the outline of the proof of Theorem~\ref{thm:qmodold.zhat} to compute the radial asymptotics of $\widehat{Z}(q)$, that is the asymptotics for $\tau\in {\rm i}\BR$ as $\tau\rightarrow 0$. Doing this, we find the following theorem.\footnote{See Appendix~\ref{sec:asymp.meth.ap} for an outline of the proof.}
\begin{Theorem}\label{thm:zhatasym}
For $\tau\rightarrow {\rm i}\infty$ with $\tau\in {\rm i}\BR$ and $\tq=\e(-1/\tau)$ and $q=\e(\tau)$, we have
\[
 Z_{0,0}(\tq)
 \sim
 \e(3/8)\,\sqrt{\frac{\tau}{\delta_{\rho_3}}}\,\e\left(\frac{\VC_{\rho_3}-4\pi^2}{(2\pi {\rm i})^2}\tau\right)\Phi^{(\rho_{3})}\biggl(\frac{2\pi {\rm i}}{\tau}\biggr).
\]
\end{Theorem}
Indeed, we can compare with some numerics and find
\begin{align*}
 2\pi\frac{\log\bigl(Z_{0,0}\bigl(\e\bigl(10^{-7}{\rm i}\bigr)\bigr)\bigr)}{10^7}
 =0.11620\ldots,\qquad 4\pi^2-\VC_{\rho_3}
 =0.11620\ldots.
\end{align*}
We see that $1=o\bigl(q^{-1}\e\bigl(-\frac{\VC_{\rho_3}}{(2\pi {\rm i})^2}\tau\bigr)\bigr)$ and therefore we cannot\footnote{This is known to some experts and while this was given in detail here, numerical observations were also noted and communicated to me by S.~Gukov and M.~Mari\~no.} take the radial limit of $\widehat{Z}$.

\begin{figure}[t]\centering
\includegraphics{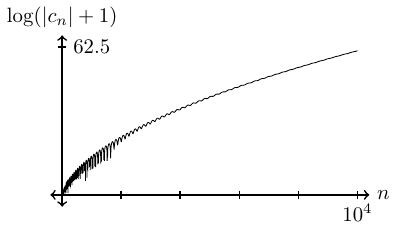}$\quad$
\includegraphics{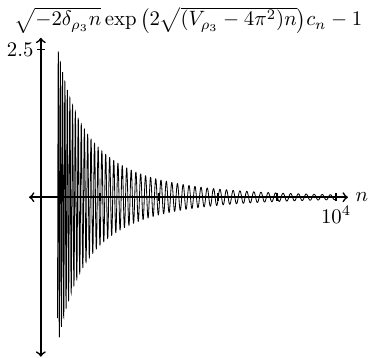}
\caption{Plots of the logarithm of the coefficients of $\widehat{Z}$ on the left and the error in the asymptotics on the right.}
\label{fig:zhatcoeffs}
\end{figure}

These asymptotics can be used to understand the asymptotics of the coefficients $\widehat{Z}(q)$ as a~$q$-series. Indeed, letting
\[
 Z_{0,0}(q)
  =
 \sum_{n=0}^{\infty}c_{n}q^{n}
\]
we can explore the behaviour of $c_{n}$. One can take the log of the first ten thousand coefficients and find behaviour like $\sqrt{n}$ as seen on the left of Figure~\ref{fig:zhatcoeffs}.
This behaviour is expected for $q$-series with asymptotics at roots of unity of the shape of Theorem~\ref{thm:zhatasym} and one can quickly find that
\[
 c_{n} \sim \frac{1}{\sqrt{-2\delta_{\rho_3}n}}\exp\Bigl(2\sqrt{(4\pi^2-V_{\rho_3})n}\Bigr)\bigl(1+O\bigl(n^{-1/2}\bigr)\bigr) ,
\]
as can be somewhat seen on the right of Figure~\ref{fig:zhatcoeffs}. This kind of statement can be proved using the asymptotics of $Z_{0,0}(q)$ at roots of unity and the circle method of Hardy--Ramanujan. However, we won't give details here. For those interested, see, for example,~\cite[Theorem 6.1]{GZ:asympnahm} or~\cite{Folsom}.

\section{Quantum modularity}\label{sec:qmod}

The volume conjecture of Kashaev~\cite{Kashaev:VC}, has been refined in many ways. One of the important refinements was the introduction of modularity~\cite{Zagier:Qmod}. Here the remarkable property was found that no matter how you approach a root of unity, quantum invariants will have the same asymptotic series however only after dividing by the same invariant evaluated at a root of unity obtained from a M\"obius transformation. More recently, this was refined~\cite{GZ:RQMOD} to the point that conjectural analytic lifts of the asymptotic series were proposed and related to Borel resummation~\cite{GGM:I}. This is well understood for knots and we can transfer much of this into the case of closed manifolds.

\subsection{Quantum modularity of the WRT invariant}
\label{qmodWRT}
\begin{figure}
\begin{center}
\includegraphics{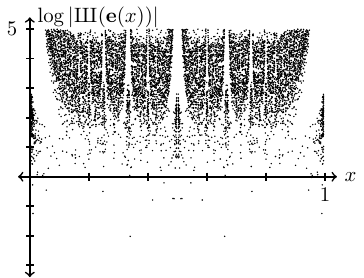}$\quad$
\includegraphics{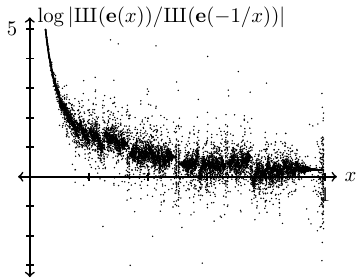}
\end{center}
\caption{Plots of the logarithm of the WRT invariant and its modular quotient for the first ten thousand rational numbers in $(0,1)$ and ordered by denominators.}
\label{fig:qmod.wrt.old}
\end{figure}

The asymptotics of the Kashaev invariant of a knot are best understood via the quantum modularity conjecture~\cite{Zagier:Qmod}. Therefore, let us start in a similar manner to~\cite{Zagier:Qmod} and list of a few values of this not so mysterious function $X(r)$ from equation~\eqref{qmodnum}. We have the following equalities:
\begin{gather*}
 X(1/2) =\WRT(\e(-1/2))  =2 ,\\
 X(1/3) =\WRT(\e(-1/3))  =1-\e(-1/3) ,\\
 X(1/4) =\WRT(\e(-1/4))  =1-\e(-1/4) ,\\
 X(1/5) =\WRT(\e(-1/5))  =2-\e(-1/5)+2\e(-2/5)+2\e(-3/5) ,\\
 X(1/6) =\WRT(\e(-1/6))  =1-\e(-1/6) ,\\
 X(1/7) =\WRT(\e(-1/7))  =2-2\e(-1/7)-\e(-2/7)+\e(-3/7) ,\\
 X(1/8) =\WRT(\e(-1/8))  =3-3\e(-1/8) .
\end{gather*}
This equality persists for $r$ with larger denominator and in fact we have the following theorem.
\begin{Theorem}\label{thm:qmodold}
The WRT invariant for $4_{1}(-1,2)$ is a quantum modular form in the sense of~{\rm \cite{Zagier:Qmod}}. That is, for $r\in\BQ\smallsetminus \BZ$ and $k\in\BZ$ as $k\rightarrow\infty$,
\[
 \WRT(\e(1/(k+r)))
 \sim
 \WRT(\e(-(k+r)))\e(3/8) \sqrt{\frac{k+r}{\delta_{\rho_7}}} \e\left(\frac{\VC_{\rho_7}}{(2\pi {\rm i})^2}(k+r)\right)\Phi^{(\rho_{7})}\biggl(\frac{2\pi {\rm i}}{k+r}\biggr) .
\]
\end{Theorem}
This theorem is proved in Section~\ref{statphase}. To get a feeling for this result, we can form similar plots to those in~\cite{Zagier:Qmod}. These are given in Figure~\ref{fig:qmod.wrt.old}.
Notice that compared with the plots in~\mbox{\cite[Figures~3 and~4]{Zagier:Qmod}} the second is behaving much worse than one could potentially hope. This is related to the fact that the volume of $4_1$ is larger than the volume of $4_1(-1,2)$. Notice that around $x=0$ the plot becomes better behaved. With a more detailed plot it would become clear that the right plot is smooth from the left and right of each rational point while also being discontinuous at each rational. We will see that the refined modularity~\cite{GZ:RQMOD} discussed in Section~\ref{sec:RQMOD} will give rise to functions that satisfy much better analytic properties and one can compare with Figure~\ref{fig:cocycleplots.wrt}, which is a similar figure to that given in~\cite[Figures~1 and~2]{GZ:RQMOD}.

Notice that Theorem~\ref{thm:qmodold} actually extends to the case that $r\in\BZ$. We know from Witten's asymptotic expansion conjecture that the left-hand side grows polynomially while the right-hand side vanishes as $\WRT(1)=0$. Therefore, in the $r=0$ case, Theorem~\ref{thm:qmodold} above states that the left-hand side is dominated by $\e\bigl(\VC_{\rho_7}/(2\pi {\rm i})^2k\bigr)$. Already the form of Witten's asymptotic expansion conjecture and this link suggests that there could be an improvement to this result by including sub-exponential corrections.

Following~\cite{GZ:RQMOD}, we can indeed develop a conjectural improvement of this result. To do this, we would like to turn the asymptotic series $\Phi^{(\rho)}$ into an analytic function. To do this we can use the sequence of Borel transformation, Pad\'e approximation and Laplace transformation~\cite{Caliceti}. Let~$s_{2N}(\Phi)$ be the Laplace transform of the $[N/N]$ Pad\'e approximate of the Borel transform. With $200$ coefficients of $\Phi^{(\rho)}$ we can numerically compute to around order $10^{-40}$ at $2\pi {\rm i}/100$. Indeed, the worst convergence is from $\Phi^{(\triv)}$ and the difference between the numerical values using~$200$ and $198$ coefficients is
\be\label{Bor.err}
 s_{200}\bigl(\Phi^{(\triv)}\bigr)(2\pi {\rm i}/100)-s_{198}\bigl(\Phi^{(\triv)}\bigr)(2\pi {\rm i}/100)
  =
 (4.8831\ldots + 3.0178\dots {\rm i})\times 10^{-40} ,
\ee
which give a lower bound for the numerical error. Numerically it appears roughly that
\begin{align*}
& \log|s_{N}\bigl(\Phi^{(\triv)}\bigr)(2\pi {\rm i}/100)-s_{500}\bigl(\Phi^{(\triv)}\bigr)(2\pi {\rm i}/100)|\sim -6.3\sqrt{N},\\
& \log|s_{N}\bigl(\Phi^{(\triv)}\bigr)(2\pi {\rm i}/200)-s_{500}\bigl(\Phi^{(\triv)}\bigr)(2\pi {\rm i}/200)|\sim -8.8\sqrt{N}.
\end{align*}
The values of the various numerical Borel--Pad\'e--Laplace transforms are given by
\begin{gather*}
 s_{200}\bigl(\Phi^{(\rho_0)}\bigr)(2\pi {\rm i}/100)  = 0.022362\ldots - 0.042136\dots {\rm i} ,\\
 \e(3/8)\sqrt{100}\exp \bigl(\VC_{\rho_1}100/2\pi {\rm i} \bigr)s_{200}\bigl(\Phi^{(\rho_1)}\bigr)(2\pi {\rm i}/100)  = -0.70748\ldots - 2.8524\dots {\rm i} ,\\
 \e(-1/8)\sqrt{100}\exp \bigl(\VC_{\rho_2}100/2\pi {\rm i} \bigr)s_{200}\bigl(\Phi^{(\rho_2)}\bigr)(2\pi {\rm i}/100)  = 1.0653\ldots - 2.6036\dots {\rm i} ,\\
 \e(3/8)\sqrt{100}\exp \bigl(\VC_{\rho_3}100/2\pi {\rm i} \bigr)s_{200}\bigl(\Phi^{(\rho_3)}\bigr)(2\pi {\rm i}/100)  = -0.95412\ldots + 0.53897\dots {\rm i} ,\\
 \e(-1/8)\sqrt{100}\exp \bigl(\VC_{\rho_4}100/2\pi {\rm i} \bigr)s_{200}\bigl(\Phi^{(\rho_4)}\bigr)(2\pi {\rm i}/100)  = 2.9304\ldots + 2.3781\dots {\rm i} ,\\
 \e(3/8)\sqrt{100}\exp \bigl(\VC_{\rho_5}100/2\pi {\rm i} \bigr)s_{200}\bigl(\Phi^{(\rho_5)}\bigr)(2\pi {\rm i}/100)  = -0.55386\ldots - 4.2697\dots {\rm i} ,\\
 \e(3/8)\sqrt{100}\exp \bigl(\VC_{\rho_6}100/2\pi {\rm i} \bigr)s_{200}\bigl(\Phi^{(\rho_6)}\bigr)(2\pi {\rm i}/100)  = (1.0359\ldots + 0.31226\dots {\rm i})\times 10^{-9},\\
 \e(3/8)\sqrt{100}\exp \bigl(\VC_{\rho_7}100/2\pi {\rm i} \bigr)s_{200}\bigl(\Phi^{(\rho_7)}\bigr)(2\pi {\rm i}/100)  = (2.2618\ldots + 0.58777\dots {\rm i})\times 10^{10}.
\end{gather*}
If we sum over the $\SU(2)$ flat connections as predicted by Witten's conjecture (with $\mu_\rho$, $d_\rho$ as in Theorem~\ref{thm:WAEC}), we find that
\begin{gather*}
 \WRT(\e(1/100))
  = 2.7567\ldots - {\rm i}7.3897\dots,\\
 \WRT(\e(1/100))
 -
 \sum_{\rho\in\{\rho_0,\rho_1,\rho_2,\rho_4,\rho_5\}}  \mu_{\rho}100^{d_\rho}\exp(\VC_{\rho}100/2\pi {\rm i})s_{200}\bigl(\Phi^{(\rho)}\bigr)(2\pi {\rm i}/100)\\
 \qquad{} =
 1.0359\dots 10^{-9} + {\rm i}3.1226\dots 10^{-10} .
\end{gather*}
Therefore, we find that additionally summing over the exponentially small contribution connection,
\begin{gather}
 \WRT(\e(1/100))
 -
 \sum_{\rho\in\{\rho_0,\rho_1,\rho_2,\rho_4,\rho_5,\rho_6\}}\mu_{\rho}100^{d_\rho}\exp(\VC_{\rho}100/2\pi {\rm i})s_{200}\bigl(\Phi^{(\rho)}\bigr)(2\pi {\rm i}/100)\nonumber\\
 \qquad{} =
 6.1973\times10^{-39} + {\rm i}1.8705\times10^{-38} .\label{WAEC.Bord.2}
\end{gather}
Notice that equation~\eqref{WAEC.Bord.2} vanishes to the order of the error in equation~\eqref{Bor.err}. These numerics lead to the following conjecture where $s\bigl(\Phi^{(\rho)}\bigr)$ is the Borel resummation of $\Phi^{(\rho)}$.
\begin{Conjecture}
For $k\in\BZ_{>0}$,
\[
 \WRT(\e(1/k))
  =
 \sum_{\rho\in\{\rho_0,\rho_1,\rho_2,\rho_4,\rho_5,\rho_6\}}\mu_{\rho}k^{d_\rho}\exp\left(\VC_{\rho}k/2\pi {\rm i}\right)s\bigl(\Phi^{(\rho)}\bigr)(2\pi {\rm i}/k) .
\]
\end{Conjecture}
This seems to be in contrast to Witten's asymptotic expansion conjecture where we only expect $\mathrm{SU}(2)$ flat connections to appear. However, as already noted in~\cite{Gar:resCS} and seen in Section~\ref{asympser}, the asymptotic series associated to the $\mathrm{SU}(2)$ flat connections contain contributions from the other $\SL_{2}(\BC)$ flat connections, which leads to the series $\rho_6$ appearing as an exponentially small correction (which of course does not contradict the conjecture).

With some more experimentation one can naturally extend this conjecture as follows.
\begin{Conjecture}\label{conj.qmod.hab.Borel}
For $k\in\BZ_{>0}$ and $r\in\BQ_{\geq 0}$, there exist functions $X_{\rho}\colon \BQ\rightarrow\BC$ such that
\[
 \WRT\biggl(\e\biggl(\frac{1}{k+r}\biggr)\biggr)
  =
 \sum_{\rho}X_\rho(r)\mu_{\rho}(k+r)^{d_\rho}\exp\biggl(\frac{\VC_{\rho}(k+r)}{2\pi {\rm i}}\biggr)s\bigl(\Phi^{(\rho)}\bigr)\biggl(\frac{2\pi {\rm i}}{k+r}\biggr) ,
\]
and $X_{\rho}(r)\in\BZ[\e(-r)]$.
\end{Conjecture}
The first few values of the functions $X_{\rho}(r)$ are given in the following table:
\begin{center}
\def\arraystretch{1.2}
\begin{tabular}{|l|l|l|l|l|l|}\hline
$r$ & $1$ & $1/2$ & $1/3$ & $1/4$ & $1/5$ \\ \hline
$X_{\rho_0}(r)$ & $1$ & $1$ & $1$ & $1$ & $1$ \\ \hline
$X_{\rho_1}(r)$ & $1$ & $1$ & $-1-\e\bigl(\frac{1}{3}\bigr)$ & $1$ & $-2-2\e\bigl(\frac{2}{5}\bigr)$ \\ \hline
$X_{\rho_2}(r)$ & $1$ & $1$ & $-1-\e\bigl(\frac{1}{3}\bigr)$ & $-1-2{\rm i}$ & $-2-\e\bigl(\frac{1}{5}\bigr)-\e\bigl(\frac{3}{5}\bigr)$ \\ \hline
$X_{\rho_3}(r)$ & $0$ & $0$ & $1+2\e\bigl(\frac{1}{3}\bigr)$ & $-2{\rm i}$ & $-1-\e\bigl(\frac{1}{5}\bigr)-2\e\bigl(\frac{2}{5}\bigr)-\e\bigl(\frac{3}{5}\bigr)$ \\ \hline
$X_{\rho_4}(r)$ & $1$ & $1$ & $2+2\e\bigl(\frac{1}{3}\bigr)$ & $-1+2{\rm i}$ & $1-\e\bigl(\frac{1}{5}\bigr)+\e\bigl(\frac{2}{5}\bigr)$ \\ \hline
$X_{\rho_5}(r)$ & $1$ & $1$ & $2-\e\bigl(\frac{1}{3}\bigr)$ & $1$ & $3+\e\bigl(\frac{1}{5}\bigr)+\e\bigl(\frac{2}{5}\bigr)+\e\bigl(\frac{3}{5}\bigr)$ \\ \hline
$X_{\rho_6}(r)$ & $1$ & $-1$ & $1$ & ${\rm i}$ & $2-\e\bigl(\frac{3}{5}\bigr)$ \\ \hline
$X_{\rho_7}(r)$ & $0$ & $2$ & $2+\e\bigl(\frac{1}{3}\bigr)$ & $1+{\rm i}$ & $3+\e\bigl(\frac{1}{5}\bigr)+3\e\bigl(\frac{2}{5}\bigr)+3\e\bigl(\frac{3}{5}\bigr)$ \\ \hline
\end{tabular}
\end{center}

Assuming this conjecture, one can numerically compute the values of $X_{\rho}(r)$ for rational numbers with denominators up to a few hundred. With this data one can recognise these functions as combinations of $\WRT_{m,n}(\e(-r))$. However, these formulae can be quite complicated if not in the correct form. Indeed, it is better to study the function $\WRT$ in more detail to guess the correct formulae. This is done in the next Section~\ref{statphase}.

\subsection{Stationary phase and proof of Theorem~\ref{thm:qmodold}}
\label{statphase}
This section will describe the proof of the quantum modularity conjecture for $4_{1}(-1,2)$. For knots, Bettin and Drappeau gave a proof of the quantum modularity conjecture for a collection of simple hyperbolic knots in~\cite{BetDrap}. We can apply a similar method to prove quantum modularity of the WRT invariant. This can further be used to guess the functions appearing with the conjectural Borel resummation formulae of Conjecture~\ref{conj.qmod.hab.Borel}. Firstly, we can rewrite the WRT invariant as a sum over the positive cone as
\begin{align*}
 \WRT(q)
 &{}=
 \sum_{0\leq\ell\leq k}(-1)^{k}q^{-\frac{1}{2}k(k+1)+\ell(\ell+1)}\frac{(q;q)_{2k+1}}{(q;q)_{\ell}(q;q)_{k-\ell}}\\
 &{}=
 \sum_{\ell,j=0}^{\infty}(-1)^{j+\ell}q^{-\frac{1}{2}j(j+1)-j\ell+\ell(\ell+1)/2}\frac{(q;q)_{2j+2\ell+1}}{(q;q)_{\ell}(q;q)_{j}} .
\end{align*}
Then consider the elements
\[
 w_{m,n,p}(q)
  =
 \sum_{\ell,j=0}^{\infty}(-1)^{j+\ell}q^{-\frac{1}{2}j(j+1)-j\ell+\ell(\ell+1)/2+mj+n\ell}\frac{(q;q)_{2j+2\ell+p}}{(q;q)_{\ell}(q;q)_{j}} ,
\]
where
\[
 w_{m,n+m,1}(q) =\WRT_{m,n}(q) .
\]
Using the quantum modular properties of the $q$-Pochhammer symbol, i.e., the analytic properties of the Faddeev quantum dilogarithm discussed a little in Appendix~\ref{sec:mod.poch}, we can analyse the asymptotic behaviour of $\WRT(q)$. The main idea is to rewrite the Pochhammer at $(\tq;\tq)_{n}$ as $(q;q)_{k}$ times a functions with good analytic properties and known asymptotics. In particular, for $r=a/c\in\BQ$ and $f$ from equation~\eqref{eq:f}, we have
\begin{align*}
 &\WRT(\e(-1/r))
 =\sum_{0\leq\ell,0\leq j}^{2\ell+2j\leq a}(-1)^{j+\ell}
 \frac{\e\bigl(\frac{1}{2r}j(j+1)+\frac{j\ell}{r}-\frac{1}{2r}\ell(\ell+1)\bigr)
 \bigl(\e\bigl(-\frac{1}{r}\bigr);\e\bigl(-\frac{1}{r}\bigr)\bigr)_{2j+2\ell+1}}
 {\bigl(\e\bigl(-\frac{1}{r}\bigr);\e\bigl(-\frac{1}{r}\bigr)\bigr)_{\ell}\bigl(\e\bigl(-\frac{1}{r}\bigr);\e\bigl(-\frac{1}{r}\bigr)\bigr)_{j}}\\
 & =\sum_{0\leq\ell,0\leq j}^{2\ell+2j\leq a} \!\!\!(-1)^{j+\ell}
 \frac{\e\bigl(\frac{r}{2}\bigl(\frac{j}{r}\bigr)^2\!+\!r\frac{j}{r}\frac{\ell}{r}\!-\!\frac{r}{2}\bigl(\frac{\ell}{r}\bigr)^2\!+\!\frac{j}{2r}\!-\!\frac{\ell}{2r}\bigr)
 (\e(r);\e(r))_{\lfloor\frac{2j+2\ell+1}{r}\rfloor}
 \e(f(2j+2\ell+1,r))}
 {(\e(r);\e(r))_{\lfloor\frac{\ell}{r}\rfloor}(\e(r);\e(r))_{\lfloor\frac{ j}{r}\rfloor}
 \e(f(\ell,r))\e(f(j,r))}\\
 & =\sum_{0\leq\ell,0\leq j}^{2\ell+2j\leq a}
 \e\Biggl(\frac{r}{2}\left(\frac{j}{r}-\left\lfloor \frac{j}{r}\right\rfloor\right)^2+r\left(\frac{j}{r}-\left\lfloor \frac{j}{r}\right\rfloor\right)\left(\frac{\ell}{r}-\left\lfloor \frac{\ell}{r}\right\rfloor\right)-\frac{r}{2}\left(\frac{\ell}{r}-\left\lfloor \frac{\ell}{r}\right\rfloor\right)^2\\
& \quad{}+\frac{r}{2}\left(\frac{j}{r}-\left\lfloor\frac{j}{r}\right\rfloor\right)
 -\frac{r}{2}\left(\frac{\ell}{r}-\left\lfloor\frac{\ell}{r}\right\rfloor\right)\Biggr)\\
 &\quad{}\times(-1)^{\lfloor\frac{j}{r}\rfloor+\lfloor\frac{\ell}{r}\rfloor}\e\left(\frac{-r}{2}\left\lfloor\frac{j}{r}\right\rfloor^2-r\left\lfloor \frac{j}{r}\right\rfloor\left\lfloor\frac{\ell}{r}\right\rfloor+\frac{r}{2}\left\lfloor\frac{\ell}{r}\right\rfloor^2+\frac{r}{2}\left\lfloor\frac{j}{r}\right\rfloor-\frac{r}{2}\left\lfloor\frac{\ell}{r}\right\rfloor\right)\\
 &\quad{}\times \e\left(\frac{1}{2}\left(\frac{j}{r}-\left\lfloor\frac{j}{r}\right\rfloor\right)-\frac{1}{2}\left(\frac{\ell}{r}-\left\lfloor\frac{\ell}{r}\right\rfloor\right)\right)
 \frac{(\e(r);\e(r))_{\lfloor\frac{2j+2\ell+1}{r}\rfloor}
 \e(f(2j+2\ell+1,r))}
 {(\e(r);\e(r))_{\lfloor\frac{\ell}{r}\rfloor}(\e(r);\e(r))_{\lfloor\frac{ j}{r}\rfloor}
 \e(f(\ell,r))\e(f(j,r))}\\
 & =\sum_{R=0}^{3}\sum_{0\leq L,0\leq J}^{2L+2J+R<c}(-1)^{J+L}\e\left(\frac{-r}{2}J^2-rJL+\frac{r}{2}L^2+\frac{r}{2}J-\frac{r}{2}L\right) \frac{(\e(r);\e(r))_{2J+2L+R}}{(\e(r);\e(r))_{L}(\e(r);\e(r))_{J}}\\
 &\quad{}\times\sum_{(x,y)\in[0,1)^2\cap\left(\frac{1}{r}\BZ^2-(J,L)\right)}^{R\leq2x+2y\leq R+1}
 \frac{\e(f(2rx+2ry+1,r))}{\e(f(rx,r))\e(f(ry,r))}\\
 &\quad{}\times \e\left(\frac{r}{2}x^2+rxy-\frac{r}{2}y^2+\frac{r}{2}x-\frac{r}{2}y+\frac{1}{2}x-\frac{1}{2}y\right),
\end{align*}
where we used the substitution
\[
 J =\biggl\lfloor\frac{j}{r}\biggr\rfloor,\qquad
 L =\biggl\lfloor\frac{\ell}{r}\biggr\rfloor,\qquad
 x =\frac{j}{r}-\biggl\lfloor \frac{j}{r}\biggr\rfloor,\qquad
 y =\frac{\ell}{r}-\biggl\lfloor \frac{\ell}{r}\biggr\rfloor ,
\]
and the fact, from equation~\eqref{eq:f.per}, that for $n\in\BZ$ we have
\[
 f(rx+nr,r) =f(rx,r).
\]
Importantly,
\[
 \e(rx) =\e\biggl(-r\bigg\lfloor\frac{j}{r}\bigg\rfloor\biggr) =\e(-rJ),
 \qquad
 \e(ry) =\e\biggl(-r\bigg\lfloor\frac{\ell}{r}\bigg\rfloor\biggr) =\e(-rL) .
\]
Therefore, we can exchange these terms between the sums over $L$, $J$ and the sums over~$x$,~$y$. In~particular, for $(\alpha,\beta)\in\BZ^2$ we have
\begin{align*}
\WRT(\e(-1/r)) ={}&\sum_{R=0}^{3}\sum_{0\leq L,0\leq J}^{2L+2J+R<c}(-1)^{J+L}\e\left(\frac{-r}{2}J^2-rJL+\frac{r}{2}L^2+\frac{r}{2}J-\frac{r}{2}L+arJ+brL\right)\\
 & \times \frac{(\e(r);\e(r))_{2J+2L+R}}{(\e(r);\e(r))_{L}(\e(r);\e(r))_{J}}
 \sum_{(x,y)\in[0,1)^2\cap\left(\frac{1}{r}\BZ^2-(J,L)\right)}^{R\leq2x+2y\leq R+1}
 \frac{\e(f(2rx+2ry+1,r))}{\e(f(rx,r))\e(f(ry,r))}\\
& \times \e\left(\frac{r}{2}x^2+rxy-\frac{r}{2}y^2+\frac{r}{2}x-\frac{r}{2}y+\frac{1}{2}x-\frac{1}{2}y+\alpha rx+\beta ry\right) .
\end{align*}
This form of the expression makes it clear the asymptotic expansions come from a stationary phase approximation applied to the sum
\begin{align}
& \sum_{(x,y)\in[0,1)^2\cap\left(\frac{1}{r}\BZ^2-(J,L)\right)}^{R\leq2x+2y\leq R+1}
 \e\biggl(f(2rx+2ry+1,r)-f(rx,r)-f(ry,r)+\frac{r}{2}x^2+rxy\nonumber\\
& \qquad{}-\frac{r}{2}y^2+\frac{r}{2}x-\frac{r}{2}y+\frac{1}{2}x-\frac{1}{2}y+\alpha rx+\beta ry\biggr).\label{eq:dostat}
\end{align}
As $r$ tends to infinity, this sum is approximated\footnote{See the proof of Theorem~\ref{thm:qmodold}.} by the integral
\begin{align*}
 &\int\int_{R}\e\biggl(f(2rx+2ry+1,r)-f(rx,r)-f(ry,r)+\frac{r}{2}x^2+rxy\\
& \qquad {}-\frac{r}{2}y^2+\frac{r}{2}x-\frac{r}{2}y+\frac{1}{2}x-\frac{1}{2}y+\alpha rx+\beta ry\biggr)  {\rm d}x{\rm d}y ,
\end{align*}
where $R$ indexes the regions depicted in Figure~\ref{fig:R.regions}.
\begin{figure}[t]\centering
\begin{tikzpicture}[scale=0.8,baseline=-3]
\draw[<->,thick] (-0.5,0) -- (4.5,0);
\draw[<->,thick] (0,-0.5) -- (0,4.5);
\filldraw (0,4) circle (2pt);
\filldraw (0,0) circle (2pt);
\filldraw (4,0) circle (2pt);
\filldraw (4,4) circle (2pt);
\fill[blue,opacity=0.2] (0,0) -- (2,0) -- (0,2) -- cycle;
\fill[purple,opacity=0.2] (2,0) -- (0,2) -- (0,4) -- (4,0) -- cycle;
\fill[red,opacity=0.2] (0,4) -- (4,0) -- (4,2) -- (2,4) -- cycle;
\fill[orange,opacity=0.2] (4,2) -- (2,4) -- (4,4) -- cycle;
\draw (0,0) -- (4,0) -- (4,4) -- (0,4) -- cycle;
\draw (2,0) -- (0,2);
\draw (4,0) -- (0,4);
\draw (4,2) -- (2,4);
\draw (0.5,0.5) node {$R_0$};
\draw (1.5,1.5) node {$R_1$};
\draw (2.5,2.5) node {$R_2$};
\draw (3.5,3.5) node {$R_3$};
\end{tikzpicture}
\caption{The regions $R$ that arise from the stationary phase approximation.}\label{fig:R.regions}
\end{figure}
Noting that from equation~\eqref{eq:asymp.f} we have
\[
 f(rx,r) =\frac{r}{(2\pi {\rm i})^2}\Li_2(\e(-x))+\frac{r}{24}+O\bigl(r^0\bigr) ,
\]
we see that the critical points of the integrand are determined by the stationary points of
\begin{align}
 V_{\alpha,\beta}(x,y)&{}=\frac{1}{(2\pi {\rm i})^2}\Li_2(\e(-2x-2y))-\frac{1}{(2\pi {\rm i})^2}\Li_2(\e(-x))-\frac{1}{(2\pi {\rm i})^2}\Li_2(\e(-y))\nonumber\\
 &\quad{}-\frac{1}{24}+\frac{1}{2}x^2+xy-\frac{1}{2}y^2+\frac{1}{2}x-\frac{1}{2}y+\alpha x+\beta y .\label{eq:Vab}
\end{align}
These satisfy equations
\begin{align}
 &0 =\frac{2}{2\pi {\rm i}}\log(1-\e(-2x-2y))-\frac{1}{2\pi {\rm i}}\log(1-\e(-x))+x+y+\frac{1}{2}+\alpha ,\nonumber\\
 &0 =\frac{2}{2\pi {\rm i}}\log(1-\e(-2x-2y))-\frac{1}{2\pi {\rm i}}\log(1-\e(-y))+x-y-\frac{1}{2}+\beta .\label{eq.logstatphase}
\end{align}
Letting $X=\e(x)$, $Y=\e(y)$, we find equations
\begin{align}
 \frac{\bigl(1-X^{-2}Y^{-2}\bigr)^2}{\bigl(1-X^{-1}\bigr)}XY =-1 ,\qquad
\frac{\bigl(1-X^{-2}Y^{-2}\bigr)^2}{\bigl(1-Y^{-1}\bigr)}XY^{-1} =-1 .\label{eq:xysp}
\end{align}
The solutions to these equations are given by
\begin{align*}
 &X =10-19\xi-76\xi^2+52\xi^3-20\xi^4+4\xi^5+13\xi^6,\\
 &Y =-10+24\xi+87\xi^2-59\xi^3+23\xi^4-5\xi^5-15\xi^6 .
\end{align*}
Given any logarithm for $(X,Y)$ we can choose an $\alpha$, $\beta$ which give rise to a solution to equation~\eqref{eq.logstatphase}. For example, consider the solution associated to the 3rd embedding of the trace field
\[
 (x_3,y_3) =(0.0081372\dots {\rm i}, 1.0000\ldots + 0.0090947\dots {\rm i}) .
\]
We find that
\begin{align*}
& 2+\alpha_3 =\frac{2}{2\pi {\rm i}}\log(1-\e(-2x_3-2y_3))-\frac{1}{2\pi {\rm i}}\log(1-\e(-x_3))+x_3+y_3+\frac{1}{2}+\alpha_3,\\
& -1+\beta_3 =\frac{2}{2\pi {\rm i}}\log(1-\e(-2x_3-2y_3))-\frac{1}{2\pi {\rm i}}\log(1-\e(-y_3))+x_3-y_3-\frac{1}{2}+\beta_3 .
\end{align*}
Choosing $(\alpha_3,\beta_3)=(-2,1)$, we find that
\[
 V_{-2,1}(x_3,y_3)-\frac{1}{4\pi^2}\VC_3 =-1 .
\]
Therefore, associated to this solution to equation~\eqref{eq.logstatphase} we get four elements of the Habiro ring
\[
 w^{(x_3,y_3,R)}(q)
  =
 q^{V_{-2,1}(x_3,y_3)-\frac{1}{4\pi^2}\VC_3}w_{\alpha_3,\beta_3,R}(q)
  =
 q^{-1}w_{-2,1,R}(q) .
\]
Considering small solutions to equation~\eqref{eq.logstatphase}, we can numerically find formulae for $X_{\rho}$ summarised in the following conjecture.

\begin{Conjecture}
The functions $X_{\rho}$ in Conjecture~{\rm \ref{conj.qmod.hab.Borel}} are given by the following elements of the Habiro ring
\begin{gather}
 X_{\rho_1}(q)  =w_{0,2,0}(q) ,\nonumber\\
 X_{\rho_2}(q)  =w_{-1,0,0}(q) ,\nonumber\\
 X_{\rho_3}(q)  =q^{-1}w_{-1,0,1}(q)+w_{-1,2,1}(q)+q^{-1}w_{-2,1,3}(q) ,\nonumber\\
 X_{\rho_4}(q)  =w_{0,1,0}(q) ,\nonumber\\
 X_{\rho_5}(q)  =w_{-1,1,0}(q) ,\nonumber\\
 X_{\rho_6}(q)  =w_{-1,1,0}(q)+q^{-1}w_{-2,0,1}(q) ,\nonumber\\
 X_{\rho_7}(q)  =w_{-1,1,1}(q) .\label{eq:can.basis}
\end{gather}
\end{Conjecture}

These elements can all be written in terms the $\WRT_{m}(q)=\WRT_{m,0}(q)$. Firstly, let{\setlength{\arraycolsep}{2.5pt}
\begin{align}
&P(q)(1-q)q^4 =\label{eq:can.basis.P}\\
&\begin{pmatrix}
1 & q^3 + 2q^2 & \begin{array}{@{}c@{}} -q^5 + 2q^4 \\  - q^3 - q^2\end{array} & \begin{array}{@{}c@{}} -q^7 - q^6 \\   - 2q^5 - q^4\end{array} & q^8 + 2q^5 & q^{10} + q^9 + q^7 & -q^{10} - q^8 & -q^{12}\vspace{0.5mm}\\
1 & q^3 + 2q^2 & q^4 - q^3 - q^2 & -q^6 - 2q^5 - 2q^4 & q^7 + 2q^5 & q^9 + q^8 + q^7 & -q^9 - q^7 & -q^{11}\vspace{0.5mm}\\
q & q^4 + q^3 + q^2 & -q^3 & \begin{array}{@{}c@{}} -q^7 - 2q^6 \\  - q^5 - q^4\end{array} & q^8 + 2q^6 & q^{10} + q^9 + q^8 & -q^{10} - q^8 & -q^{12}\vspace{0.5mm}\\
0 & -q^2 + q & -q^4 + q^3 & q^5 - q^3 & -q^5 + q^4 & -q^7 + q^6 & 0 & 0\vspace{0.5mm}\\
1 & 2q^3 + q^2 & q^5 - q^3 - q^2 & -q^7 - q^6 - 3q^5 & 2q^6 + q^5 & \begin{array}{@{}c@{}} q^{10} + q^9 + q^8 \\ {}+ q^7 - q^6\end{array} & -q^9 - q^8 & -q^{12}\vspace{0.5mm}\\
1 & 2q^3 + q^2 & \begin{array}{@{}c@{}} q^5 + q^4\\   - 2q^3 - q^2\end{array} & -q^7 - q^6 - 3q^5 & q^6 + 2q^5 & \begin{array}{@{}c@{}}q^{10} + q^9 + q^8 \\  + q^7 - q^6\end{array} & -q^9 - q^8 & -q^{12}\vspace{0.5mm}\\
0 & 0 & 0 & -q^5 + q^4 & 0 & 0 & 0 & 0\vspace{0.5mm}\\
1 & 2q^3 + q^2 & q^5 - q^3 - q^2 & -q^7 - q^6 - 3q^5 & q^6 + 2q^5 & \begin{array}{@{}c@{}} q^{10} + q^9 + q^8 \\  + q^7 - q^6\end{array} & -q^9 - q^8 & -q^{12}
\end{pmatrix}\!.\nonumber
\end{align}}%
Then we have
\begin{gather}
\begin{pmatrix}
 X_{\rho_0}(q) &
 X_{\rho_1}(q) &
 X_{\rho_2}(q) &
 X_{\rho_3}(q) &
 X_{\rho_4}(q) &
 X_{\rho_5}(q) &
 X_{\rho_6}(q) &
 X_{\rho_7}(q)
\end{pmatrix}^{t}\label{eq:can.basis.Q}\\
\qquad{} =
P(q)
\begin{pmatrix}
 \WRT_{-3}(q) &
 \WRT_{-2}(q) &
 \WRT_{-1}(q) &
 \WRT_{0}(q) &
 \WRT_{1}(q) &
 \WRT_{2}(q) &
 \WRT_{3}(q) &
 \WRT_{4}(q)
\end{pmatrix}^{t}.\nonumber
\end{gather}
\begin{proof}[Proof of Theorem~\ref{thm:qmodold}]

We need to use standard analysis of asymptotic properties of sums. The first step is to apply a summation method to express the sum as an integral. Then one applies a stationary phase approximation to the integral. The main thing that needs to be checked is that the contributions from the boundaries are negligible and that one can deform the original contour of integration to a contour suitable for the stationary phase approximation.

Firstly, recall $V_{\alpha,\beta}(x,y)$ from equation~\eqref{eq:Vab}. Then, using standard properties of the Bloch--Wigner dilogarithm, as $x\rightarrow\infty$
\[
 \Im(\Li_2(\e(x)))
  =
 \arg(1-\e(x))\log|\e(x)|+O(1) .
\]
Therefore, breaking into the cases of $\Im(x)$, $\Im(y)$ tending to $\pm\infty$, we can show that for various limits $\e(V_{\alpha,\beta}(x,y)\pm x\pm y)\rightarrow 0$, where the extra $x$, $y$ come with signs so that they have positive imaginary part. Therefore, letting
\begin{align*}
& s_{\alpha,\beta}(x,y,r)  =
 \e\biggl(f(2rx+2ry+1,r)-f(rx,r)-f(ry,r)+\frac{r}{2}x^2+rxy\\
 &\hphantom{s_{\alpha,\beta}(x,y,r)  =\biggl(}{}
 -\frac{r}{2}y^2+\frac{r}{2}x-\frac{r}{2}y+\frac{1}{2}x-\frac{1}{2}y+\alpha rx+\beta ry\biggr) ,
\end{align*}
we can apply the Abel--Plana summation method~\cite[Section 8.3]{Olver} to the sum from equation~\eqref{eq:dostat} to find that for $i=0,1,2,3$
\be\label{eq:int.of.sum}
\begin{aligned}
 \sum_{(x,y)\in[0,1)^2\cap\left(\frac{1}{r}\BZ^2-(J,L)\right)}^{i\leq2x+2y\leq i+1}s_{\alpha,\beta}(x,y,r)
 & =
 \iint_{R_i} s_{\alpha,\beta}(x,y,r)\,{\rm d}x{\rm d}y
 +T_{\alpha,\beta}(J,L,r) ,
\end{aligned}
\ee
where $T$ is given by integrals over a subset of the boundary $\partial R_i$ and the imaginary axis above it with integrand $s_{\alpha,\beta}$ divided by ${\rm e}^{2\pi x}-1$ for various arguments $x$.

We will see that the asymptotics of $\int\!\!\int_R s_{-1,1}(x,y,r)\,{\rm d}x{\rm d}y$ can be determined by the saddle point method. However, first we need to see that the boundary makes a negligible contribution to the integral. In order to see this notice that at the vertices the imaginary part of $V_{\alpha,\beta}(x,y)$ vanishes, i.e., for all $x,y\in\{0,\tfrac{1}{2},1\}$ we have $\Im(V_{\alpha,\beta}(x,y))=0$. The initial integrals over~$\partial R_i$ have integrands that take values that are larger than $O(\exp(-\VC_{\rho_6} r/2\pi {\rm i}))$ however there is a~cancellation. See Figure~\ref{fig:oht.plot}.

To show that the integrals in a neighbourhood of the boundary are in fact exponentially small, we simply choose a parametrisation with endpoints in the regions where the $V_{-1,1}$ is exponentially small. The we deform the contours so they become exponentially small. See Lemma~\ref{lem:BCDsmall}.
We see that these integrals are $o(\exp(-\VC_{\rho_6} r/2\pi {\rm i}))$ as $r\rightarrow\infty$. We can apply similar analysis for the other boundaries and this can also be used to see that the boundary terms $T_{-1,1}(J,L,r)$ coming from Abel--Plana summation will all be $o(\exp(-\VC_{\rho_6} r/2\pi {\rm i}))$ as $r\rightarrow\infty$.

One can similarly show that by deformations of the contours
\[
 \iint_{R_i} s_{-1,1}(x,y,r)\,{\rm d}x{\rm d}y
  =o(\exp(-\VC_{\rho_6} r/2\pi {\rm i}))
 \qquad\text{for}\quad
 i =0,2,3 .
\]
Again see Lemma~\ref{lem:BCDsmall}.
Then finally we see that above region $R_{1}$ the integrand $s_{-1,1}$ has a critical point given by
$(x_6,y_6)
=(0.33675\ldots - 0.0027718\dots {\rm i},
0.55471\ldots + 0.014062\dots {\rm i})$
given by the sixth embedding of equation~\eqref{eq:xysp}.
We can deform the contour over this region by steepest descent.
This can be done everywhere until the contour is exponentially small or hits a critical point. The only critical point over this region is $(x_6,y_6)$.
Then applying stationary phase to the integral of these integrals we find that they have critical point at exactly $(x_6,y_6)$ and the resulting asymptotics are determined by a formal Gaussian integration around this critical point as claimed.
\end{proof}

\subsection[Quantum modularity of $\widehat{Z}$]{Quantum modularity of $\boldsymbol{\widehat{Z}}$}

We can describe a similar version of quantum modularity for $Z_{0,0}(q)$ from equation~\eqref{eq:Z00} as for the WRT invariant given in Theorem~\ref{thm:qmodold}. Here, instead of taking radial asymptotics, we take asymptotics that behave in a similar way to the WRT invariant. That is, we take asymptotics as~$\tau$ tends to $\infty$ horizontally. This has the affect that for all $n\in\BZ$ we have $q^{n}=o(\tau)$. Therefore, the asymptotics behaviour of $Z_{0,0}(\tq)$ could contain corrections with powers of $q$ that would have been exponentially small when $\tau$ went to ${\rm i}\infty$ on some angle or vertically. Indeed, this can be observed numerically. We have the following numerics:\footnote{We use the numerical Borel resummation here, however, one could replace this with a truncation or an optimal truncation.}
\begin{gather*}
 Z_{0,0}(\e(-1/(1000+{\rm i})))
  =
 (3.3121\ldots - 0.69001\dots {\rm i})\times 10^{97},\\
 \e(-1/8)\sqrt{\frac{1000+{\rm i}}{\delta_{\rho_6}}}\e\biggl(-\frac{\VC_{\rho_6}}{(2\pi {\rm i})^2}(1000+{\rm i})\biggr)s_{100}\bigl(\Phi^{(\rho_{6})}\bigr)\biggl(-\frac{2\pi {\rm i}}{1000+{\rm i}}\biggr)\\
  \qquad{}=(3.3183\ldots - 0.69130\dots {\rm i})\times 10^{97}.
\end{gather*}
Then we find that
\begin{gather*}
 Z_{0,0}(\e(-1/(1000+{\rm i})))\e(1/8) \sqrt{\frac{\delta_{\rho_6}}{1000+{\rm i}}} \e\biggl(\frac{\VC_{\rho_6}}{(2\pi {\rm i})^2}(1000+{\rm i})\biggr)s_{100}\bigl(\Phi^{(\rho_{6})}\bigr)\biggl(-\frac{2\pi {\rm i}}{1000+{\rm i}}\biggr)^{-1}\\
 \qquad{} =0.99813\ldots\\
 \qquad{} =1-0.0018674\ldots\\
 \qquad{} =1-\e(1000+{\rm i})+1.3025\ldots\times10^{-8}\\
 \qquad{} =1-\e(1000+{\rm i})+2\e(1000+{\rm i})^3-8.4823\ldots\times10^{-17}\\
 \qquad{} =1-\e(1000+{\rm i})+2\e(1000+{\rm i})^3-2\e(1000+{\rm i})^6+2.7775\ldots\times10^{-25}.
\end{gather*}
Continuing one can find the right-hand side is, up to exponentially small corrections, given by the $q$-series $Z_{0,0}(\e(1000+{\rm i}))$. We can apply a similar analysis as we did in the proof of Theorem~\ref{thm:qmodold} to obtain the following theorem.
\begin{Theorem}\label{thm:qmodold.zhat}
The $\widehat{Z}$--invariant computed in~{\rm \cite{GukovMan}} for $4_{1}(-1,2)$ is a quantum modular form in the sense of~{\rm \cite{Zagier:Qmod}}. That is, for $y\in\BR_{>0}$ and $x\in\BR$ as $x\rightarrow\infty$,
\[
 \frac{Z_{0,0}(\e(-1/(x+{\rm i}y)))}{Z_{0,0}(\e(x+{\rm i}y))}
 \sim
 \e(-1/8) \sqrt{\frac{x+{\rm i}y}{\delta_{\rho_6}}} \e\biggl(-\frac{\VC_{\rho_6}}{(2\pi {\rm i})^2}(x+{\rm i}y)\biggr)\Phi^{(\rho_{6})}\biggl(-\frac{2\pi {\rm i}}{x+{\rm i}y}\biggr) .
\]
\end{Theorem}
\begin{proof}[Outline of proof]
Here we describe the adjustments one needs to make to the analysis done in Section~\ref{statphase} to the case when $\tau$ is in the upper half plane as opposed to a rational number. The main difference here is that the sum now has infinitely many terms for finite $\tau$. Firstly, writing
\[
 Z_{0,0}(q)
  =
 \sum_{k,j=0}^{\infty}(-1)^{j}q^{k(2k+1)+jk+\frac{1}{2}j(j+1)}\frac{(q;q)_{k+j}}{(q;q)_{2k}(q;q)_{j}} ,
\]
we see that for $k,j>|\tau|^{1+\epsilon}$ for any fixed $\epsilon>0$ then the summand is $o(1)$ as $\e\bigl(-|\tau|^{2+2\epsilon}/\tau\bigr)=o(1)$ for $\Im(\tau)>0$. Therefore, We can replace the sum by a finite sum with the same asymptotics for fixed $\tau$ restricting $k,j<|\tau|^{1+\epsilon}$. Then one can apply the same arguments of Section~\ref{statphase} while using the same function $f$ from equation~\eqref{eq:f}, which is also defined for $\Im(\tau)>0$ and satisfies the important equation~\eqref{eq:f.qp}.
\end{proof}

While we saw in Section~\ref{radzhat} that the radial asymptotics of the $Z_{0,0}$ are given by the asymptotic series $\Phi^{(\rho_{3})}$, we can find the WRT invariant or in more detail the asymptotic series $\Phi^{(\rho_{0})}$ as an exponentially small correction. Indeed, we find that
\begin{gather*}
 \widehat{Z}(\e({\rm i}/1000))
  =3.7315\ldots\times10^8,\\
 \widehat{Z}(\e({\rm i}/1000))-\sqrt{100} \e(-1000 {\rm i})\exp\left(\VC_{\rho_3}100/2\pi\right)s_{200}(\Phi^{(\rho_3)})(2\pi/100)
  =0.011693\ldots,\\
 2s_{200}(\Phi^{(\rho_0)})(2\pi/100)
  =0.011693\ldots .
\end{gather*}
Therefore, while this example disagrees with the original conjecture on the radial asymptotics of~$Z_{0,0}$, the series $Z_{0,0}$ sees contributions from the trivial connection,\footnote{Note that there is an additional factor of $2$ appearing in front of the asymptotics of the trivial connection, which indicates that it would be natural to halve $Z_{0,0}$. Integrability is often wanted for the function $\widehat{Z}$ as it is hoped this will lead to a categorification however these factors of $2$ could be related to those found in~\cite{GPPV} related to $-1/2=\zeta(0)=1+1+1+\cdots$.} which can be exponentially small.

\subsection{Refined quantum modularity}\label{sec:RQMOD}

The quantum modularity conjecture~\cite{Zagier:Qmod} for hyperbolic knots was refined in~\cite{GZ:RQMOD}. That work showed how these examples of quantum modular forms could be viewed as giving rise to cocycles in some non-commutative cohomology of $\SL_{2}(\BZ)$. Indeed, for usual modular forms $f\colon \UHP\rightarrow\BC^{N}$ we have a transformation formula such that for $\gamma=[a,b;c,d]\in\SL_{2}(\BZ)$
\[
 f\biggl(\frac{a\tau+b}{c\tau+d}\biggr)
  =
 \rho_{\gamma} f(\tau) \mathbf{j}_{\gamma}(\tau) ,
\]
where $\rho\colon \SL_{2}(\BZ)\rightarrow\GL_{N}(\BC)$ a representation and $\mathbf{j}_{\gamma}(\tau)$ an automorphy factor. The representation $\rho$ can be viewed as an element of the non-commutative group cohomology where $\SL_{2}(\BZ)$ acts trivially on $\GL_{N}(\BC)$
\begin{gather*}
 H^{1}(\SL_{2}(\BZ),\GL_{N}(\BC))
  =
 \{\rho\colon \SL_{2}(\BZ)\rightarrow\GL_{N}(\BC)\mid
 \rho(\gamma_{1}\gamma_{2})
  =
 \rho(\gamma_{1})\rho(\gamma_{2})\}.
\end{gather*}
While examples of quantum invariants coming from hyperbolic manifolds will never give rise to functions $f$ with this structure they do give something very similar. In particular, Garoufalidis--Zagier noticed that quantum invariants of hyperbolic knots gave rise to $f$ and $\Omega$ such that
\[
 f\biggl(\frac{a\tau+b}{c\tau+d}\biggr)
  =
 \Om_{\gamma}(\tau) f(\tau) \mathbf{j}_{\gamma}(\tau) ,
\]
where $f$ is a function from the upper and lower half planes and the rationals to $\BC^{N}$, and for each $\gamma=[a,b;c,d]\in\SL_{2}(\BZ)$, we have $\Om_{\gamma}\colon \BC_{\gamma}\rightarrow\GL_{N}(\BC)$ is holomorphic where
$\BC_{\gamma}=\BC\smallsetminus \frac{1}{c}\BR_{\leq-d}$
if $c\neq0$ and $\BC_{\gamma}=\BC\smallsetminus \BR$ if $c=0$, and
\[
 \Om_{\gamma_1\gamma_2}(\tau)
  =
 \Om_{\gamma_1}(\gamma_2\cdot\tau)\Om_{\gamma_2}(\tau) ,
\]
where $\SL_{2}(\BZ)$ acts as usual via M\"obius transformations. We see that this $\Om$ now gives rise to an element of a more interesting twisted group cohomology valued in $N\times N$ matrices of analytic functions $\omega_N$ analytic on a simply connected dense domain
\[
 H^{1}(\SL_{2}(\BZ),\omega_N)
 :=
 \frac{\{\Om\colon \SL_{2}(\BZ)\rightarrow\omega_N \mid
 \Om_{\gamma_1\gamma_2}(\tau)
  =
 \Om_{\gamma_1}(\gamma_2\cdot\tau)\Om_{\gamma_2}(\tau)\}}
 {\{\Om_{\gamma}(\tau)\sim g(\gamma\cdot\tau)^{-1}\Omega_{\gamma}(\tau)g(\tau)\mid
 g\in\omega_{N}\}} .
\]
Often $\mathbf{j}_{\gamma}(\tau)$ can be written as $h(\tau)^{-1}h(\gamma\tau)$ and therefore it appears that $\Omega$ would often be trivial in the first cohomology group. It does become a coboundary when we restrict the domains of~$\Omega$, however, as defined it is not in general a coboundary. This is because $f$ is only defined for away from the irrational reals. This is the most important property of quantum modularity. Our examples will be functions that give rise to an analytic continuation to a simply connected dense domain of $\BC$. This is summarised in the following definition.
\begin{Definition}[matrix valued quantum modular form]
We say that a function $f\colon \UHP\cup\BQ\cup\overline{\UHP}\rightarrow\GL_{N}(\BC)$ is a matrix valued quantum modular form with automorphy factor $\mathbf{j}$ if for all $\gamma=[a,b;c,d]\in\SL_{2}(\BZ)$ the function
\be\label{eq:cocyc}
 \Om_{\gamma}(\tau):=f\biggl(\frac{a\tau+b}{c\tau+d}\biggr)\,\mathbf{j}_{\gamma}(\tau)^{-1}f(\tau)^{-1}
\ee
extends to an analytic function for $\tau\in\BC_{\gamma}$.
\end{Definition}
This definition has a variety of examples. The first interesting example came from the figure-eight knot.
\begin{Example}[figure-eight knot]\label{ex:41}
There is a matrix valued function associated to the figure-eight knot. The function at rationals appeared in the work of Garoufalidis--Zagier~\cite{GK:rat,GZ:RQMOD} and for the upper and lower half plane various other works~\cite{GGM:I,GGMW,GK:qser,GZ:qser}. This was also discussed in~\cite{Wh:thesis}. Consider the, functions for $|q|\neq 1$
\begin{gather*}
 g_{m}(q) =\sum_{k=0}^{\infty}(-1)^{k}\frac{q^{k(k+1)/2+km}}{(q;q)_{k}^2} ,\\
 G_{m}(q) =\sum_{k=0}^{\infty}(-1)^{k}\frac{q^{k(k+1)/2+km}}{(q;q)_{k}^2}\Biggl(m-2E_{1}(q)+\sum_{\ell=1}^{k}\frac{1+q^j}{1-q^j}\Biggr) ,\\
 \mathfrak{G}_{m}(q) =\sum_{k=0}^{\infty}(-1)^{k}\frac{q^{k(k+1)/2+km}}{(q;q)_{k}^2}\Bigg(\frac{1}{8}\Biggl(2m-4E_{1}(q)+2\sum_{\ell=1}^{k}\frac{1+q^j}{1-q^j}\Biggr)^2\!\!-\frac{1}{24}+ \! \sum_{\ell=1}^{k}\frac{q^{\ell}}{(1-q^{\ell})^2}\Bigg) .
\end{gather*}
Combining these into a matrix, we can define
\[
 \mathfrak{g}_{m}(q)
  =
 \begin{pmatrix}
 g_{m}(q) & G_{m}(q) & \mathfrak{G}_{m}(q)\\
 g_{m+1}(q) & G_{m+1}(q) & \mathfrak{G}_{m+1}(q)\\
 g_{m+2}(q) & G_{m+2}(q) & \mathfrak{G}_{m+2}(q)
 \end{pmatrix} .
\]
For $q=\e(N/M)$ with $N$ and $M$ coprime, and for $j=1,2$
\[
 X_{j}
  =\frac{1}{2}+(-1)^{j}\frac{\sqrt{-3}}{2}
\]
consider the functions
\[
 J_{j,m}(q)
  =
 \sum_{k\in\BZ/M\BZ}\frac{(-1)^{k}q^{k(k+1)/2+km}X_{j}^{(k-m)/M+1/2M}}{\prod_{\ell=0}^{M-1}\bigl(1-q^{1+k+\ell}X_{j}^{1/M}\bigr)^{2(1+\ell+k)/M-1}}
\]
and
\[
 J_{0,m}(q)
  =
 \sum_{k=0}^{\infty}(-1)^{k}q^{-k(k+1)/2-km-m}(q;q)_{k}^{2} .
\]
Combining these into a matrix we define
\[
 \mathfrak{g}_{m}(q)
  =
 \begin{pmatrix}
 J_{0,m}(q) & J_{1,m}(q) & J_{2,m}(q)\\
 J_{0,m+1}(q) & J_{1,m+1}(q) & J_{2,m+1}(q)\\
 J_{0,m+2}(q) & J_{1,m+2}(q) & J_{2,m+2}(q)
 \end{pmatrix} .
\]
Finally, take the automorphy factor for $\gamma=[a,b;c,d]\in\SL_{2}(\BZ)$
\[
 \mathbf{j}(\tau;\gamma)
  =
 \begin{cases}
 \begin{pmatrix}
 1 & 0 & 0\\
 0 & (c\tau+d) & 0\\
 0 & 0 & (c\tau+d)^{2}
 \end{pmatrix} & \text{if }\tau\in\UHP\cup\overline{\UHP} ,\\
 \begin{pmatrix}
 (c\tau+d)^{3/2} & 0 & 0\\
 0 & \e\bigl(\lambda_{\gamma}(\tau)\VC_{4_1}/(2\pi {\rm i})^2\bigr) & 0\\
 0 & 0 & \e\bigl(-\lambda_{\gamma}(\tau)\VC_{4_1}/(2\pi {\rm i})^2\bigr)
 \end{pmatrix} & \text{if }\tau\in\BQ,
 \end{cases}
\]
where $\VC_{4_1}=2{\rm i}\Im(\Li_2(\e(1/6)))=2.02988\dots {\rm i}$ is the complexified volume of the figure-eight knot and $\lambda_{\gamma}(r/s)=c/(s(cr+ds))$ is the cocycle defined in~\cite[equation (24)]{GZ:RQMOD}. Then $\mathfrak{g}$ is a matrix valued quantum modular form, a theorem proved in the previously cited works~\cite{GGM:I,GGMW,GK:rat,GK:qser,GZ:qser,GZ:RQMOD,Wh:thesis}.
\end{Example}
One of the main results of this current work is that the same holds for a matrix associated to half surgery on the figure-eight knot. Firstly, we define the following $q^{\pm1}$-series with $\bigl(q^{-1};q^{-1}\bigr)_{\infty}=(q;q)_{\infty}^{-1}$, $\bigl(x;q^{-1}\bigr)_{\infty}=(qx;q)_{\infty}^{-1}$ and $E_{k}\bigl(q^{-1}\bigr)=-E_{k}(q)$ defined in equation~\eqref{Eell}:
\begin{gather}
 Z^{(1)}_{m}(q)
   =
 (q;q)_{\infty}^2\sum_{k,j=0}^{\infty}\frac{q^{k(2k+1)+jk+j-mk-m}}{(q;q)_{j}(q;q)_{2k}(q;q)_{k+j}},
\nonumber\\
 Z^{(2)}_{m}(q)
   =
 (q;q)_{\infty}^2\nonumber\\ \hphantom{Z^{(2)}_{m}(q)  =}{}
 \times\sum_{k,j=0}^{\infty}
 \frac{q^{k(2k+1)+jk+j-mk-m}}{(q;q)_{j}(q;q)_{2k}(q;q)_{k+j}}\Biggl(-k-\frac{1}{2}+2E_{1}(q)-\sum_{n=1}^{j}\frac{q^n}{1-q^n}-\sum_{n=1}^{k+j}\frac{q^n}{1-q^n}\Biggr)\nonumber\\ \hphantom{Z^{(2)}_{m}(q)  =}{}
  +(q;q)_{\infty}^2\sum_{k=0}^{\infty}\sum_{j=-k}^{-1}
 \frac{q^{k(2k+1)+jk+j-mk-m}(q^{-1};q^{-1})_{-j-1}}{(q;q)_{2k}(q;q)_{k+j}},\nonumber\\
 Z^{(3)}_{m}(q)
   =
 (q;q)_{\infty}^2\sum_{k,j=0}^{\infty}
 \frac{q^{k(2k+1)+jk+j-mk-m}}{(q;q)_{j}(q;q)_{2k}(q;q)_{k+j}}\nonumber\\ \hphantom{ Z^{(3)}_{m}(q)  = (q;q)_{\infty}^2\sum_{k,j=0}^{\infty}}{}
  \times\Biggl(\frac{1}{2}\Biggl(7k+2j-2m-4E_1(q)+\sum_{n=j+1}^{k+j}\frac{q^n}{1-q^n}+4\sum_{n=1}^{2k}\frac{q^n}{1-q^n}\Biggr)\nonumber\\ \hphantom{ Z^{(3)}_{m}(q)  = (q;q)_{\infty}^2\sum_{k,j=0}^{\infty}}{}
  \times\Biggl(k+\frac{1}{2}-2E_1(q)+\sum_{n=1}^{k+j}\frac{q^n}{1-q^n}+\sum_{n=1}^{j}\frac{q^n}{1-q^n}\Biggr)\nonumber\\ \hphantom{ Z^{(3)}_{m}(q)  = (q;q)_{\infty}^2\sum_{k,j=0}^{\infty}\times\Biggl(}{}
  +2E_2(q)+\frac{1}{2}\sum_{n=1}^{k+j}\frac{q^n}{(1-q^n)^2}-\frac{1}{2}\sum_{n=1}^{j}\frac{q^n}{(1-q^n)^2}\Biggr)\nonumber\\ \hphantom{ Z^{(3)}_{m}(q)  = }{}
  +(q;q)_{\infty}^2\sum_{k=0}^{\infty}\sum_{j=-k}^{-1}
 \frac{q^{k(2k+1)+jk+j-mk-m}(q^{-1};q^{-1})_{-j-1}}{(q;q)_{2k}(q;q)_{k+j}}\nonumber\\ \hphantom{ Z^{(3)}_{m}(q)  = (q;q)_{\infty}^2\sum_{k=0}^{\infty}\sum_{j=-k}^{-1}}{}
  \quad\times\Biggl(3k+j-m-\frac{3}{4}-E_1(q)+2\sum_{n=1}^{2k}\frac{q^n}{1-q^n}+\sum_{n=1}^{-j-1}\frac{q^{-n}}{1-q^{-n}}\Biggr),\nonumber\\
 Z^{(4)}_{m}(q)
   =
 (q;q)_{\infty}\bigl(q^{3/2};q\bigr)_{\infty}\sum_{k,j=0}^{\infty}\frac{q^{(2k+1)(2k+2)/2+(k+1/2)j-m(k+1/2)+j-m}}{(q;q)_{j}(q;q)_{2k+1}\bigl(q^{\frac{3}{2}};q\bigr)_{k+j}},\nonumber\\
 Z^{(5)}_{m}(q)
   =
 (q;q)_{\infty}(-q;q)_{\infty}\sum_{k,j=0}^{\infty}(-1)^{j+m}\frac{q^{k(2k+1)+jk+j-mk-m}}{(q;q)_{j}(q;q)_{2k}(-q;q)_{k+j}},\nonumber\\
 Z^{(6)}_{m}(q)
   =
 (q;q)_{\infty}\bigl(-q^{3/2};q\bigr)_{\infty}\sum_{k,j=0}^{\infty}(-1)^{j+m}\frac{q^{(2k+1)(2k+2)/2+(k+1/2)j-m(k+1/2)+j-m}}{(q;q)_{j}(q;q)_{2k+1}\bigl(-q^{3/2};q\bigr)_{k+j}},\nonumber\\
 Z^{(7)}_{m}(q)
   =
 (q;q)_{\infty}\bigl(q^{1/2};q\bigr)_{\infty}\sum_{k,j=0}^{\infty}\frac{q^{(2k+1)(2k+2)/2+(j-k-1/2)(k+1/2)+(j-k-1/2)-m(k+1/2)-m}}{\bigl(q^{1/2};q\bigr)_{j-k}(q;q)_{2k+1}(q;q)_{j}},\nonumber\\
 Z^{(8)}_{m}(q)
   =
 (q;q)_{\infty}(-q;q)_{\infty}\sum_{k,j=0}^{\infty}(-1)^{j+m}\frac{q^{k(2k+1)+(j-k)k+(j-k)-mk-m}}{(-q;q)_{j-k}(q;q)_{2k}(q;q)_{j}}, \nonumber\\
 Z^{(9)}_{m}(q)
   =
 (q;q)_{\infty}\bigl(-q^{1/2};q\bigr)_{\infty}\nonumber\\ \hphantom{Z^{(9)}_{m}(q)=}{}
 \times\sum_{k,j=0}^{\infty}(-1)^{j+m}\frac{q^{(2k+1)(2k+2)/2+(j-k-1/2)(k+1/2)+(j-k-1/2)-m(k+1/2)-m}}{\bigl(-q^{1/2};q\bigr)_{j-k}(q;q)_{2k+1}(q;q)_{j}}.
\label{eq:theqser}
\end{gather}
With these $q$-series define
\[
 (\mathbf{Z}_m(q))_{i,j}=\begin{cases} Z_{m+i-1}^{(j+1-\delta_{j,1}-\delta_{j,2})}(q) & \text{if }|q|<1 ,\\
  Z_{m+i-1}^{(j+1-\delta_{j,1})}(q) & \text{if }|q|>1 ,\end{cases}
\]
where $\bigl(q^{-1};q^{-1}\bigr)_{\infty}=(q;q)_{\infty}^{-1}$ and $\bigl(x;q^{-1}\bigr)_{\infty}=(qx;q)^{-1}$ for $|q|<1$ and $x\notin q^{\BZ}$.
Let $\mathbf{j}_{\gamma}(\tau)$ be the automorphy factor generated by
\begin{gather*}
 \mathbf{j}_{T}(\tau)
  =
 \begin{pmatrix}
 1 & 0 & 0 & 0 & 0 & 0 & 0 & 0\\
 0 & 1 & 0 & 0 & 0 & 0 & 0 & 0\\
 0 & 0 & 0 & 0 & 1 & 0 & 0 & 0\\
 0 & 0 & 0 & 1 & 0 & 0 & 0 & 0\\
 0 & 0 & 1 & 0 & 0 & 0 & 0 & 0\\
 0 & 0 & 0 & 0 & 0 & 0 & 0 & 1\\
 0 & 0 & 0 & 0 & 0 & 0 & 1 & 0\\
 0 & 0 & 0 & 0 & 0 & 1 & 0 & 0\\
 \end{pmatrix},\\
 \mathbf{j}_{S}(\tau)
  =
 \tau^{-\chi_{|q|>1}}
 \begin{pmatrix}
 \tau^{\chi_{|q|>1}} & 0 & 0 & 0 & 0 & 0 & 0 & 0\\
 0 & \tau^{1-2\chi_{|q|>1}} & 0 & 0 & 0 & 0 & 0 & 0\\
 0 & 0 & 0 & 1 & 0 & 0 & 0 & 0\\
 0 & 0 & 1 & 0 & 0 & 0 & 0 & 0\\
 0 & 0 & 0 & 0 & 1 & 0 & 0 & 0\\
 0 & 0 & 0 & 0 & 0 & 0 & 1 & 0\\
 0 & 0 & 0 & 0 & 0 & 1 & 0 & 0\\
 0 & 0 & 0 & 0 & 0 & 0 & 0 & 1\\
 \end{pmatrix}
\end{gather*}
when $\tau\notin\BR$ and $\chi$ is the characteristic function. For $\tau\in\BQ$, consider the equations~\eqref{eq:Xequs} with solutions~\eqref{eq:Xsols}. Recalling $\Delta$ and $X$ from equations~\eqref{eq:Delta} and~\eqref{eq:Xsols}, we define functions at roots of unity $q=\e(N/M)$
\begin{gather}\label{eq:per.fun}
 \WRT^{(j)}_{m}(q) =\frac{-{\rm i}}{M\etaroots(q)}\sqrt{\frac{\bigl(1-X_{1,j}^2\bigr)^2(1-X_{2,j})\bigl(1-X_{1,j}^{-1}X_{2,j}\bigr)}{\Delta_{\rho_j}}}\\
 \times\!\!\sum_{k,\ell\in\BZ/M\BZ}
 \frac{q^{k^2+k\ell-mk+\ell}X_{1,j}^{\frac{2k+\ell-m}{M}}X_{2,j}^{\frac{k+1}{M}}  \prod_{i=0}^{M-1}\!\bigl(1-q^{i+1+\ell-k}X_{1,j}^{-1/M}X_{2,j}^{1/M}\bigr)^{-(i+1+\ell-k)/M-1/2}}{ \prod_{i=0}^{M-1}\!\bigl(1-q^{i+1+\ell}X_{2,j}^{1/M}\bigr)^{(i+1+\ell)/M-1/2}\prod_{i=0}^{M-1}\!\bigl(1-q^{i+1+2k}X_{1,j}^{2/M}\bigr)^{(i+1+2k)/M-1/2}} ,\nonumber
\end{gather}
which agree with the functions~\eqref{eq:phicont}. Finally, we take
\[
 \WRT^{(0)}_{m}(q)
  =
 \sum_{0\leq\ell\leq k}(-1)^{k}q^{-\frac{1}{2}k(k+1)+\ell(\ell+1)+mk}\frac{(q;q)_{2k+1}}{(q;q)_{\ell}(q;q)_{k-\ell}} .
\]
Then for roots of unity $q=\e(N/M)$, we define
\begin{gather*}
 (\mathbf{Z}_m(q))_{i,j} =\WRT_{m+i-1}^{(j-1)}(q) .
\end{gather*}
Let $\mathbf{j}_{\gamma}(\tau)$ be the automorphy factor given by
\[
 \begin{pmatrix}
 1 & 0 & \cdots & 0\\
 0 & (c\tau+d)^{-1/2}\e\bigl(-\lambda_{\gamma}(\tau)\VC_{\rho_1}/(2\pi {\rm i})^2\bigr) & \cdots & 0\\
 \vdots & \vdots & \ddots & \vdots\\
 0 & 0 & \cdots & (c\tau+d)^{-1/2}\e\bigl(-\lambda_{\gamma}(\tau)\VC_{\rho_7}/(2\pi {\rm i})^2\bigr)
 \end{pmatrix}
\]
when $\tau\in\BQ$.

\begin{Theorem}\label{thm:zisqmod}
The function $\mathbf{Z}_{m}$ is a matrix valued quantum modular form with automorphy factor $\mathbf{j}_{\gamma}(\tau)$.
\end{Theorem}

The proof of this theorem uses state integrals given in Section~\ref{sec:triv.statint} and their factorisation using Lemmas~\ref{lem:7x7stateint},~\ref{lem:vanishingqseries},~\ref{lem:quadrels7x7},~\ref{lem:8x8stateint},~\ref{lem:inhom.vanishingqseries}. Before, giving these details we can reflect on the pictures in Figure~\ref{fig:qmod.wrt.old} illustrating the old version of quantum modularity. We see that by refining the original quantum modularity to a matrix valued function, we now have functions which are analytic as opposed to smooth from the left and right of every rational however discontinuous at each rational. We can plot the cocycle $\Om_{S}$ of equation~\eqref{eq:cocyc} at reals and see the analytic structure explicitly in the following Figure~\ref{fig:cocycleplots.wrt}.
\begin{figure}[t]\centering

\caption{Plots of the real part of the first row of the cocycle associated to $4_{1}(-1,2)$. Each colour is the plot for an entry of this row in order \textcolor{red}{red}, \textcolor{blue}{blue}, \textcolor{green}{green}, \textcolor{purple}{purple}, \textcolor{orange}{orange}, \textcolor{yellow}{yellow}, \textcolor{pink}{pink}, \textcolor{magenta}{magenta}. Here we cut the plot off when it gets close to $0$ as there is an exponential singularity of the form $\exp(\VC_{\rho}/2\pi {\rm i}\tau)$ and the oscillations becomes large requiring many more data points to see the smooth behaviour.}
\label{fig:cocycleplots.wrt}
\end{figure}
This quantum modularity conjecture also gives us the answer to what the functions $X^{(\rho_{j})}(q)$ from the observation~\eqref{eq:phiqmod} should be. Indeed, we see they should equal to \smash{$\WRT^{(j)}_{0}(\e(x))=\Phi^{(\rho_j)}_{x,0}(0)$}.

\section{The state integral}

Since work of Garoufalidis--Kashaev~\cite{GK:rat,GK:qser}, we have known how to factorise certain ``untrapped'' state integrals at rationals and in the upper and lower half plane. Their work makes a clear connection between $q$-hypergeometric sums and state integrals. In particular, to get a state integral associated to a $q$-hypergeometric sum one needs only replace the sum by an integral and all the $q$-Pochhammer symbols by Faddeev's quantum dilogarithm. For example, for the figure-eight knot in Example~\ref{ex:41}, we have $q$-series
\[
 \sum_{k=0}^{\infty}(-1)^{k}\frac{q^{k(k+1)/2}}{(q;q)_{k}^{2}}
\]
which corresponds to a state integral
\[
 \int_{\BR+{{\rm i}0}}\e\bigl(-x^2/2\bigr)\Phi_{\mathsf{b}}(x)^2\,{\rm d}x .
\]
Here the correspondence is obtained by replacing the sum over $k$ by an integral over $x$, the $(-1)^{k}q^{k(k+1)}$ by $\e\bigl(-x^2/2\bigr)$, and $(q;q)_{k}^{-1}$ by $\Phi_{\mathsf{b}}(x)$.
Sometimes state integrals cannot be factorised by pushing a contour to infinity such as
\[
 \int_{\BR+{\rm i}0}\e\bigl(-x^2\bigr)\Phi_{\mathsf{b}}(x)\, {\rm d}x ,
\]
which would come from the sum $\sum_{k=0}^{\infty}q^{k(k+1)}(q;q)_{k}^{-1}$.
For the $q$-series, this amounts to the fact the series does not converge when $|q|>1$. This can be remedied by an ``untrapping'' procedure, which was found in the context of state integrals in unpublished work of Garoufalidis--Kashaev~\cite{Gar:untraptalk}. This was then found in the context of $q$-series in~\cite{Wh:thesis} and using both methods factorisation can be done for a much larger class of state integrals. This uses the $q$-binomial theorem and the Fourier transform of the Faddeev quantum dilogarithm.

After factorising the state integral, there is a final step to proving quantum modularity. This step involves expressing the factorisation as a quotient of two matrices. This was originally discovered in work of Garoufalidis--Zagier~\cite{GZ:qser,GZ:RQMOD} under the name quadratic relations. In terms of $q$-difference equations, this can be described in terms of a certain duality discussed in~\cite{GW:qmod}. Once the identity is found, this step can be done automatically using $q$-holonomic methods of~\cite{WilfZeil,Zeilberger} and checking boundary conditions.

\subsection[Untrapping the $q$-series]{Untrapping the $\boldsymbol{q}$-series}

For quantum modularity, we want functions defined for both $|q|<1$ and $|q|>1$. To do this, we can rewrite our functions in a way that makes a definition for $|q|\neq1$ manifest. While there is still choice, the quantum modularity will fix this modulo multiplication by certain modular functions. In particular, to prove quantum modularity we need to construct a state integral that factorises in terms of $\widehat{Z}_{m}(q)$ and similar functions. To construct such a state integral, we start with another identity for $Z_{m,n}$, which follows from the $q$-binomial theorem,
\begin{align*}
 &Z_{m,n}(q)
 =
 \sum_{0\leq k\leq\ell}(-1)^{k+\ell}q^{\frac{1}{2}3k(k+1)+\frac{1}{2}\ell(\ell+1)-(m+1)k-m-n\ell-n}\frac{(q;q)_{\ell}}{(q;q)_{2k}(q;q)_{\ell-k}}\\
 &\quad{}=
 (q;q)_{\infty}\sum_{\ell,k,j=0}^{\infty}(-1)^{\ell}q^{2k(k+1)+k\ell+\frac{1}{2}\ell(\ell+1)+jk+j\ell+j-(m+n+1)k-m-n\ell-n}\frac{1}{(q;q)_{j}(q;q)_{2k}(q;q)_{\ell}}\\
 &\quad{}=
 (q;q)_{\infty}^2\sum_{k,j=0}^{\infty}q^{2k(k+1)+jk+j-(m+n+1)k-m-n}\frac{1}{(q;q)_{j}(q;q)_{2k}(q;q)_{k+j-n}}\\
 &\quad{}=
 (q;q)_{\infty}^2\sum_{k,j=0}^{\infty}q^{k^2+jk+j-(m+n)k-m-n}\frac{1}{(q;q)_{j-n}(q;q)_{2k}(q;q)_{j-k}} .
\end{align*}
The form of this expression suggests considering the following state integral
\begin{gather}\label{eq:state.int}
 \calS_{m,m'}(\tau) =q^{-m}\tq^{m'}\int_{\BR+{\rm i}0}\int_{\BR+{\rm i}0} \Phi_{\mathsf{b}}(x_2+c_{\mathsf{b}})\Phi_{\mathsf{b}}(2x_1+c_{\mathsf{b}})\Phi_{\mathsf{b}}(x_2-x_1+c_{\mathsf{b}})\\
 \hphantom{\calS_{m,m'}(\tau) =q^{-m}\tq^{m'}\int_{\BR+{\rm i}0}\int}{}
 \times\e(-x_1^2-x_1x_2+ix_1(m\mathsf{b}+m'\mathsf{b}^{-1})-{\rm i}x_2(\mathsf{b}+\mathsf{b}^{-1}))\,{\rm d}x_1{\rm d}x_2 .\nonumber
\end{gather}
This state integral factorises into bilinear combinations of the following $q$-series and $\tq$ series. Indeed, it is ``untrapped'' in the sense that pushing the contour of integration to $({\rm i}\infty,{\rm i}\infty)$ the integral tends to zero.

\subsection{Factorisation and the non-abelian connections}

In a series of works~\cite{GGM:I,GK:qser,GZ:qser}, the asymptotic series of non-abelian connections associated to the figure-eight knot were considered. These works gave a $2\times 2$ matrix valued quantum modular form. This missed a contribution from the trivial connection, which was included in~\cite{GGMW,Wh:thesis} giving a~${3\times 3}$ matrix valued quantum modular form. This structure is paralleled for the invariants associated to the manifold $4_1(-1,2)$. Therefore, our first piece of quantum modularity gives rise to a $7\times 7$ matrix valued quantum modular form.
\begin{Lemma}\label{lem:7x7stateint}
 We have the following identity for $q=\e(\tau)$, $\tq=\e(-1/\tau)$:
 \begin{align*}
 \calS_{m,m'}(\tau)&{}=
 Z_{m}^{(1)}(q)Z_{m'}^{(2)}\bigl(\tq^{-1}\bigr)+\tau Z_{m}^{(2)}(q)Z_{m'}^{(1)}\bigl(\tq^{-1}\bigr)+Z_{m}^{(4)}(q)Z_{m'}^{(5)}\bigl(\tq^{-1}\bigr)+Z_{m}^{(5)}(q)Z_{m'}^{(4)}\bigl(\tq^{-1}\bigr)\\
 &{}\quad\!+Z_{m}^{(6)}(q)Z_{m'}^{(6)}\bigl(\tq^{-1}\bigr)+Z_{m}^{(7)}(q)Z_{m'}^{(8)}\bigl(\tq^{-1}\bigr)+Z_{m}^{(8)}(q)Z_{m'}^{(7)}\bigl(\tq^{-1}\bigr)+Z_{m}^{(9)}(q)Z_{m'}^{(9)}\bigl(\tq^{-1}\bigr) .
 \end{align*}
\end{Lemma}
This lemma can be proved using the techniques of~\cite{GK:qser} and some details are given in Appendix~\ref{sec:qserfac}. This factorisation can be simplified further. In particular, we have the following lemma.
\begin{Lemma}\label{lem:vanishingqseries}
For $|q|>1$, we have
\[
 Z^{(2)}_{m}(q)
  =0 ,
\]
and for $|q|\neq0$, excluding the case of $j=3$ and $|q|<1$, we have
\begin{gather}
 q^{2m+2}Z^{(j)}_{m}(q)+\bigl(q^{2m+4}+q^{m+1}+q^{m+2}\bigr)Z^{(j)}_{m+1}(q)\nonumber\\
 \qquad\quad{}+\bigl(-q^{2m+7}-q^{2m+5}-q^{2m+4}+q^{m+3}+1\bigr)Z^{(j)}_{m+2}(q)\nonumber\\
 \qquad\quad{}+\bigl(-q^{2m+9}-q^{2m+7}-q^{2m+6}+q^{m+3}-q^{m+5}-q^{m+4}-q^{m+3}\bigr)Z^{(j)}_{m+3}(q)\nonumber\\
 \qquad\quad{}+\bigl(q^{2m+10}+q^{2m+9}+q^{2m+7}+q^{m+4}-q^{m+6}\bigr)Z^{(j)}_{m+4}(q)\nonumber\\
 \qquad\quad{}+\bigl(q^{2m+12}+q^{2m+11}+q^{2m+9}-q^{m+6}+q^{m+6}\bigr)Z^{(j)}_{m+5}(q)\nonumber\\
 \qquad\quad{}+\bigl(-q^{2m+12}-q^{m+7}\bigr)Z^{(j)}_{m+6}(q)-q^{2m+14}Z^{(j)}_{m+7}(q)\nonumber\\
 \qquad{}=2\delta_{j,1}\chi_{|q|<1}(1-q)+2\delta_{j,3}\chi_{|q|>1}(1-q) .\label{eq:z.qdif}
\end{gather}
\end{Lemma}
Note that equation~\eqref{eq:z.qdif} is up to a factor of $2$ the same as the equation~\eqref{7thorderqdiffinhom}. The homogeneous order eight equation is given in equation~\eqref{8thorderqdiff} with $n=0$. For $Z^{(1)}(q)$ with $|q|<1$, this was utilised in the proof of Proposition~\ref{prop:zhat}.
The proof of this lemma uses standard $q$-holonomic methods and essentially reproduces the basis one would obtain from the Frobenius method.
Therefore, we see that the state integral satisfies a homogeneous $7$-th order equation in~$m$ and~$m'$. It is then natural to collect these series into a $7\times7$ Wronskian matrix, where for $i,j=1,\dots,7$, we have
\[
 \bigl(\mathbf{Z}_m^{(7\times 7)}(q)\bigr)_{i,j}=\begin{cases} Z_{m+i-1}^{(j+2-\delta_{j,1})}(q) & \text{if }|q|<1, \\
 Z_{m+i-1}^{(j+2-2\delta_{j,1})}(q) & \text{if }|q|>1.\end{cases}
\]
Then for the companion matrix{\setlength{\arraycolsep}{2.5pt}
\begin{align*}
 &A^{(7\times 7)}(x,q)=\\
 &
 \begin{pmatrix}
 0 & 1 & 0 & 0 & 0 & 0 & 0\\
 0 & 0 & 1 & 0 & 0 & 0 & 0\\
 0 & 0 & 0 & 1 & 0 & 0 & 0\\
 0 & 0 & 0 & 0 & 1 & 0 & 0\\
 0 & 0 & 0 & 0 & 0 & 1 & 0\\
 0 & 0 & 0 & 0 & 0 & 0 & 1\\
 \frac{1}{q^{12}} & \frac{q^3x + q + 1}{q^{13}x} & \frac{(-q^7 - q^5 - q^4)x^2 + q^3x + 1}{q^{14}x^2} & \frac{(-q^5 - q^3 - q^2)x - q - 1}{q^{10}x} & \frac{(q^6 + q^5 + q^3)x - q^2 + 1}{q^{10}x} & \frac{q^3 + q^2 + 1}{q^5} & \frac{-q^5x - 1}{q^7x}
 \end{pmatrix}
\end{align*}}%
we have
\[
 \mathbf{Z}_{m+1}^{(7\times 7)}(q)=A^{(7\times 7)}(q^m,q)\mathbf{Z}_m^{(7\times 7)}(q) .
\]
This module is self dual in the sense of~\cite{GW:qmod}. In particular, we have a gauge transformation
$P(x;q)$ whose columns respectively are given by
\begin{gather*}
 \begin{pmatrix}
 \frac{x^4q^{10}(2q^4 + 4q^3 + 6q^2 + 4q + 2) + 2x^3q^9 + x^2q^3(2q^6 + 4q^5 + 4q^4 + 6q^3 + 4q^2 + 4q + 2)+2)}{x^7q^{21}} \\[1mm]
 \frac{-2x^4q^{10} +x^2q^3(- 2q^4 - 4q^3 - 4q^2 - 4q - 2) - 2}{x^6q^{15}} \\[1mm]
 \frac{x^2q^3(2q^2 + 4q + 2) + 2}{x^5q^{10}} \\[1mm]
 \frac{-2x^2q^3 - 2}{x^4q^6} \\[1mm]
 \frac{2}{x^3q^3} \\[1mm]
 \frac{-2}{x^2q} \\
 0
 \end{pmatrix} ,\\
 \begin{pmatrix}
 \frac{-2x^4q^{14} + x^2q^5(- 2q^4 - 4q^3 - 4q^2 - 4q - 2) - 2}{x^6q^{21}} \\[1mm]
 \frac{x^2q^5(2q^2 + 4q + 2) + 2}{x^5q^{15}} \\[1mm]
 \frac{-2x^2q^{5} - 2}{x^4q^{10}} \\[1mm]
 \frac{2}{x^3q^6} \\[1mm]
 \frac{-2}{x^2q^3} \\
 0 \\
 \frac{-2}{x^2q}
 \end{pmatrix} ,\qquad
 \begin{pmatrix}
 \frac{(x^2q^7(2q^2 + 4q + 2) + 2}{x^5q^{20}} \\[1mm]
 \frac{-2x^2q^7 - 2}{x^4q^{14}} \\[1mm]
 \frac{2}{x^3q^9} \\[1mm]
 \frac{-2}{x^2q^5} \\
 0 \\
 \frac{-2}{x^2q^3} \\[1mm]
 \frac{2}{x^3q^3}
 \end{pmatrix} ,\qquad
 \begin{pmatrix}
 \frac{-2x^2q^9 - 2}{x^4q^{18}} \\[1mm]
 \frac{2}{x^3q^{12}} \\[1mm]
 \frac{-2}{x^2q^7} \\
 0 \\
 \frac{-2}{x^2q^5} \\[1mm]
 \frac{2}{x^3q^6} \\[1mm]
 \frac{-2x^2q^3 - 2}{x^4q^6}
 \end{pmatrix} ,\\
 \begin{pmatrix}
 \frac{2}{x^3q^{15}} \\[1mm]
 \frac{-2}{x^2q^9} \\
 0 \\
 \frac{-2}{x^2q^7} \\[1mm]
 \frac{2}{x^3q^9} \\[1mm]
 \frac{-2x^2q^5 - 2}{x^4q^{10}} \\[1mm]
 \frac{x^2q^3(2q^2 + 4q + 2) + 2xq^2 + 2}{x^5q^{10}}
 \end{pmatrix} ,\qquad
 \begin{pmatrix}
 \frac{-2}{x^2q^{11}} \\
 0 \\
 \frac{-2}{x^2q^9} \\[1mm]
 \frac{2}{x^3q^{12}} \\[1mm]
 \frac{-2x^2q^7 - 2}{x^4q^{14}} \\[1mm]
 \frac{x^2q^5(2q^2 + 4q + 2) + 2xq^3 + 2}{x^5q^{15}} \\[1mm]
 \frac{-2x^4q^{10} +x^2q^3(- 2q^4 - 4q^3 - 4q^2 - 4q - 2) +xq^2(- 2q - 2) - 2}{x^6q^{15}}
 \end{pmatrix} ,\\
 \begin{pmatrix}
 0 \\
 \frac{-2}{x^2q^{11}} \\[1mm]
 \frac{2}{x^3q^{15}} \\[1mm]
 \frac{-2x^2q^9 - 2}{x^4q^{18}} \\[1mm]
 \frac{x^2q^7(2q^2 + 4q + 2) + 2xq^4 + 2}{x^5q^{20}} \\[1mm]
 \frac{-2x^4q^{14} + x^2q^5(- 2q^4 - 4q^3 - 4q^2 - 4q - 2) +xq^3(- 2q - 2) - 2}{x^6q^{21}} \\
 \frac{\text{\scriptsize $\begin{array}{c} x^4q^{10}(2q^4 + 4q^3 + 6q^2 + 4q + 2) + x^3q^7(2q^4 + 2q^3 + 2q^2 + 2q + 2) \\ + x^2q^3(2q^6 + 4q^5 + 4q^4 + 6q^3 + 4q^2 + 4q + 2) + xq^2(2q^2 + 2q + 2) + 2\end{array}$}}{x^7q^{21}}
 \end{pmatrix} ,
\end{gather*}
which satisfies the equation
\[
 A^{(7\times 7)}(x,q)
  =
 P(qx,q)A^{(7\times 7)}\bigl(q^7x,q^{-1}\bigr)^{\mathsf{T}}P(x,q)^{-1} .
\]
This duality can be used to give rise to quadratic relations between the $q$-series that is most conveniently given in matrix form.
\begin{Lemma}[quadratic relations]\label{lem:quadrels7x7}
We have the following identity:
\[
 \mathbf{Z}_{m}^{(7\times 7)}(q)\mathbf{Z}_{-m-6}^{(7\times 7)}\bigl(q^{-1}\bigr)^{\mathsf{T}}
  =
 P(q^{m},q)
 \qquad\text{or}\qquad
 \mathbf{Z}_{m}^{(7\times 7)}(q)^{-1}
  =
 \mathbf{Z}_{-m-6}^{(7\times 7)}\bigl(q^{-1}\bigr)^{\mathsf{T}}P(q^{m},q)^{-1} .
\]
\end{Lemma}

Combining the previous Lemmas~\ref{lem:7x7stateint},~\ref{lem:vanishingqseries}~and~\ref{lem:quadrels7x7}, we have the following theorem.
\begin{Theorem}\label{thm:7x7qser}
$\mathbf{Z}_{m}^{(7\times 7)}$ is a matrix valued quantum modular form. In particular,
\begin{alignat*}{3}
 &\mathbf{Z}_{m}^{(7\times 7)}\bigl|_{T}
  =
 \Om_{T}(\tau)\mathbf{Z}_{m}^{(7\times 7)}\ms_{T}(\tau;T) ,\quad
 &&\mathbf{Z}_{m}^{(7\times 7)}\bigl|_{S}
  =
 \Om_{S}(\tau)
 \mathbf{Z}_{6-m'}^{(7\times 7)}
 \ms_{S}(\tau;T) ,&\\
 &\Om_{T}(\tau) =\mathrm{I} ,\quad
 &&\Om_{S}(\tau) =\calS_{m,m'}P\bigl(q^{6-m'},q\bigr)^{-1} ,&\\
 &\ms_{T}(\tau;T)
  =
 \begin{pmatrix}
 1 & 0 & 0 & 0 & 0 & 0 & 0\\
 0 & 0 & 0 & 1 & 0 & 0 & 0\\
 0 & 0 & 1 & 0 & 0 & 0 & 0\\
 0 & 1 & 0 & 0 & 0 & 0 & 0\\
 0 & 0 & 0 & 0 & 0 & 0 & 1\\
 0 & 0 & 0 & 0 & 0 & 1 & 0\\
 0 & 0 & 0 & 0 & 1 & 0 & 0
 \end{pmatrix} ,\qquad
 &&\ms_{S}(\tau;S)
  =
 \begin{pmatrix}
 \tau & 0 & 0 & 0 & 0 & 0 & 0\\
 0 & 0 & 1 & 0 & 0 & 0 & 0\\
 0 & 1 & 0 & 0 & 0 & 0 & 0\\
 0 & 0 & 0 & 1 & 0 & 0 & 0\\
 0 & 0 & 0 & 0 & 0 & 1 & 0\\
 0 & 0 & 0 & 0 & 1 & 0 & 0\\
 0 & 0 & 0 & 0 & 0 & 0 & 1
 \end{pmatrix}, &
\end{alignat*}
where $\Om_S(\tau)=\calS_{m,m'}P\bigl(q^{6-m'},q\bigr)^{-1}$ is holomorphic on $\BC\smallsetminus\BR_{\leq0}$.
\end{Theorem}
\begin{proof}
Lemma~\ref{lem:7x7stateint} expresses the state integral as a bilinear combination of $q$ and $\tq$ series. This combination is therefore analytic for $\tau\in\BC\smallsetminus \BR_{\leq0}$ as the state integrals are manifestly analytic on this domain. This follows from the domain of Faddeev's quantum dilogarithm. Lemma~\ref{lem:vanishingqseries} removes one term from the bilinear combination. Then the remaining terms in the bilinear combination can be expressed as an entry of a quotient of matrices, which follows from Lemma~\ref{lem:quadrels7x7}.
\end{proof}

\subsection{Including the trivial connection}\label{sec:triv.statint}

As was previously mentioned, it was shown for the figure-eight knot~\cite{GGMW,Wh:thesis} there is a natural extension of a $2\times 2$ matrix valued quantum modular form to a $3\times 3$ matrix valued quantum modular form. In the example of $4_1(-1,2)$, there is a completely analogous statement. Here we see that for $|q|<1$ the $\widehat{Z}$ invariant appears as a component in the factorisation. This agrees with an aspect of previous conjectures~\cite{GukovMan} on the relation between the WRT invariant and the~$\widehat{Z}$ invariant.

As was illustrated in~\cite{GGMW}, Appell--Lerch type sums appear and we have a state integral with a~simple additional factor of the inverse of a hyperbolic function. Here however, the computation is a little more complicated as the inverse of $\mathbf{Z}_m(q)$ is more complicated than those cases. This comes from the fact that for the inhomogeneous equations satisfied by
\[
 \mathbf{Z}^{\mathrm{inhom}}_{m}(q)
  =
 \begin{pmatrix}
 1 & 0 & \cdots & 0\\
 Z^{(1+2\chi_{|q|>1})}_{m}(q) & \bigl(\mathbf{Z}_{m}^{(7\times 7)}\bigr)_{1,1} & \cdots & \bigl(\mathbf{Z}_{m}^{(7\times 7)}\bigr)_{1,7}\\
 \vdots & \vdots & \cdots & \vdots\\
 Z^{(1+2\chi_{|q|>1})}_{m}(q) & \bigl(\mathbf{Z}_{m}^{(7\times 7)}\bigr)_{7,1} & \cdots & \bigl(\mathbf{Z}_{m}^{(7\times 7)}\bigr)_{7,7}
 \end{pmatrix}
\]
is given by
\[
 \mathbf{Z}^{\mathrm{inhom}}_{m+1}(q)
  =
 A(x;q)\mathbf{Z}^{\mathrm{inhom}}_{m}(q)
\]
with
\[
 A(x;q)^{-1}
  =\begin{pmatrix}
 1&\frac{-2q + 2}{q^2x^2}&0&0&0&0&0&0\\[1mm]
 0&\frac{-q^3x + (-q - 1)}{qx}&1&0&0&0&0&0\\[1mm]
 0&\frac{(q^7 + q^5 + q^4)x^2 - q^3x - 1}{q^2x^2}&0&1&0&0&0&0\\[1mm]
 0&\frac{(q^7 + q^5 + q^4)x + (q^3 + q^2)}{x}&0&0&1&0&0&0\\[1mm]
 0&\frac{(-q^8 - q^7 - q^5)x + (q^4 - q^2)}{x}&0&0&0&1&0&0\\[1mm]
 0&-q^{10} - q^9 - q^7&0&0&0&0&1&0\\[1mm]
 0&\frac{q^{10}x + q^5}{x}&0&0&0&0&0&1\\[1mm]
 0&q^{12}&0&0&0&0&0&0
\end{pmatrix}^{\mathsf{T}}.
\]
Therefore,
\[
 \mathbf{Z}^{\mathrm{inhom}}_{m+1}(q)^{-1}
  =
 \mathbf{Z}^{\mathrm{inhom}}_{m}(q)^{-1}A(x;q)^{-1} .
\]
Let $L^{(j)}$ be given by
\begin{align*}
 &\mathbf{Z}^{\mathrm{inhom}}_{m}(q)^{-1} =\\
 &
 \begin{pmatrix}
 1 & 0 & \cdots & 0\\
 L^{(1+\chi_{|q|<1})}_{m}\bigl(q^{-1}\bigr) & \bigl(\mathbf{Z}_{-m-6}^{(7\times 7)}\bigl(q^{-1}\bigr)^{\mathsf{T}}P\bigl(q^{m},q\bigr)^{-1}\bigr)_{1,1} & \cdots & \bigl(\mathbf{Z}_{-m-6}^{(7\times 7)}\bigl(q^{-1}\bigr)^{\mathsf{T}}P\bigl(q^{m},q\bigr)^{-1}\bigr)_{1,7}\\[1mm]
 L^{(4)}_{m}\bigl(q^{-1}\bigr) & \bigl(\mathbf{Z}_{-m-6}^{(7\times 7)}\bigl(q^{-1}\bigr)^{\mathsf{T}}P(q^{m},q)^{-1}\bigr)_{2,1} & \cdots & \bigl(\mathbf{Z}_{-m-6}^{(7\times 7)}\bigl(q^{-1}\bigr)^{\mathsf{T}}P\bigl(q^{m},q\bigr)^{-1}\bigr)_{2,7}\\
 \vdots & \vdots & \cdots & \vdots\\
 L^{(9)}_{m}\bigl(q^{-1}\bigr) & \bigl(\mathbf{Z}_{-m-6}^{(7\times 7)}\bigl(q^{-1}\bigr)^{\mathsf{T}}P\bigl(q^{m},q\bigr)^{-1}\bigr)_{7,1} & \cdots & \bigl(\mathbf{Z}_{-m-6}^{(7\times 7)}\bigl(q^{-1}\bigr)^{\mathsf{T}}P\bigl(q^{m},q\bigr)^{-1}\bigr)_{7,7}
 \end{pmatrix} \!.
\end{align*}
The matrix $P(x,q)^{-1}$ has first column
\[
 \begin{pmatrix}
 q^{15}x^5\\[1mm]
 q^{13}x^5 + \bigl(\frac{1}{2}q^{10} + q^9\bigr)x^4\\[1mm]
 \bigl(-q^{13} - q^{12}\bigr)x^5 + q^8x^4 + \frac{1}{2}q^5x^3\\[1mm]
 \bigl(-q^{11} - q^{10}\bigr)x^5 + \frac{1}{2}q^8x^4\\[1mm]
 q^{10}x^5 + q^7x^4\\[1mm]
 q^8x^5\\
 0
 \end{pmatrix}.
\]
Then we have
\begin{align*}
 L^{(j)}_{m+1}\bigl(q^{-1}\bigr)
 ={}&
 L^{(j)}_{m}\bigl(q^{-1}\bigr)
 +\frac{2-2q}{q^{2m+2}}\bigl(\mathbf{Z}_{-m-6}^{(7\times 7)}\bigl(q^{-1}\bigr)^{\mathsf{T}}P\bigl(q^{m},q\bigr)^{-1}\bigr)_{1,1}\\
= {} &
 L^{(j)}_{m}\bigl(q^{-1}\bigr)
 +
 q^{13+3m}(2-2q)Z_{-m-6}^{(j)}\bigl(q^{-1}\bigr)\\
&
 +\biggl(q^{11+3m} + \biggl(\frac{1}{2}q^{8} + q^7\biggr)q^{2m}\biggr)(2-2q)Z_{-m-5}^{(j)}\bigl(q^{-1}\bigr)\\
 & +\biggl(\bigl(-q^{11} - q^{10}\bigr)q^{3m} + q^{6+2m} + \frac{1}{2}q^{3+m}\biggr)(2-2q)Z_{-m-4}^{(j)}\bigl(q^{-1}\bigr)\\
 & +\biggl(\bigl(-q^{9} - q^{8}\bigr)q^{3m} + \frac{1}{2}q^{6+2m}\biggr)(2-2q)Z_{-m-3}^{(j)}\bigl(q^{-1}\bigr)\\
 & +\bigl(q^{8+3m} + q^{5+2m}\bigr)(2-2q)Z_{-m-2}^{(j)}\bigl(q^{-1}\bigr)\\
 & +q^{6+3m}(2-2q)Z_{-m-1}^{(j)}\bigl(q^{-1}\bigr).
\end{align*}
The solution to such an equation is unique up to a constant in $m$.
For example, if we let
\[
 K^{(4)}_{m,n}(q)
  =
 (q;q)_{\infty}(q^{3/2};q)_{\infty}\sum_{k,j=0}^{\infty}\frac{q^{(2k+1)(2k+2)/2+(k+1/2)j-m(k+1/2)+j+(n-1)m}}{(q;q)_{j}(q;q)_{2k+1}(q^{3/2};q)_{k+j}(1-q^{3/2-n+k})} ,
\]
then
\[
 K^{(4)}_{m,n}(q)-K^{(4)}_{m-1,n}(q)
  =
 q^{m\,n}Z^{(4)}_{m}(q) .
\]
Then using this function, we can compute that
\begin{align}
 \frac{1}{2q(1-q)}L^{(4)}_{m}(q)
 ={}&q^{3}K^{(4)}_{-m-6,3}(q)
 +q^{2}K^{(4)}_{-m-5,3}(q) + \left(\frac{1}{2} + q\right)K^{(4)}_{-m-5,2}(q)\nonumber\\
 & +\bigl(-q^{-1} - 1\bigr)K^{(4)}_{-m-4,3}(q) + K^{(4)}_{-m-4,2}(q) + \frac{1}{2}q^{-1}K^{(4)}_{-m-4,1}(q)\nonumber\\
 & +\bigl(-q^{-2} - q^{-1}\bigr)K^{(4)}_{-m-3,3}(q) + \frac{1}{2}q^{-2}K^{(4)}_{-m-3,2}(q)\nonumber\\
 & +q^{-4}K^{(4)}_{-m-2,3}(q) + q^{-3}K^{(4)}_{-m-2,2}(q)
 +q^{-5}K^{(4)}_{-m-1,3}(q).\label{eq:LinK}
\end{align}
This can be done for every series in equation~\eqref{eq:theqser} completely analogously and more details are given in equations~\eqref{eq:theqser.ler}--\eqref{eq:LinK.gen}.
Taking this into account, we take the state integral
\begin{align}
 \calS^{(0)}_{m,m'}(\tau) ={}&\frac{2-2\tq}{q^{m}\tq^{2-m'}}\int_{\BR-{\rm i}0}\int_{\BR-{\rm i}0}\Phi_{\mathsf{b}}(x_2+c_{\mathsf{b}}) \Phi_{\mathsf{b}}(2x_1+c_{\mathsf{b}})\Phi_{\mathsf{b}}(x_2-x_1+c_{\mathsf{b}})\nonumber\\
 & \times\e\bigl(-x_1^2-x_1x_2+{\rm i}x_1\bigl(m\mathsf{b}-m'\mathsf{b}^{-1}\bigr)-{\rm i}x_2\bigl(\mathsf{b}+\mathsf{b}^{-1}\bigr)\bigr)\nonumber\\
 &\times\biggl(\tq^{-3}\frac{\e\bigl(-6{\rm i}x_1\mathsf{b}^{-1}\bigr)}{1-\e\bigl(-{\rm i}x_1\mathsf{b}^{-1}\bigr)\tq^{2}}
 +\tq^{-2}\frac{\e\bigl(-5{\rm i}x_1\mathsf{b}^{-1}\bigr)}{1-\e\bigl(-{\rm i}x_1\mathsf{b}^{-1}\bigr)\tq^{2}} \nonumber\\
& + \biggl(\frac{1}{2} + \tq^{-1}\biggr)\frac{\e\bigl(-5{\rm i}x_1\mathsf{b}^{-1}\bigr)}{1-\e\bigl(-{\rm i}x_1\mathsf{b}^{-1}\bigr)\tq}
  +(-\tq - 1)\frac{\e\bigl(-4{\rm i}x_1\mathsf{b}^{-1}\bigr)}{1-\e\bigl(-{\rm i}x_1\mathsf{b}^{-1}\bigr)\tq^{2}}\nonumber\\
   & + \frac{\e\bigl(-4{\rm i}x_1\mathsf{b}^{-1}\bigr)}{1-\e\bigl(-{\rm i}x_1\mathsf{b}^{-1}\bigr)\tq} + \frac{1}{2}\tq \frac{\e\bigl(-4{\rm i}x_1\mathsf{b}^{-1}\bigr)}{1-\e\bigl(-{\rm i}x_1\mathsf{b}^{-1}\bigr)}
  +\bigl(-\tq^2 - \tq\bigr)\frac{\e\bigl(-3{\rm i}x_1\mathsf{b}^{-1}\bigr)}{1-\e\bigl(-{\rm i}x_1\mathsf{b}^{-1}\bigr)\tq^{2}} \nonumber\\
  & + \frac{1}{2}\tq^2\frac{\e\bigl(-3{\rm i}x_1\mathsf{b}^{-1}\bigr)}{1-\e\bigl(-{\rm i}x_1\mathsf{b}^{-1}\bigr)\tq}
  +\tq^4\frac{\e\bigl(-2{\rm i}x_1\mathsf{b}^{-1}\bigr)}{1-\e\bigl(-{\rm i}x_1\mathsf{b}^{-1}\bigr)\tq^{2}} + \tq^3\frac{\e\bigl(-2{\rm i}x_1\mathsf{b}^{-1}\bigr)}{1-\e\bigl(-{\rm i}x_1\mathsf{b}^{-1}\bigr)\tq}\nonumber\\
  &
 +\tq^5\frac{\e\bigl(-{\rm i}x_1\mathsf{b}^{-1}\bigr)}{1-\e\bigl(-{\rm i}x_1\mathsf{b}^{-1}\bigr)\tq^{2}}\biggr)\,{\rm d}x_1{\rm d}x_2 .\label{eq:inhom.statint}
\end{align}
The following lemma can be proved using the techniques of~\cite{GK:qser} and some details are given in Appendix~\ref{sec:qserfac}. We take
\[
 K^{(0)}_{m,n}(q) =(q;q)_{\infty}^2\sum_{j=0}^{\infty}\frac{q^{(n-1)(2n-1)+nj}}{(q;q)_{j}(q;q)_{2n-2}(q;q)_{n-1+j}}
\]
and define $L^{(0)}_{m}(q)$ as in equation~\eqref{eq:LinK}.
\begin{Lemma}\label{lem:8x8stateint}
 We have the following identity:
 \begin{align*}
\calS^{(0)}_{m,m'}(\tau)
  ={}&
 Z_{m}^{(1)}(q)L_{m'}^{(2)}\bigl(\tq^{-1}\bigr)+\tau Z_{m}^{(2)}(q)L_{m'}^{(1)}\bigl(\tq^{-1}\bigr)+\tau^2 Z_{m}^{(3)}(q)L_{m'}^{(0)}\bigl(\tq^{-1}\bigr)\\
 &+Z_{m}^{(4)}(q)L_{m'}^{(5)}\bigl(\tq^{-1}\bigr)+Z_{m}^{(5)}(q)L_{m'}^{(4)}\bigl(\tq^{-1}\bigr)+Z_{m}^{(6)}(q)L_{m'}^{(6)}\bigl(\tq^{-1}\bigr)
\\
 & +Z_{m}^{(7)}(q)L_{m'}^{(8)}\bigl(\tq^{-1}\bigr) +Z_{m}^{(8)}(q)L_{m'}^{(7)}\bigl(\tq^{-1}\bigr)+Z_{m}^{(9)}(q)L_{m'}^{(9)}\bigl(\tq^{-1}\bigr) .
 \end{align*}
\end{Lemma}
This can be simplified further using the following identities again following from simple $q$-holonomic methods.
\begin{Lemma}\label{lem:inhom.vanishingqseries}
For $|q|\neq1$,
\[
 L^{(0)}_{m}(q)
  =\chi_{|q|<1} ,
\]
and for $|q|>1$
\[
 L^{(2)}_m(q) =1 .
\]
\end{Lemma}
Using the Lemmas~\ref{lem:7x7stateint},~\ref{lem:vanishingqseries},~\ref{lem:quadrels7x7},~\ref{lem:8x8stateint},~\ref{lem:inhom.vanishingqseries}, and an analogous description of $L^{(0)}_{m}$ when $q$ is a root of unity gives the proof of Theorem~\ref{thm:zisqmod}.

\section{Resurgence}
\label{res}

The study of resurgence in relation to quantum invariants of three-manifolds was initiated in~\cite{Gar:resCS}. Recently, a proposal for a complete description of the resurgent structure was given in~\cite{GGM:I,GGMW}. That work gives a conjectural, but complete, understanding of both the Borel resummation and the Stokes constants. The first being given by the state integrals of Andersen--Kashaev and the second by the 3d index of Dimofte--Gaiotto--Gukov~\cite{DGG}. For the case of closed manifolds neither of these invariants are known but there are conjectural proposals~\cite{Gang3d}. Regardless, the methods in~\cite{GGM:I,GGMW} can be applied and one can construct a state integral that conjecturally gives the Borel resummation of the asymptotic series $\Phi^{(\rho)}$ and give $q$-series as generating functions of Strokes constants. Both of these invariants are given as a product of two pieces. These are all built from a matrix $\mathbf{Z}$. The Borel resummation and Stokes constants respectively are conjectured~\cite{GGM:I,GPPV} to be roughly of the shape
``$\mathbf{Z}(\tq)\mathbf{Z}(q)^{-1}$''  and ``$\mathbf{Z}(q)\mathbf{Z}\bigl(q^{-1}\bigr)^{\mathsf{T}}$''.
In the example of $4_1(-1,2)$, we illustrate how these conjectures should be given precisely.

\subsection{Borel resummation}
\label{canbas}

We can find the quantum modular behaviour for $\tau\in\BR_{>0}$ and $\tau\in\BR_{<0}$ numerically by using Borel--Pad\'e--Laplace. This can be done with $q$-series, functions at roots of unity, or $q^{-1}$ series. There are no Stokes lines on the reals, or equivalently they have zero Stokes constants there (as can be numerically seen in the bottom right $2\times2$ block of $M$ in equation~\eqref{eq:closeststokes}), so we only need to do this once in this example. Numerical computations lead to the functions $X_{\rho}$ of equation~\eqref{eq:can.basis}. We can write this in terms of the module associated to equation~\eqref{8thorderqdiff} with $n=0$. This is stored in the matrix $P(q)$ of equations~\eqref{eq:can.basis.P} and~\eqref{eq:can.basis.Q}. Assuming the conjectures~\cite{Gar:resCS} that the asymptotic series $\Phi^{(\rho)}_{m}$ are Borel resummable, we let $s\bigl(\widehat{\Phi}\bigr)$ represent the Borel resummation of the $8\times 8$ matrix of asymptotic series indexed by $m=0,\dots,7$ by $\rho_0,\dots,\rho_7$. Considering the eight by eight matrix of $q$ series, we find numerically the conjectural identities.
\begin{Conjecture}\label{conj:borel}\setlength{\arraycolsep}{4pt}
For $\tau$ with small just above the positive reals $($region {\rm I}$)$,
\[
\mathbf{Z}_{0}(\tq) =s\bigl(\Phih\bigr)(-1/\tau)S_{\rm I}(q),
\]
where
\begin{align*}
S_{\rm I}(q) ={}&\begin{pmatrix}
 -1 & 0 & 0 & 0 & 0 & 0 & 0 & 0\\
 0 & -1 & 0 & 0 & 0 & 0 & 0 & 0\\
 0 & 0 & -1 & 0 & 0 & 0 & 0 & 0\\
 0 & 0 & 0 & -1 & 0 & 0 & 0 & 0\\
 0 & 0 & 0 & 0 & -1 & 0 & 0 & 0\\
 0 & 0 & 0 & 0 & 0 & -1 & 0 & 0\\
 0 & 0 & 0 & 0 & 0 & 0 & 1 & 0\\
 0 & 0 & 0 & 0 & 0 & 0 & 0 & 1
\end{pmatrix}P(q)\mathbf{Z}_{-3}(q)\\
&\times
\begin{pmatrix}
 -1 & 0 & 0 & 0 & 0 & 0 & 0 & 0\\
 0 & -1 & 0 & 0 & 0 & 0 & 0 & 0\\
 0 & 0 & 1 & 0 & 0 & 0 & 0 & 0\\
 0 & 0 & 0 & 1 & 0 & 0 & 0 & 0\\
 0 & 0 & 0 & 0 & -1 & 0 & 0 & 0\\
 0 & 0 & 0 & 0 & 0 & -1 & 0 & 0\\
 0 & 0 & 0 & 0 & 0 & 0 & -1 & 0\\
 0 & 0 & 0 & 0 & 0 & 0 & 0 & 1
\end{pmatrix}
\begin{pmatrix}
 1 & 0 & 0 & 0 & 0 & 0 & 0 & 0\\
 0 & \tau & 0 & 0 & 0 & 0 & 0 & 0\\
 0 & 0 & 0 & 1 & 0 & 0 & 0 & 0\\
 0 & 0 & 1 & 0 & 0 & 0 & 0 & 0\\
 0 & 0 & 0 & 0 & 1 & 0 & 0 & 0\\
 0 & 0 & 0 & 0 & 0 & 0 & 1 & 0\\
 0 & 0 & 0 & 0 & 0 & 1 & 0 & 0\\
 0 & 0 & 0 & 0 & 0 & 0 & 0 & 1
\end{pmatrix},
\end{align*}
as $\tau$ on a small angle just above the negative reals $($region {\rm II}$)$,
\[
\mathbf{Z}_{0}(\tq) =s\bigl(\Phih\bigr)(-1/\tau)S_{\rm II}(q),
\]
where
\begin{align*}
S_{\rm II}(q) ={}&\begin{pmatrix}
 1 & 0 & 0 & 0 & 0 & 0 & 0 & 0\\
 0 & 1 & 0 & 0 & 0 & 0 & 0 & 0\\
 0 & 0 & 1 & 0 & 0 & 0 & 0 & 0\\
 0 & 0 & 0 & 1 & 0 & 0 & 0 & 0\\
 0 & 0 & 0 & 0 & 1 & 0 & 0 & 0\\
 0 & 0 & 0 & 0 & 0 & 1 & 0 & 0\\
 0 & 0 & 0 & 0 & 0 & 0 & 0 & 1\\
 0 & 0 & 0 & 0 & 0 & 0 & 1 & 0
\end{pmatrix}P(q)\mathbf{Z}_{-3}(q)\\
&\times
\begin{pmatrix}
 1 & 0 & 0 & 0 & 0 & 0 & 0 & 0\\
 0 & -1 & 0 & 0 & 0 & 0 & 0 & 0\\
 0 & 0 & -1 & 0 & 0 & 0 & 0 & 0\\
 0 & 0 & 0 & -1 & 0 & 0 & 0 & 0\\
 0 & 0 & 0 & 0 & 1 & 0 & 0 & 0\\
 0 & 0 & 0 & 0 & 0 & 1 & 0 & 0\\
 0 & 0 & 0 & 0 & 0 & 0 & 1 & 0\\
 0 & 0 & 0 & 0 & 0 & 0 & 0 & -1
\end{pmatrix}
\begin{pmatrix}
 1 & 0 & 0 & 0 & 0 & 0 & 0 & 0\\
 0 & \tau & 0 & 0 & 0 & 0 & 0 & 0\\
 0 & 0 & 0 & 1 & 0 & 0 & 0 & 0\\
 0 & 0 & 1 & 0 & 0 & 0 & 0 & 0\\
 0 & 0 & 0 & 0 & 1 & 0 & 0 & 0\\
 0 & 0 & 0 & 0 & 0 & 0 & 1 & 0\\
 0 & 0 & 0 & 0 & 0 & 1 & 0 & 0\\
 0 & 0 & 0 & 0 & 0 & 0 & 0 & 1
\end{pmatrix},
\end{align*}
as $\tau$ on a small angle just below the negative reals $($region {\rm III}$)$,
\[
\mathbf{Z}_{0}(\tq) =s\bigl(\Phih\bigr)(-1/\tau)S_{\rm III}(q),
\]
where
\begin{align*}
S_{\rm III}(q) ={} &\begin{pmatrix}
 1 & 0 & 0 & 0 & 0 & 0 & 0 & 0\\
 0 & 1 & 0 & 0 & 0 & 0 & 0 & 0\\
 0 & 0 & 1 & 0 & 0 & 0 & 0 & 0\\
 0 & 0 & 0 & 1 & 0 & 0 & 0 & 0\\
 0 & 0 & 0 & 0 & 1 & 0 & 0 & 0\\
 0 & 0 & 0 & 0 & 0 & 1 & 0 & 0\\
 0 & 0 & 0 & 0 & 0 & 0 & 0 & 1\\
 0 & 0 & 0 & 0 & 0 & 0 & 1 & 0
\end{pmatrix}P(q)\mathbf{Z}_{-3}(q)\\
&\times
\begin{pmatrix}
 1 & 0 & 0 & 0 & 0 & 0 & 0 & 0\\
 0 & 1 & 0 & 0 & 0 & 0 & 0 & 0\\
 0 & 0 & 1 & 0 & 0 & 0 & 0 & 0\\
 0 & 0 & 0 & 1 & 0 & 0 & 0 & 0\\
 0 & 0 & 0 & 0 & -1 & 0 & 0 & 0\\
 0 & 0 & 0 & 0 & 0 & -1 & 0 & 0\\
 0 & 0 & 0 & 0 & 0 & 0 & -1 & 0\\
 0 & 0 & 0 & 0 & 0 & 0 & 0 & 1
\end{pmatrix}
\begin{pmatrix}
 1 & 0 & 0 & 0 & 0 & 0 & 0 & 0\\
 0 & \tau^{-2} & 0 & 0 & 0 & 0 & 0 & 0\\
 0 & 0 & 0 & \tau^{-1} & 0 & 0 & 0 & 0\\
 0 & 0 & \tau^{-1} & 0 & 0 & 0 & 0 & 0\\
 0 & 0 & 0 & 0 & \tau^{-1} & 0 & 0 & 0\\
 0 & 0 & 0 & 0 & 0 & 0 & \tau^{-1} & 0\\
 0 & 0 & 0 & 0 & 0 & \tau^{-1} & 0 & 0\\
 0 & 0 & 0 & 0 & 0 & 0 & 0 & \tau^{-1}
\end{pmatrix},
\end{align*}
as $\tau$ on a small angle just below the positive reals $($region {\rm IV}$)$,
\[
\mathbf{Z}_{0}(\tq) =s\bigl(\Phih\bigr)(-1/\tau)S_{\rm IV}(q),
\]
where
\begin{align*}
S_{\rm IV}(q) ={}&\begin{pmatrix}
 -1 & 0 & 0 & 0 & 0 & 0 & 0 & 0\\
 0 & 1 & 0 & 0 & 0 & 0 & 0 & 0\\
 0 & 0 & 1 & 0 & 0 & 0 & 0 & 0\\
 0 & 0 & 0 & 1 & 0 & 0 & 0 & 0\\
 0 & 0 & 0 & 0 & 1 & 0 & 0 & 0\\
 0 & 0 & 0 & 0 & 0 & 1 & 0 & 0\\
 0 & 0 & 0 & 0 & 0 & 0 & -1 & 0\\
 0 & 0 & 0 & 0 & 0 & 0 & 0 & -1
\end{pmatrix}P(q)\mathbf{Z}_{-3}(q)\\
&\times
\begin{pmatrix}
 -1 & 0 & 0 & 0 & 0 & 0 & 0 & 0\\
 0 & -1 & 0 & 0 & 0 & 0 & 0 & 0\\
 0 & 0 & 1 & 0 & 0 & 0 & 0 & 0\\
 0 & 0 & 0 & 1 & 0 & 0 & 0 & 0\\
 0 & 0 & 0 & 0 & -1 & 0 & 0 & 0\\
 0 & 0 & 0 & 0 & 0 & -1 & 0 & 0\\
 0 & 0 & 0 & 0 & 0 & 0 & -1 & 0\\
 0 & 0 & 0 & 0 & 0 & 0 & 0 & 1
\end{pmatrix}
\begin{pmatrix}
 1 & 0 & 0 & 0 & 0 & 0 & 0 & 0\\
 0 & \tau^{-2} & 0 & 0 & 0 & 0 & 0 & 0\\
 0 & 0 & 0 & \tau^{-1} & 0 & 0 & 0 & 0\\
 0 & 0 & \tau^{-1} & 0 & 0 & 0 & 0 & 0\\
 0 & 0 & 0 & 0 & \tau^{-1} & 0 & 0 & 0\\
 0 & 0 & 0 & 0 & 0 & 0 & \tau^{-1} & 0\\
 0 & 0 & 0 & 0 & 0 & \tau^{-1} & 0 & 0\\
 0 & 0 & 0 & 0 & 0 & 0 & 0 & \tau^{-1}
\end{pmatrix}\!.
\end{align*}
\end{Conjecture}
Of course, these equations are between matrices and therefore we can invert them. This allows us to use these equations to give conjectural formulae for the Borel resummation as a product of~$q$ and~$\tq$ series weighted by powers of $\tau$. These conjectures are analogous to those in~\cite{GGM:I,GGMW} on the Borel transform. Moreover, using the analytic properties of the state integrals, we see that the Borel resummation on various rays have analytic continuation to a~cut plane. This means that they can be compared on the upper half plane or the lower half plane.\looseness=1

\subsection{Stokes matrices}
\label{stomat}

Using the previous Conjecture~\ref{conj:borel} for the Borel resummation, we can compute the Stokes phenomenon between these regions and find agreement with equation~\eqref{eq:closeststokes}. Assuming various conjectures on the behaviour of the Stokes phenomenon~\cite{GGM:I}, we can compute the generating series of Stokes matrices as quotients of these matrices. Indeed, to go from region $\rm II$ to region $\rm I$ we must multiply $s\bigl(\widehat{\Phi}\bigr)$ by $S_{\rm II}(q)S_{I}(q)^{-1}$. This then stores the Stokes constants the upper half plane in the matrix of $q$-series{\setlength{\arraycolsep}{1.5pt}
\begin{align*}
 &\mathsf{S}_{-}(q) =S_{\rm II}(q)S_{\rm I}(q)^{-1}-I\\
 &
 =\begin{pmatrix}
 0 & 0 & 0 & 0 & 0 & 0 & 0 & 0\\
 -1\!+\!q\!+\!3q^2 & -q\!-\!2q^2 & 1\!+\!q & -q^2 & q\!+\!2q^2 & -1\!-\!q\!+\!q^2 & -1\!+\!3q^2 & -1\!+\!3q^2 \\
 q\!-\!q^2 & q^2 & -q\!-\!q^2 & 0 & -q^2 & q\!+\!q^2 & q & q \\
 -1\!+\!2q\!+\!q^2 & -q\!-\!q^2 & 1 & -q^2 & q\!+\!q^2 & -1\!+\!q^2 & -1\!+\!q\!+\!2q^2 & -1\!+\!q\!+\!2q^2 \\
 1\!-\!2q\!-\!q^2 & q\!+\!q^2 & -1 & q^2 & -q\!-\!q^2 & 1\!-\!q^2 & 1\!-\!q\!-\!2q^2 & 1\!-\!q\!-\!2q^2 \\
 -q\!+\!q^2 & -q^2 & q\!+\!q^2 & 0 & q^2 & -q\!-\!q^2 & -q & -q \\
 1\!-\!3q\!-\!q^2 & q & -1\!+\!q\!+\!2q^2 & q^2 & -q & 1\!-\!q\!-\!3q^2 & -2q\!-\!3q^2 & -2q\!-\!3q^2 \\
 1\!-\!3q\!-\!q^2 & q & -1\!+\!q\!+\!2q^2 & q^2 & -q & 1\!-\!q\!-\!3q^2 & -2q\!-\!3q^2 & -2q\!-\!3q^2
\end{pmatrix}\\
&+O\bigl(q^3\bigr).
\end{align*}}%
The Stokes constants in the lower half plane are analogously stored by{\setlength{\arraycolsep}{1.5pt}
\begin{align*}
 &\mathsf{S}_{+}(q) =S_{\rm IV}\bigl(q^{-1}\bigr)S_{\rm III}\bigl(q^{-1}\bigr)^{-1}-I\\
 &
 =\begin{pmatrix}
 0 & 0 & 0 & 0 & 0 & 0 & 0 & 0\\
 0 & -q\!-\!2q^2 & q^2 & -q\!-\!q^2 & q\!+\!q^2 & -q^2 & -q & -q\\
 -1 & 1\!+\!q & -q\!-\!q^2 & 1 & -1 & q\!+\!q^2 & 1\!-\!q\!-\!2q^2 & 1\!-\!q\!-\!2q^2\\
 -q\!+\!q^2 & -q^2 & 0 & -q^2 & q^2 & 0 & -q^2 & -q^2\\
 -q\!-\!q^2 & q\!+\!2q^2 & -q^2 & q\!+\!q^2 & -q\!-\!q^2 & q^2 & q & q\\
 -2q^2 & -1\!-\!q\!+\!q^2 & q\!+\!q^2 & -1\!+\!q^2 & 1\!-\!q^2 & -q\!-\!q^2 & -1\!+\!q\!+\!3q^2 & -1\!+\!q\!+\!3q^2\\
 2q^2 & 1\!-\!3q^2 & -q & 1\!-\!q\!-\!2q^2 & -1\!+\!q\!+\!2q^2 & q & -2q\!-\!3q^2 & -2q\!-\!3q^2\\
 q\!+\!2q^2 & 1\!-\!3q^2 & -q & 1\!-\!q\!-\!2q^2 & -1\!+\!q\!+\!2q^2 & q & -2q\!-\!3q^2 & -2q\!-\!3q^2
 \end{pmatrix}\\
 &  +O\bigl(q^3\bigr) .
\end{align*}}%
More precisely, to compute the Stokes constants from these matrices one must factorise the above $q$-series into elementary Stokes automorphism~\cite[equations~(66)--(68)]{GGM:I} associated to each argument and take the logarithm, which can be done to any order in $q$ given that these elementary pieces are unipotent.
Factorising these matrices, we find that{\setlength{\arraycolsep}{1.5pt}
\begin{align*}
&I+\mathsf{S}_{-}(q) =\\
& \begin{pmatrix}
1 & 0 & 0 & 0 & 0 & 0 & 0 & 0\\
0 & 1 & 0 & 0 & 0 & 0 & 0 & 0\\
0 & 0 & 1 & 0 & 0 & 0 & 0 & 0\\
0 & 0 & 0 & 1 & 0 & 0 & 0 & 0\\
0 & 0 & 0 & 0 & 1 & 0 & 0 & 0\\
0 & 0 & 0 & 0 & 0 & 1 & 0 & 0\\
0 & 0 & 0 & 0 & 0 & 0 & 1 & 0\\
1\!-\!q & q\!+\!2q^2 & -1\!-\!q\!-\!3q^2 & q^2 & -q\!-\!2q^2 & 1\!+\!q\!+\!2q^2 & 2q\!+\!7q^2 & 1\!+\!2q\!+\!7q^2
\end{pmatrix}\\[-0.5mm]
&\begin{pmatrix}
1 & 0 & 0 & 0 & 0 & 0 & 0 & 0\\
0 & 1 & 0 & 0 & 0 & 0 & -1\!-\!2q\!-\!4q^2 & -1\!-\!2q\!-\!4q^2\\
0 & 0 & 1 & 0 & 0 & 0 & q\!+\!2q^2\!+\!4q^3 & q\!+\!2q^2\\
0 & 0 & 0 & 1 & 0 & 0 & -1\!-\!q\!-\!3q^2 & -1\!-\!q\!-\!3q^2\\
0 & 0 & 0 & 0 & 1 & 0 & 1\!+\!q\!+\!3q^2 & 1\!+\!q\!+\!3q^2\\
0 & 0 & 0 & 0 & 0 & 1 & -q\!-\!2q^2 & -q\!-\!2q^2\\
0 & 0 & 0 & 0 & 0 & 0 & 1 & 0\\
0 & 0 & 0 & 0 & 0 & 0 & 0 & 1
\end{pmatrix}\\[-0.5mm]
&\begin{pmatrix}
1 & 0 & 0 & 0 & 0 & 0 & 0 & 0\\
0 & 1 & 0 & 0 & 0 & 0 & 0 & 0\\
0 & 0 & 1 & 0 & 0 & 0 & 0 & 0\\
0 & 0 & 0 & 1 & 0 & 0 & 0 & 0\\
0 & 0 & 0 & 0 & 1 & 0 & 0 & 0\\
0 & 0 & 0 & 0 & 0 & 1 & 0 & 0\\
1\!-\!3q\!-\!q^2 & q & -1\!+\!q\!+\!2q^2 & q^2 & -q & 1\!-\!q\!-\!3q^2 & -2q\!-\!3q^2\\
0 & 0 & 0 & 0 & 0 & 0 & 0 & 1
\end{pmatrix}
+O\bigl(q^3\bigr) ,\\[-0.5mm]
&I+\mathsf{S}_{+}(q) =\\[-0.5mm]
&\begin{pmatrix}
1 & 0 & 0 & 0 & 0 & 0 & 0 & 0\\
0 & 1 & 0 & 0 & 0 & 0 & 0 & 0\\
0 & 0 & 1 & 0 & 0 & 0 & 0 & 0\\
0 & 0 & 0 & 1 & 0 & 0 & 0 & 0\\
0 & 0 & 0 & 0 & 1 & 0 & 0 & 0\\
0 & 0 & 0 & 0 & 0 & 1 & 0 & 0\\
-q\!-\!4q^2 & 1\!+\!2q\!+\!4q^2 & -q\!-\!2q^2 & 1\!+\!q\!+\!3q^2 & -1\!-\!q\!-\!3q^2 & q\!+\!2q^2 & 1\!+\!2q\!+\!7q^2 & 2q\!+\!7q^2\\
0 & 0 & 0 & 0 & 0 & 0 & 0 & 1\\
\end{pmatrix}\\[-0.5mm]
&\begin{pmatrix}
1 & 0 & 0 & 0 & 0 & 0 & 0 & 0\\
q^2 & 1 & 0 & 0 & 0 & 0 & -q\!-\!2q^2 & -q\!-\!2q^2\\
-1\!-\!q\!-\!3q^2 & 0 & 1 & 0 & 0 & 0 & 1\!+\!q\!+\!3q^2 & 1\!+\!q\!+\!3q^2\\
-q\!+\!q^2 & 0 & 0 & 1 & 0 & 0 & -q^2\!+\!0 & -q^2\\
-q\!-\!2q^2 & 0 & 0 & 0 & 1 & 0 & q\!+\!2q^2\!+\!0 & q\!+\!2q^2\\
q\!+\!q^2 & 0 & 0 & 0 & 0 & 1 & -1\!-\!q\!-\!2q^2 & -1\!-\!q\!-\!2q^2\\
\end{pmatrix}\\[-0.5mm]
&\begin{pmatrix}
1 & 0 & 0 & 0 & 0 & 0 & 0 & 0\\
0 & 1 & 0 & 0 & 0 & 0 & 0 & 0\\
0 & 0 & 1 & 0 & 0 & 0 & 0 & 0\\
0 & 0 & 0 & 1 & 0 & 0 & 0 & 0\\
0 & 0 & 0 & 0 & 1 & 0 & 0 & 0\\
0 & 0 & 0 & 0 & 0 & 1 & 0 & 0\\
0 & 0 & 0 & 0 & 0 & 0 & 1 & 0\\
q\!+\!2q^2 & 1\!-\!3q^2 & -q & 1\!-\!q\!-\!2q^2 & -1\!+\!q\!+\!2q^2 & q & -2q\!-\!3q^2 & 1\!-\!2q\!-\!3q^2
\end{pmatrix} .
\end{align*}}%
We can read off the Stokes constants from the entries of the matrices appearing in this factorisation.

\subsection{The 3d index}
In~\cite{DGG}, remarkable $q$-series invariants of three-manifolds were proposed. These were then shown invariant under certain 2--3 moves when the associated $q$-series were convergent. Using state integrals, \cite{GK:3dind} showed that this is a true topological invariant of the manifold. As with many three-manifold invariants, the 3d index is built from gluing together pieces associated to tetrahedra. This piece is the so called tetrahedral index and given in a certain basis by the $q$-series
\[
 \mathcal{I}_{\Delta}(m,e;q)
  =
 \sum_{n=\max(-e,0)}^{\infty}(-1)^{n}\frac{q^{n(n+1)/2-(n+\frac{1}{2}e)m}}{(q;q)_{n}(q;q)_{n+e}} .
\]
Using Neumann--Zagier data of an ideal triangulation~\cite{NeuZag} of a manifold $M$ that supports an angle structure~\cite{Gar3d}, the 3d index of~\cite{DGG} can be defined by taking
\[
 \mathcal{I}_{M}(m,e;q)
  =
 \sum_{k\in\BZ^{r}}q^{\frac{1}{2}\nu\cdot k}\prod_{j=1}^{N}\mathcal{I}_{\Delta}((dm-be)_j-B_j\cdot k,(-cm+ae)_j+A_j\cdot k;q),
\]
where $r$ is the number of tetrahedra minus the number of cusps, $m,e\in\BZ^{N-r}$, $\alpha,\beta,\gamma,\delta\in M_{N-r,N}(\BZ)$ give the gluing equations for the cusps, and $A,B\in M_{r,N}(\BZ)$ and $\nu\in\BZ^{r}$ give the gluing equations for the edges. It was conjectured~\cite{GGM:I,GrassiGM,Marino} that the Stokes constants of asymptotic series associated to knots were related to the 3d index. It is natural to attempt to understand if a similar story could hold for closed $3$-manifolds.

Recently, a surgery formula was proposed for the 3d index~\cite{Gang3d}. This proposal does not enjoy a~proof of topological invariance, unlike the 3d index for knots. However, there is some numerical evidence~\cite{Gang3d} of topological invariance. Using the computation of the Stokes constants from Section~\ref{stomat} and the conjecture of~\cite[Conjecture 2]{GGM:I}, we have a natural guess for the 3d index of $4_1(-1,2)$ given by
\[
 1+\mathsf{S}_{+}^{(\rho_6,\rho_6)}(q)
  =
 1 - 2q - 3q^2 + 3q^4 + 10q^5 + 14q^6 + 22q^7 + 20q^8 + 14q^9+\cdots .
\]
Therefore, we can compare\footnote{Upon posting the original version of this work on the arXiv, I received an email from Dongmin Gang with a question of whether the Stokes constants matched their 3d index, which for this manifold was given $1-2q -3q^2 +3q^4 +10 q^5+14 q^6+\cdots$.} with the proposal of~\cite{Gang3d}. For the knot $4_1$, we take $\alpha=[1,0]$, $\beta=[0,-1]$, $\gamma=[-1,0]$, $\delta=[1,0]$, $A=[-1,-1]$, $B=[1,1]$, $\nu=0$. We will use the notation of~\cite{Gang3d}, so that for $m,e\in\BZ$ and $[a,b;c,d]\in\SL_{2}(\BZ)$,
\begin{gather*}
 \mathcal{K}^{\mathrm{SO}(3)}(m,e,c,d;q)
  :=
 \frac{1}{2}(-1)^{am+2be}\bigl(\delta_{cm+2de,0}\bigl(q^{\frac{am+2be}{2}}-q^{-\frac{am+2be}{2}}\bigr)-\delta_{|cm+2de|,2}\bigr) ,\\
 \mathcal{I}_{S^{3}\smallsetminus \mathbf{4_1}}^{(\mu,\lambda)}(m,e;q)
 =
 \sum_{e_1\in\BZ}\mathcal{I}_{\Delta}(m-e_1,m+e-e_1;q)\mathcal{I}_{\Delta}(e-e_1,-e_1;q),\\
 \mathcal{I}_{4_{1}(-1,2)}(q)
 :=\sum_{m,e\in\BZ}
 \mathcal{K}^{\mathrm{SO}(3)}(m,e,-1,2;q)\mathcal{I}_{S^{3}\smallsetminus \mathbf{4_1}}^{(\mu,\lambda)}(m,e;q)\\
\hphantom{\mathcal{I}_{4_{1}(-1,2)}(q) \ }{} = 1 - 2q - 3q^2 + 3q^4 + 10q^5 + 14q^6 + 22q^7 + 20q^8  + 14q^9 +\cdots.
\end{gather*}
As can be seen, this expansion agrees with the expansion of $1+\mathsf{S}_{+}^{(\rho_6,\rho_6)}(q)$ where we recall that $\rho_{6}$ is the geometric connection. This gives more numerical evidence for the proposals of both~\cite{Gang3d,GGM:I}.

\appendix

\section{\texorpdfstring{Modularity of the $\boldsymbol{q}$-Pochhammer symbol}{Modularity of the q-Pochhammer symbol}}\label{sec:mod.poch}

The $q$-Pochhammer symbol satisfies various modularity properties that it inherits from the Faddeev quantum dilogarithm. A fantastic description of Faddeev's quantum dilogarithm is given in~\cite[Appendix A]{AK:I} however we will include some additional formulae that we will need. There are three different formulae for the Faddeev quantum dilogarithm on various domains. Using the notation $\tau=\mathsf{b}^{2}$, $c_{\mathsf{b}}=\frac{{\rm i}}{2}(\mathsf{b}+\mathsf{b}^{-1})$ and the standard $q=\e(\tau)$ and $\tq=\e(-1/\tau)$, we have (when defined)
\begin{align*}
 \Phi_{\mathsf{b}}(x)
 ={}&
 \exp\biggl(\int_{\BR+{\rm i}0}\frac{\exp(-2{\rm i}xw)}{4\sinh(\mathsf{b}w)\sinh\bigl(\mathsf{b}^{-1}w\bigr)}\frac{{\rm d}w}{w}\!\biggr)
 =
 \bigl(-q^{1/2}\exp(2\pi\mathsf{b}x);q\bigr)_{-\lfloor\Im(\mathsf{b}^{-1}x+{\rm i}\frac{\mathsf{b}^{-2}}{2})+\frac{1}{2}\rfloor}\\
 & \times\exp\biggl(\mathsf{b}^{2}\!\int_{{\rm i}\BR+0}\frac{\log\bigl(1+\tq^{-1/2}\e(-w)\exp(2\pi\mathsf{b}^{-1}x)\bigr)}{1-\e\bigl(-\mathsf{b}^{2}w\bigr)}{\rm d}w\biggr)
 \\
 ={}&
 \frac{\bigl(-q^{1/2}\exp(2\pi\mathsf{b}x);q\bigr)_{\infty}}{\bigl(-\tq^{1/2}\exp(2\pi\mathsf{b}^{-1}x);\tq\bigr)_{\infty}}
  .
\end{align*}
These can be proved to be equivalent where their domains overlap by use of the residue theorem. For example, when $|\Im(x)|<\Im(c_{\mathsf{b}})$ and $\Im(\tau)>0$ the first and the last equation can be shown equal as
\begin{align*}
 &\int_{\BR+{\rm i}0}\frac{\exp(-2{\rm i}xw)}{4\sinh(\mathsf{b}w)\sinh\bigl(\mathsf{b}^{-1}w\bigr)}\frac{{\rm d}w}{w}
   =-\int_{{\rm i}\BR+0}\frac{\e\bigl(-{\rm i}xw+\bigl(\mathsf{b}+\mathsf{b}^{-1}\bigr)w/2\bigr)}{(1-\e(\mathsf{b}w))\bigl(1-\e\bigl(\mathsf{b}^{-1}w\bigr)\bigr)}\frac{{\rm d}w}{w}\\
 & \qquad{} =2\pi {\rm i}\sum_{k=1}^{\infty}\bigl(\mathrm{Res}_{w=k\mathsf{b}}+\mathrm{Res}_{w=k\mathsf{b}^{-1}}\bigr)\frac{\e\bigl(-{\rm i} xw+\bigl(\mathsf{b}+\mathsf{b}^{-1}\bigr)w/2\bigr)}{(1-\e(\mathsf{b}w))\bigl(1-\e\bigl(\mathsf{b}^{-1}w\bigr)\bigr)}\frac{{\rm d}w}{w}\\
 & \qquad{}=
 \sum_{k=1}^{\infty}(-1)^{k}\tq^{k/2}\frac{\e\bigl(-k{\rm i}x\mathsf{b}^{-1}\bigr)}{k(1-\tq^{k})}-(-1)^{k}q^{k/2}\frac{\e(-k{\rm i}x\mathsf{b})}{k(1-q^{k})}\\
 & \qquad{} =
 -\log\bigl(-\tq^{1/2}\exp(2\pi\mathsf{b}^{-1}x);\tq\bigr)_{\infty}+\log\bigl(-q^{1/2}\exp(2\pi\mathsf{b}x);q\bigr)_{\infty},
\end{align*}
where the last line follows from the equation~\cite[Proposition~2]{Zagier:Dilog}
\[
 -\log(x;q)_{\infty}
  =
 \sum_{k=0}^{\infty}\frac{x^k}{k(1-q^{k})} .
\]
Then for the second and last when $\Re\bigl(\mathsf{b}^{2}\bigr)>0$ and $\Im\bigl(\mathsf{b}^{2}\bigr)>0$, we use the integral for $n\in\BZ$
\[
 -\int_{-i\mathsf{b}^{-1}x+\frac{\mathsf{b}^{-2}}{2}-\frac{1}{2}-n}^{{\rm i}\infty}\frac{2\pi {\rm i}}{1-\e\bigl(-\mathsf{b}^{2}w\bigr)}{\rm d}w
  =
 \mathsf{b}^{-2}\log\bigl(1+q^{1/2+n}\exp(2\pi\mathsf{b}x)\bigr) .
\]
We see that for the principle branch of the logarithm
\begin{align*}
 &\mathsf{b}^{2} \int_{{\rm i}\BR+0}\frac{\log\bigl(1+\tq^{-1/2}\e(-w)\exp(2\pi\mathsf{b}^{-1}x)\bigl)}{1-\e(-\mathsf{b}^{2}w)}{\rm d}w\\
 & =
 -\sum_{k=1}^{\infty}\mathrm{Res}_{w=k\mathsf{b}^{-2}}\frac{\log\bigl(1+\tq^{-1/2}\e(-w)\exp(2\pi\mathsf{b}^{-1}x)\bigr)} {1-\e\bigl(-\mathsf{b}^{2}w\bigr)}\frac{{\rm d}w}{\mathsf{b}^{-2}}\\
 &\qquad{}-\sum_{n=-\lfloor\Im(\mathsf{b}^{-1}x+{\rm i} \frac{\mathsf{b}^{-2}}{2})+\frac{1}{2}\rfloor}^{\infty}\int_{-{\rm i}\mathsf{b}^{-1}x+\frac{\mathsf{b}^{-2}}{2}-\frac{1}{2}-n}^{{\rm i}\infty}\frac{2\pi {\rm i}}{1-\e\bigl(-\mathsf{b}^{2}w\bigr)}\frac{{\rm d} w}{\mathsf{b}^{-2}}\\
 & =
 -\sum_{k=0}^{\infty}\log\bigl(1+\tq^{1/2+k}\exp\bigl(2\pi\mathsf{b}^{-1}x\bigr)\bigr)
 +\sum_{n=0}^{\infty}\log\bigl(1+q^{1/2+n}\exp(2\pi\mathsf{b}x)\bigr)\\
 &\qquad{}-\sum_{n=0}^{-\lceil\Im(\mathsf{b}^{-1}x+{\rm i}\frac{\mathsf{b}^{-2}}{2})+\frac{1}{2}\rceil}\log\bigl(1+q^{1/2+n}\exp(2\pi\mathsf{b}x)\bigr)\\
 & =
 -\log\bigl(-\tq^{1/2}\exp(2\pi\mathsf{b}^{-1}x);\tq\bigr)_{\infty}
 +\log\bigl(-q^{1/2}\exp(2\pi\mathsf{b}x);q\bigr)_{\infty}\\
 &\qquad{}-\log\bigl(-q^{1/2}\exp(2\pi\mathsf{b}x);q\bigr)_{-\lfloor\Im(\mathsf{b}^{-1}x+{\rm i}\frac{\mathsf{b}^{-2}}{2})+\frac{1}{2}\rfloor}.
\end{align*}
These various descriptions can be used to show that $\Phi_{\mathsf{b}}(x)$ extends to a meromorphic function for $\mathsf{b}^{2}\in\BC\smallsetminus \BR_{\leq0}$ with simple poles and zeros determined by
\[
 \Phi_{\mathsf{b}}(x)^{\pm} =0
 \qquad\text{if and only if}\qquad
 x\in\mp(c_{\mathsf{b}}+{\rm i}\mathsf{b}\BZ_{\geq0}+{\rm i}\mathsf{b}^{-1}\BZ_{\geq0}) .
\]
Therefore, for $k\in\BZ$ and $\Im\bigl(\mathsf{b}^{2}\bigr)>0$ we have
\begin{align*}
 &\Phi_{\mathsf{b}}\biggl({\rm i}\biggl(\frac{1}{2}\biggr)\mathsf{b}+{\rm i} \biggl(-k-\frac{1}{2}\biggr)\mathsf{b}^{-1}\biggr)\Phi_{\mathsf{b}}\biggl({\rm i} \biggl(\frac{1}{2}\biggr)\mathsf{b}+{\rm i} \biggl(-\frac{1}{2}\biggr)\mathsf{b}^{-1}\biggr)^{-1}\\
 & \qquad =
 \frac{(q;q)_{\infty}}{\bigl(\tq^{k+1};\tq\bigr)_{\infty}}\frac{(\tq;\tq)_{\infty}}{(q;q)_{\infty}}
  =
 (\tq;\tq)_{k}\\
 & \qquad =
 \mathsf{b}(q;q)_{\lfloor\Re(\mathsf{b}^{-2})k\rfloor}
 \exp\biggl(\frac{2\pi {\rm i}}{24\mathsf{b}^{2}}-\frac{\pi {\rm i}}{4}+\frac{2\pi {\rm i}\mathsf{b}^{2}}{24}
 +\mathsf{b}^{2} \int_{{\rm i}\BR+0}\frac{\log\bigl(1-\tq^{k}\e(-w)\bigr)}{1-\e\bigl(-\mathsf{b}^{2}w\bigr)}{\rm d}w\biggr) ,
\end{align*}
and given the analyticity of the quantum dilogarithm and the Pochhammer symbols the first and second line hold for all $\mathsf{b}^{2}\in\BC\smallsetminus \BR_{\leq0}$. Therefore, let
\be\label{eq:f}
 f(k;\tau)
  =
 \frac{\tau}{24}-\frac{1}{8}+\frac{1}{24\tau}+\frac{\tau}{2\pi {\rm i}} \int_{{\rm i}\BR+0}\frac{\log(1-\e(-k/\tau-w))}{1-\e(-\tau w)}{\rm d}w .
\ee
This function is meromorphic in $k$ and $\tau$ for $\Re(k/\tau)\notin\BZ$ and satisfies the periodicity
\be\label{eq:f.per}
 f(k+\tau;\tau) =f(k;\tau) .
\ee
Then we can write
\be\label{eq:f.qp}
 (\tq;\tq)_{k}
  =
 (q;q)_{\lfloor\Re(k/\tau)\rfloor}
 \e(f(k;\tau)) .
\ee
We can determine the asymptotics of $f$ by considering the equation
\begin{align*}
 &\frac{\tau}{2\pi {\rm i}}\!\int_{{\rm i}\BR+0}\frac{\log(1-\e(-k/\tau-w))}{1-\e(-\tau w)}{\rm d}w
  =
 \frac{\tau}{2\pi} \int_{\BR+{\rm i}0}\frac{\log(1-\exp(-2\pi {\rm i}k/\tau+2\pi w))}{1-\exp(2\pi \tau w)}{\rm d}w\\
 &=
 \frac{\tau}{2\pi} \int_{0+{\rm i}0}^{\infty}\frac{\log(1-\exp(-2\pi {\rm i}k/\tau+2\pi w))}{1-\exp(2\pi \tau w)}{\rm d}w\\
 & \qquad
 -
 \frac{\tau}{2\pi} \int_{0-i0}^{\infty}\frac{\log(1-\exp(-2\pi {\rm i}k/\tau-2\pi w))}{1-\exp(-2\pi \tau w)}{\rm d}w\\
 &=
 \frac{\tau}{2\pi} \int_{0}^{\infty}\log(1-\exp(-2\pi {\rm i}k/\tau-2\pi w)){\rm d}w
 +\frac{1}{4\pi i}\log(1-\exp(-2\pi {\rm i}k/\tau))\\
 &\qquad+\frac{\tau}{2\pi}\!\int_{0}^{\infty}\frac{\log(1-\exp(-2\pi {\rm i}k/\tau+2\pi w))-\log(1-\exp(-2\pi {\rm i}k/\tau-2\pi w))}{1-\exp(2\pi \tau w)}{\rm d}w\\
 &=
 \frac{\tau}{(2\pi {\rm i})^2}\Li_{2}(\e(-k/\tau))
 +\frac{1}{2(2\pi {\rm i})}\log(1-\exp(-2\pi {\rm i}k/\tau))\\
 &\qquad+\frac{{\rm i}\tau}{2\pi {\rm i}} \int_{0}^{\infty}\frac{\log(1-\exp(-2\pi {\rm i}k/\tau+2\pi w))-\log(1-\exp(-2\pi {\rm i}k/\tau-2\pi w))}{1-\exp(2\pi \tau w)}{\rm d}w .
\end{align*}
This can be used to show that
\be\label{eq:asymp.f}
 f(k\tau,\tau) =\frac{\tau}{(2\pi {\rm i})^2}\Li_2(\e(-k))+\frac{\tau}{24}+O\bigl(\tau^0\bigr) .
\ee

\section{Asymptotic methods}\label{sec:asymp.meth.ap}

In this appendix, we will provide more details on the proofs of the various theorems on asymptotics. In particular we will give more details on Theorems~\ref{thm:zhatasym} and~\ref{thm:qmodold}. In upcoming work with Fantini, more refined methods are developed, which allow for even stronger results so we only outline these older methods.

\subsection{Outline of the proof of Theorem~\ref{thm:zhatasym}}

In this appendix, we will gives a few details on the proof of Theorem~\ref{thm:zhatasym}.
This involves very similar arguments as the next Appendix~\ref{sec:more.deets}.
For $\tq=\e(i/r)$ for $r\in\BR_{>0}$, we can rewrite the sum
\[
 \sum_{0\leq k\leq\ell}(-1)^{k+\ell}\frac{\tq^{\frac{1}{2}3k(k+1)+\frac{1}{2}\ell(\ell+1)-k}(\tq;\tq)_{\ell}}{(\tq;\tq)_{2k}(\tq;\tq)_{\ell-k}}
  =
 \sum_{k,j=0}^{\infty}(-1)^{j}\frac{\tq^{k(2k-1)+jk+\frac{1}{2}j(j+1)}(\tq;\tq)_{k+j}}{(\tq;\tq)_{2k}(\tq;\tq)_{j}}
\]
with the Faddeev quantum dilogarithm and explicitly $f$ from equation~\eqref{eq:f}
\begin{align*}
& \sum_{k,j=0}^{\infty}
 \exp\biggl(-\frac{2\pi}{r}\biggl(\frac{{\rm i}}{2}rj+k(2k-1)+jk+\frac{1}{2}j(j+1)\\
& \qquad{}-r{\rm i}f(k+j;{\rm i}/r)+r{\rm i}f(2k;{\rm i}/r)+r{\rm i}f(j;{\rm i}/r)\biggr)\biggr) .
\end{align*}
Notice that when $k,j>2r$ the summand is $\mathrm{O}(\exp(-2\pi r))$. Therefore, this sum is asymptotic, up to $\mathrm{O}(\exp(-2\pi r))$, to
\begin{align*}
&\sum_{x,y\in\frac{1}{r}\BZ\cap[0,2]^2}
 \exp\biggl(-2\pi\biggl(\frac{{\rm i}}{2}ry+x(2xr-1)+yxr+\frac{1}{2}y(yr+1)\\
& \qquad{}-{\rm i}f(xr+yr;{\rm i}/r)+{\rm i}f(2xr;{\rm i}/r)+{\rm i}f(yr;{\rm i}/r)\biggr)\biggr) .
\end{align*}
Using Abel--Plana summation, this sum is asymptotic to an integral with an appropriate para\-metri\-sa\-tion
\begin{align*}
 & \int_0^2\int_0^2
 \exp\biggl(-2\pi\biggl(\frac{{\rm i}}{2}ry+x(2xr-1)+yxr+\frac{1}{2}y(yr+1)\\
 &\qquad{} -{\rm i}f(xr+yr;{\rm i}/r)+{\rm i}f(2xr;{\rm i}/r)+{\rm i}f(yr;{\rm i}/r)\biggr)\biggr){\rm d}x{\rm d}y .
\end{align*}
Therefore, we are interested in the behaviour of the function
\begin{align*}
 W(x,y)
  ={}&
 2x^2+yx+\frac{1}{2}y^2-\frac{1}{2}y-\frac{1}{(2\pi {\rm i})^2}\Li_{2}(\e(-x-y))+\frac{1}{(2\pi {\rm i})^2}\Li_{2}(\e(-2x))\\
 &
 +\frac{1}{(2\pi {\rm i})^2}\Li_{2}(\e(-y))+\frac{1}{24} .
\end{align*}
This has one critical point on this branch given by $(x_3,y_3)=(0.0081371{\rm i},-0.017231{\rm i})$.
One can show that on the set $x\in[0,1]$
\[
 \Re\biggl(\frac{1}{(2\pi {\rm i})^2}\Li_{2}(\e(-x))\biggr)
\]
attains one maximum at $x=1/2$ given by
\[
 \Re\biggl(\frac{1}{(2\pi {\rm i})^2}\Li_{2}(-1)\biggr)
  =
 \frac{1}{48} ,
\]
with minima at $x=0,1$ given by
\[
 \Re\biggl(\frac{1}{(2\pi {\rm i})^2}\Li_{2}(1)\biggr)
  =
 -\frac{1}{24} .
\]
With this, we see that for $x,y\in[0,2]$
\[
 W(x,y)\geq2x^2+yx-\frac{1}{8}-\frac{1}{48}-\frac{1}{24}-\frac{1}{24}+\frac{1}{24},
\]
and for $y>1.25$, we have
\[
 W(x,y)\geq0 .
\]
Therefore, taking a parametrisation with end points beginning in the region where $W$ is greater than or equal to zero and using Fubini's theorem, we can deform the contour by the steepest descent of $W$ fixing the real parts of $x$, $y$.
In the region where $W(x,y)\geq W(x_3,y_3)=0.0029434\ldots=\bigl(\VC_{\rho_3}-4\pi^2\bigr)/(2\pi {\rm i})^2$. This deformation then passes through the critical point $(x_3,y_3)$ where the integrand attains its maximum.
Indeed, all possible critical points are located at $\BZ^2$ translates of the following numbers and besides $(x_3,y_3)$ they do not contribute:
\begin{align*}
 &(0.00813711\dots i,- 0.0172319\dots i) ,\\
 &(0.5 + 0.350548\dots i,
 0.5 - 0.141931\dots i) ,\\
 &(- 0.199662\dots i,
 0.5 + 0.0835868\dots i) ,\\
 &(0.5 - 0.042154\dots i,
 0.5 + 0.146002\dots i) ,\\
 &(- 0.111324\dots i,
 0.5 - 0.0478456\dots i) ,\\
 &(0.336754\ldots - 0.00277186\dots i,
 0.108538\ldots - 0.0112903\dots i) ,\\
 &(-0.336754\ldots - 0.00277186\dots i, -0.108538 -
 0.0112903\dots i) .
\end{align*}

\subsection{More details on Theorem~\ref{thm:qmodold}}\label{sec:more.deets}

In this appendix, we will provide some more details on the proof of Theorem~\ref{thm:qmodold}. These methods are well known and see~\cite{Ohtsuki:VCYVC} for another relevant example of these methods. We wish to apply a~stationary phase approximation to the integral of equation~\eqref{eq:int.of.sum}. The integrand has exponential behaviour as $r\to\infty$ and the behaviour is determined by the function $V_{\alpha,\beta}$ of equation~\eqref{eq:Vab}. The size of the growth is determined by $\Im(V_{\alpha,\beta})$. To begin with, we want to understand the regions where the size of $\Im(V_{\alpha,\beta})$ is larger than $\Im\bigl(-\frac{V_{\rho_6}}{(2\pi {\rm i})^2}\bigr)$. This is depicted in Figure~\ref{fig:oht.plot} and we prove explicitly this is true on the green regions.
\begin{figure}[t]\centering

\caption{The regions $R$ that arise from the stationary phase approximation. In \red{red} we have numerically computed the regions where $\Im(V_{\alpha,\beta})$ is larger than $\Im\bigl(-\frac{V_{\rho_6}}{(2\pi {\rm i})^2}\bigr)$ and in \green{green} is the region with an explicit proof that it is smaller given in Lemma~\ref{lem:greenregion}. The four regions $A$, $B$, $C$, $D$ give the connected components of the complement of the set $U$ depicted in green. Note the very small red bumps on the boundary in the regions $C$ and $D$.}\label{fig:oht.plot}
\end{figure}

\begin{Lemma}\label{lem:greenregion}
The region $U$ depicted in green in Figure~{\rm \ref{fig:oht.plot}} given by
\[
 U
  =
 \left\{\begin{array}{@{}c|l@{}}
 &2x+2y\in[0.05,0.33]+\BZ\quad \text{or}\\
 &x\in[0.67000,0.95000]+\BZ\quad \text{or}\\
 &x,y\in[0,0.025]\quad \text{or}\\
 &x\in[0.95,1] ,\quad y\in[0,0.025]\quad \text{or}\\
 (x,y)\in[0,1]& x\in[0,0.025] ,\quad y\in[0.95,1]\quad \text{or}\\
 &x\in[0.95,1] ,\quad y\in[0.025,0.075] ,\quad 2x+2y\leq 2.05\quad \text{or}\\
 &x,y\in[0.49250,1]\quad \text{or}\\
 &x\in[0.5,1] ,\quad 2x+2y\in[-0.15,0.5]+\BZ\quad \text{or}\\
 &y\in[0.5,1] ,\quad 2x+2y\in[-0.15,0.5]+\BZ
 \end{array}\right\}
\]
satisfies the property that for $(x,y)\in U$
\[
 \Im\biggl(\frac{V_{-1,1}(x,y)}{(2\pi {\rm i})^2}+\frac{V_{\rho_6}}{(2\pi {\rm i})^2}\biggr)>0 .
\]
\end{Lemma}
\begin{proof}
Firstly, we have
\[
 -\frac{V_{\rho_6}}{(2\pi {\rm i})^2}
  =
 0.12330 - 0.035425{\rm i} .
\]
Let
\[
 D(x)
  =
 \Im\biggl(\frac{\Li_{2}(\e(-x))}{(2\pi {\rm i})^2}\biggr) .
\]
Then
\[
 \frac{{\rm d}}{{\rm d}x}D(x)
  =
 -\frac{\log(2\sin(\pi x))}{2\pi} ,\qquad
 \frac{{\rm d}^2}{{\rm d}x^2}D(z)
  =
 -\frac{\cot(\pi x)}{2}
\]
has two stationary points at $1/6,5/6$ for $x\in[0,1]$ and achieves maximum and minimum values
\[
 \Im(\Li_{2}(\e(\pm1/6))
  =
 \pm1.0149\dots .
\]
Therefore, one can compute
\begin{align*}
& D(0.05)+0.035425
 =0.052608\ldots>0.051418\ldots =D(1/6),\\
& D(0.33)+0.035425
 =0.052854\ldots>0.051418\ldots =D(1/6) .
\end{align*}
Therefore, for $2x+2y\in[0.05,0.33]+\BZ$ or $x,y\in[0.67000,0.95000]$, we see that
\[
 D(2x+2y)-D(x)-D(y)>-0.035425 .
\]
For $x,y\in[0,0.025]$, we see that
\[
 D(2x+2y)-D(x)-D(y)+0.035425
 >
 0.035425-2D(0.025)
  =
 0.012735 .
\]
For $x\in[0.95,1]$ and $y\in[0,0.025]$, we see that
\[
 D(2x+2y)-D(x)-D(y)+0.035425
 >
 0.035425-D(0.1)-D(0.025)
  =
 0.00068087 .
\]
For $x\in[0.95,1]$, $y\in[0.025,0.075]$ and $2x+2y\leq 2.05$, we have
\[
 D(2x+2y)-D(x)-D(y)+0.035425
 >
 0.035425-D(0.0500)-D(0.025)
  =
 0.0068973 .
\]
For $x,y\in[0.49250,1]$, we see that
\begin{gather*}
 D(2x+2y)-D(x)-D(y)+0.035425
 >
 0.035425-D(1/6)-2D(0.0075)
 \\
 \qquad{} =
 3.5615\times 10^{-5} .
\end{gather*}
If $y\in[0.5,1]$ and $2x+2y\in[-0.15,0.5]+\BZ$, then
\begin{gather*}
 D(2x+2y)-D(x)-D(y)+0.035425
 >
 0.035425-D(1/6)
  =
 0.0097162 .\tag*{\qed}
\end{gather*}\renewcommand{\qed}{}
\end{proof}

Let four regions $A$, $B$, $C$, $D$ give the connected components of the complement of the set $U$ depicted in green.
\begin{Lemma}\label{lem:BCDsmall}
The integral of equation~\eqref{eq:int.of.sum} with $\alpha=-1$, $\beta=1$ over the region $B$, $C$, or $D$ are $o\bigl(-V_{\rho_6}r/(2\pi {\rm i})^2\bigr)$.
\end{Lemma}
\begin{proof}
We will prove this by illustrating that the contours can be explicitly deformed,
\begin{align*}
 V_{-1,1}(x,y) ={}&\frac{1}{(2\pi {\rm i})^2}\Li_2(\e(-2x-2y))-\frac{1}{(2\pi {\rm i})^2}\Li_2(\e(-x))-\frac{1}{(2\pi {\rm i})^2}\Li_2(\e(-y))\\
 & -\frac{1}{24}+\frac{1}{2}x^2+xy-\frac{1}{2}y^2-\frac{1}{2}x+\frac{1}{2}y.
\end{align*}
Therefore,
\[
\begin{aligned}
 \frac{\partial V_{-1,1}}{\partial x}
 & =
 \frac{2}{2\pi {\rm i}}\log(1-\e(-2x-2y))-\frac{1}{2\pi {\rm i}}\log(1-\e(-y))+x+y-\frac{1}{2} ,\\
 \frac{\partial V_{-1,1}}{\partial y}
 & =
 \frac{2}{2\pi {\rm i}}\log(1-\e(-2x-2y))-\frac{1}{2\pi {\rm i}}\log(1-\e(-x))+x-y+\frac{1}{2} .
\end{aligned}
\]
Firstly, we will check the conditions for region $B$.
We see that
\[
 \Im\biggl({\rm i}\frac{\partial V_{-1,1}}{\partial y}-{\rm i}\frac{\partial V_{-1,1}}{\partial x}\biggr)
  =
 \frac{1}{2\pi}\arg(1-\e(-x)))
 -\frac{1}{2\pi}\arg(1-\e(-y)))
 -2\Re(y)+1 .
\]
For $\Re(x),\Re(y)\in(0,1/2)$, $\Im(x)\leq0$, and $\Im(y)\geq0$, we find that
\[
\begin{aligned}
 0&\leq\frac{1}{2\pi}\arg(1-\e(-x))) ,\\
 -1&<-\frac{1}{2\pi}\arg(1-\e(-y)))
 -2\Re(y) .
\end{aligned}
\]
Therefore, for these conditions $x,y\in[0,1/2]$ deforming the contour in the direction the $(-{\rm i},{\rm i})$ direction will always increase the imaginary part of $V_{-1,1}(x,y)$, and therefore we see that the integral over the region $B$ can be deformed to an integral that is $o(-V_{\rho_6}r/(2\pi {\rm i})^2)$.

Secondly, for region $C$, we find for $\Re(x)\in[0.05,1/2)$ and $\Im(x)\geq0$
\[
 -\frac{1}{2\pi}\arg(1-\e(-x)))
 \geq
 -0.45
\]
Moreover, for $y\in[0.95,1)$ and $\Re(x)\in[0.05,1/2)$ and $\Im(x)\geq0$, we find
\[
 \frac{2}{2\pi}\arg(1-\e(-2x-2y)))>0,
\]
therefore, we find
\[
 \Im\biggl({\rm i}\frac{\partial V}{\partial x}\biggr)
  =
 \frac{2}{2\pi}\arg(1-\e(-2x-2y)))
 -\frac{1}{2\pi}\arg(1-\e(-x)))
 +\Re(x)+y-\frac{1}{2}
 \geq
 0.05.
\]
This implies that flowing a contour near the region $C$ by increasing the imaginary part of $x$ in this region will always increase the imaginary part of $V_{-1,1}$. Therefore, we see that we can deform the contour in this direction.
Therefore, for large enough $\Im(x)$, we see this will have positive imaginary part hence this contributes nothing to the asymptotics. Almost identical computation go through for the region $D$ with $x$ and $y$ interchanged.
\end{proof}

\begin{Lemma}\label{lem:Agood}
The integral of equation~\eqref{eq:int.of.sum} with $\alpha=-1$, $\beta=1$ over the region $A$ can be computed by stationary phase approximation at the critical point $(x_6,y_6)$.
\end{Lemma}
\begin{proof}
We deform the contour over region $A$ by steepest ascent with respect to the function $\Im(V_{-1,1}(x,y))$ fixing the real parts of $(x,y)$. In that region, there is only one critical point with real part contained in $A$ given by $(x_6,y_6)$. Therefore, on deforming the contour, we find that is passes through the critical point, where it attains it minimum value given by $\Im\bigl(-\frac{V_{\rho_6}}{(2\pi {\rm i})^2}\bigr)$. Therefore, we can apply stationary phase to the integral of equation~\eqref{eq:int.of.sum} at this point.
\end{proof}

\section{Details of state integral factorisations}
To deal with the state integrals, we use the notation
\[
 \Phi(z;\tau)
  =
 \Phi_{\mathsf{b}}\bigl({\rm i}z\mathsf{b}^{-1}+c_{\mathsf{b}}\bigr)
  =
 \frac{(q\e(z);q)_{\infty}}{(\e(z/\tau);\tq)_{\infty}} .
\]
Recall that for a positive integer $\ell$, the Eisenstein series are defined by
\be
\label{Eell}
\begin{aligned}
 E_{\ell}(q)
 =
 \frac{\zeta(1-\ell)}{2}+\sum_{s=1}^{\infty}s^{\ell-1}\frac{q^{s}}{1-q^s} ,
\end{aligned}
\ee
where $\zeta(s)$ is the Riemann zeta function, analytic for $|q|<1$ and extended to $|q|>1$ by the symmetry $E_{\ell}\bigl(q^{-1}\bigr) = -E_{\ell}(q)$. These Eisenstein series can be used with polylogarithms to describe an expansion of the infinite Pochhammer symbol.
\begin{Lemma}[\cite{GK:qser,GZ:qser}]
\label{lem.epexpansions}
We have identities in $\BQ[\![q,\epsilon]\!]$, where we use $\bigl(q^{-1};q^{-1}\bigr)_\infty=(q;q)_\infty^{-1}$ and $\bigl(x;q^{-1}\bigr)_\infty=(qx;q)_\infty^{-1}$,
\begin{align*}
&\frac{(qe^{\ve};q)_{\infty}}{(q;q)_{\infty}}\frac{1}{(qe^{\ve};q)_{m}} =
\frac{1}{(q;q)_{m}}\sqrt{\frac{-\ve}{1-\exp(\ve)}}
\exp\left(-\sum_{\ell=1}^{\infty}
 \left(E_{\ell}(q)-\sum_{n=1}^{m}\Li_{1-\ell}(q^{n})\right)
 \frac{\ve^\ell}{\ell!}\right) , \\
&\frac{\bigl(q^{-1}e^{\ve};q^{-1}\bigr)_{\infty}}{\bigl(q^{-1};q^{-1}\bigr)_{\infty}}
\frac{1}{\bigl(q^{-1}e^{\ve};q^{-1}\bigr)_{m}} \\
& \qquad =
\frac{1}{\bigl(q^{-1};q^{-1}\bigr)_{m}}\frac{-1}{\ve}\sqrt{\frac{-\ve}{1-\exp(\ve)}}
\exp\left(-\sum_{\ell=1}^{\infty}
 \left(E_{\ell}\bigl(q^{-1}\bigr)-\sum_{n=1}^{m}\Li_{1-\ell}(q^{-n})\right)
 \frac{\ve^{\ell}}{\ell!}\right) .
\end{align*}
\end{Lemma}

\subsection{At rationals}

We can take the formula for the WRT invariant
\[
 \WRT_{m}(q)
  =
 \sum_{0\leq\ell\leq k}(-1)^{k}q^{-\frac{1}{2}k(k+1)+\ell(\ell+1)+mk}\frac{(q;q)_{2k+1}}{(q;q)_{\ell}(q;q)_{k-\ell}} ,
\]
and naturally take the state integral
\[
 \int_{-{\rm i}\mathsf{b}(\BR-{\rm i}0)} \int_{-{\rm i}\mathsf{b}(\BR-{\rm i}0)}
 \frac{\e(-z_1(z_1+1+\tau)/2\tau+z_3(z_3+1+\tau)/\tau+z_1(m+m'/\tau))}{\Phi(z_3;\tau)^{-1}\Phi(z_1-z_3;\tau)^{-1}\Phi(2z_1+1+\tau;\tau)}{\rm d}z_1{\rm d}z_3 .
\]
We can rewrite this integral by first sending $z_1\mapsto -z_1$, $z_3\mapsto -z_3$ and using the inversion relation for the quantum dilogarithm~\cite[Appendix~A]{AK:I},
\begin{align}\label{eq:stateint.4112.v2}
& q^{-1-m}\tq^{1+m'}\\
&  \nonumber\times \iint_{(-{\rm i}\mathsf{b}(\BR-{\rm i}0))^2} \!\!\!\!\!\!\!\! \frac{\e(3z_1(z_1+1+\tau)/2\tau+z_3(z_3+1+\tau)/2\tau-z_1(m+1+(m'+1)/\tau))}{\Phi(z_3-z_1;\tau)^{-1}\Phi(2z_1;\tau)^{-1}\Phi(z_3;\tau)}{\rm d}z_1{\rm d}z_3.
\end{align}
Suppose that $\tau=N/M\in\BQ_{>0}$. Then we can use the fundamental lemma of~\cite[Lemma~2.1]{GK:rat} to factorise this state integral. However, we can use the Fourier transform formula for the quantum dilogarithm~\cite[Appendix~A]{AK:I} to expand the Faddeev dilogarithm in $\Phi(z_{3};\tau)$, and then again to contract the integral over $z_3$. With these steps, one can show that this integral is given by
\begin{align*}
&\calS_{m,m'}
=
 q^{-m}\tq^{m'} \\
& \times \int_{-{\rm i}\mathsf{b}(\BR-{\rm i}0)} \int_{-{\rm i}\mathsf{b}(\BR-{\rm i}0)} \frac{\e\bigl(z_1^2/\tau+z_1z_2/\tau+z_1(-m-m'/\tau)+z_2(1+1/\tau)\bigr)}{\Phi(z_2;\tau)^{-1}\Phi(2z_1;\tau)^{-1}\Phi(z_2-z_1;\tau)^{-1}}{\rm d}z_1{\rm d}z_2 .
\end{align*}
This integral expression has various benefits and in particular has a direction at infinity where the two dimensional contour can be pushed with vanishing contribution. However, we do not need this for the factorisation at rationals. We can explicitly use the fundamental lemma of~\cite[Lemma~2.1]{GK:rat} to factorise this state integral. This gives the following lemma analogous to Lemma~\ref{lem:7x7stateint} but for $\tau\in\BQ$. Lemma~\ref{lem:quadrels7x7} also holds for $\tau\in\BQ$ and together they give Theorem~\ref{thm:7x7qser} for $\tau\in\BQ$.
\begin{Lemma}
 We have the following identity for $q=\e(\tau)$, $\tq=\e(-1/\tau)$ with $\tau=N/M$ with coprime $M,N\in\BZ$,
 \be\label{eq:lem.bilin.Q}
 \calS_{m,m'}(\tau)
  =
 \sum_{j=1}^7\e\biggl(-\frac{\VC_{\rho_j}}{(2\pi {\rm i})^2MN}\biggr)\WRT^{(j)}_{m}(q)\WRT^{(j)}_{m'}\bigl(\tq^{-1}\bigr) .
 \ee
\end{Lemma}
\begin{proof}
Indeed, if we take $f_i(z_1,z_2)$ to be quotients of the integrand with $z_{i}\mapsto z_i+N$ in the numerators, then
\begin{align*}
 f_{1}(z_{1},z_{2})
 ={}&
 \frac{\e((z_1+N)^2/\tau+(z_1+N)z_2/\tau+(z_1+N)(-m-m'/\tau)+z_2(1+1/\tau))}{\Phi(z_2;\tau)^{-1}\Phi(2z_1+2N;\tau)^{-1}\Phi(z_2-z_1-N;\tau)^{-1}}\\
 & \times\frac{\Phi(z_2;\tau)^{-1}\Phi(2z_1;\tau)^{-1}\Phi(z_2-z_1;\tau)^{-1}}{\e\bigl(z_1^2/\tau+z_1z_2/\tau+z_1(-m-m'/\tau)+z_2(1+1/\tau)\bigr)}\\
={} &
 \frac{\e(2Mz_{1}+Mz_{2})(1-\e(-Mz_1+Mz_2))}{(1-\e(2Mz_1))^2}.
\end{align*}
Then
\begin{align*}
 f_{2}(z_{1},z_{2})
 ={}&
 \frac{\e(z_1^2/\tau+z_1(z_2+N)/\tau+z_1(-m-m'/\tau)+(z_2+N)(1+1/\tau))}{\Phi(z_2+N;\tau)^{-1}\Phi(2z_1;\tau)^{-1}\Phi(z_2+N-z_1;\tau)^{-1}}\\
 & \times\frac{\Phi(z_2;\tau)^{-1}\Phi(2z_1;\tau)^{-1}\Phi(z_2-z_1;\tau)^{-1}}{\e\bigl(z_1^2/\tau+z_1z_2/\tau+z_1(-m-m'/\tau)+z_2(1+1/\tau)\bigr)}\\
={} &
 \frac{\e(Mz_{1})}{(1-\e(Mz_2))(1-\e(-Mz_1+Mz_2))} .
\end{align*}
Notice that
\[
 f_{j}(z_1+kN,z_2+\ell N)
  =
 f_{j}(z_1,z_2) ,
\]
and therefore from a generalisation of~\cite[Lemma~2.1]{GK:rat} to multiple dimensions,
\begin{align}
 &\int_{-{\rm i}\mathsf{b}(\BR-{\rm i}0)}\int_{-{\rm i}\mathsf{b}(\BR-{\rm i}0)}\frac{\e\bigl(z_1^2/\tau+z_1z_2/\tau+z_1(-m-m'/\tau) +z_2(1+1/\tau)\bigr)}{\Phi(z_2;\tau)^{-1}\Phi(2z_1;\tau)^{-1}\Phi(z_2-z_1;\tau)^{-1}}{\rm d}z_1{\rm d}z_2\nonumber\\
 & =\biggl(\int_{-{\rm i}\mathsf{b}(\BR-{\rm i}0)}-\int_{N-{\rm i}\mathsf{b}(\BR-{\rm i}0)}\biggr)^2\frac{\e\bigl(z_1^2/\tau+z_1z_2/\tau +z_1(-m-m'/\tau)+z_2(1+1/\tau)\bigr){\rm d} z_1{\rm d} z_2}{\Phi(z_2;\tau)^{-1}\Phi(2z_1;\tau)^{-1}\Phi(z_2-z_1;\tau)^{-1} }\nonumber\\
 & \quad \times\frac{1}{\bigl(1-\frac{\e(2Mz_{1}+Mz_{2})(1-\e(-Mz_1+Mz_2))}{(1-\e(2Mz_1))^2}\bigr)\bigl(1-\frac{\e(Mz_{1})}{(1-\e(Mz_2))(1-\e(-Mz_1+Mz_2))}\bigr)}\nonumber\\
 & =\biggl(\int_{-{\rm i}\mathsf{b}(\BR-{\rm i}0)}-\int_{N-{\rm i}\mathsf{b}(\BR-{\rm i}0)}\biggr)^2 \frac{\e(z_1^2/\tau+z_1z_2/\tau+z_1(-m-m'/\tau)+z_2(1+1/\tau))}{\Phi(z_2;\tau)^{-1}\Phi(2z_1;\tau)^{-1}\Phi(z_2-z_1;\tau)^{-1}}\nonumber\\
 & \quad \times \frac{(1-\e(2Mz_1))^2(1-\e(Mz_2))(1-\e(-Mz_1+Mz_2)){\rm d}z_1{\rm d} z_2}{\bigl((1-\e(2Mz_1))^2-\e(2Mz_{1}+Mz_{2})(1-\e(-Mz_1+Mz_2))\bigr) }\nonumber\\
 & \quad \times \frac{1}{\bigl((1-\e(Mz_2))(1-\e(-Mz_1+Mz_2))-\e(Mz_{1})\bigr)}  .\label{eq:integral.4112.atrat}
\end{align}
We see that the zeros of the quantum dilogarithms now cancel with the zeros in the numerator and that the integrand has poles when
\begin{gather*}
 0 =(1-\e(Mz_2))(1-\e(-Mz_1+Mz_2))-\e(Mz_{1}) ,\\
 0 =(1-\e(2Mz_1))^2-\e(2Mz_{1}+Mz_{2})(1-\e(-Mz_1+Mz_2)) .
\end{gather*}
Therefore, considering the algebraic variety
\begin{gather*}
 0 =(1-X_2)\bigl(1-X_1^{-1}X_2\bigr)-X_1 ,\\
 0 =\bigl(1-X_1^2\bigr)^2-X_1^2X_2\bigl(1-X_1^{-1}X_2\bigr) ,
\end{gather*}
this can be solved exactly as notice that
\[
 \frac{1-X_2}{X_{1}} =\frac{X_1^2X_2}{\bigl(1-X_1^2\bigr)^2}
\]
and so
\[
 X_{2} =\frac{\bigl(1-X_1^2\bigr)^2}{\bigl(1-X_1^2\bigr)^2+X_{1}^3} ,
\]
and taking
\[
 \xi^7 - \xi^6 - 2 \xi^5 + 6 \xi^4 - 11 \xi^3 + 6 \xi^2 + 3 \xi - 1 =0 ,
\]
we have
\begin{gather*}
 X_{1,j} =-3+11\xi_j+20\xi_j^2-15\xi_j^3+6\xi_j^4-2\xi_j^5-4\xi_j^6 ,\\
 X_{2,j} =-9+19\xi_j+76\xi_j^2-52\xi_j^3+20\xi_j^4-4\xi_j^5-13\xi_j^6 .\\
\end{gather*}
There is one additional solution when $X_{1}=1$ and $X_{2}=0$, however $z_2=\infty$ for such a solution. Therefore, letting
\[
 x_{i,j,k} =\frac{1}{2\pi {\rm i}}\log(X_{i,j})+k
\]
notice that
\begin{align*}
 &(1-\e(2M(z_1+x_{1}/M)))^2-\e(2M(z_1+x_{1}/M)+M(z_2+x_2/M))\\
 &\qquad\quad\times(1-\e(-M(z_1+x_{1}/M)+M(z_2+x_2/M)))\\
 & \qquad =
 \bigl(1-X_{1}^2\e(2Mz_1)\bigr)^2-X_{1}^2X_{2}\e(2Mz_1+Mz_2)\bigl(1-X_{1}^{-1}X_{2}\e(-Mz_1+Mz_2)\bigr)\\
 & \qquad =
 \bigl(4X_1^4 + (-2X_2 - 4)X_1^2 + X_2^2X_1\bigr)2\pi {\rm i}Mz_1+\bigl(-X_2X_1^2 + 2X_2^2X_1\bigr)2\pi {\rm i}Mz_{2}+\cdots
\end{align*}
and
\begin{align*}
&\bigl((1-\e(M(z_2+x_2/M)))(1-\e(-M(z_1+x_{1}/M)+M(z_2+x_2/M)))-\e(M(z_1+x_{1}/M))\bigr)\\
& \qquad=\bigl((1-X_2\e(Mz_2))\bigl(1-X_1^{-1}X_2\e(-Mz_1+Mz_2)\bigr)-X_{1}\e(Mz_1)\bigr)\\
& \qquad =
\bigl(-X_1 + \bigl(-X_2^2 + X_2\bigr)X_1^{-1}\bigr)2\pi {\rm i}Mz_1+\bigl(-X_2 + \bigl(2X_2^2 - X_2\bigr)X_1^{-1}\bigr)2\pi {\rm i}Mz_{2}+\cdots .
\end{align*}
Then we can compute
\begin{align*}
 &\det
 \begin{pmatrix}
 4X_1^4 + (-2X_2 - 4)X_1^2 + X_2^2X_1 & -X_2X_1^2 + 2X_2^2X_1\\
 -X_1 + \bigl(-X_2^2 + X_2\bigr)X_1^{-1} & -X_2 + \bigl(2X_2^2 - X_2\bigr)X_1^{-1}
 \end{pmatrix}\\
 &\qquad  =
 -257+806\xi+947\xi^2-749\xi^3+331\xi^4-133\xi^5-213\xi^6
  =\Delta .
\end{align*}
The last ingredient we need is the quantum dilogarithm at rationals. For this, we use the cyclic dilogarithm
\[
 \CYCDS_{M}(x;q)
  =
 \prod_{j=1}^{M-1}\bigl(1-q^{j}x\bigr)^{\frac{j}{M}}
\]
and have the formula~\cite[Theorem 1.9]{GK:rat}
\[
 \Phi(z;\tau)^{-1}
 =\;
 \exp\biggl(\frac{1}{2\pi {\rm i}NM}\Li_{2}(\e(Mz))\biggr)(1-\e(Mz))^{\frac{z}{N}-1}\CYCDS_{N}\bigl(\e(z/\tau);\tq^{-1}\bigr)\CYCDS_{M}(\e(z);q) .
\]
Therefore, the double integral deformed near $(x_{1,j,k}/M,x_{2,j,\ell}/M)$ gives
\begin{align*}
 &\frac{\e\bigl(x_{1,j,k}^2/NM+x_{1,j,k}x_{2,j,\ell}/NM+x_{1,j,k}(-m/M-m'/N)+x_{2,j,\ell}(1/M+1/N)\bigr)}{\Phi(x_{2,j,\ell}/M+1+\tau;\tau)^{-1}\Phi(2x_{1,j,k}/M+1+\tau;\tau)^{-1}\Phi(x_{2,j,\ell}/M-x_{1,j,k}/M+1+\tau;\tau)^{-1}}\\
 & \quad \times\frac{\bigl(1-X_{1,j}^2\bigr)^2(1-X_{2,j})\bigl(1-X_{1,j}^{-1}X_{2,j}\bigr)}{M^2\Delta}\\
 & =
 \frac{\e(x_{1,j,k}^2/NM+x_{1,j,k}x_{2,j,\ell}/NM+x_{1,j,k}(-m/M-m'/N)+x_{2,j,\ell}(1/M+1/N))} {\bigl(1-X_{1,j}^2\bigr)^{-2}(1-X_{2,j})^{-1}\bigl(1-X_{1,j}^{-1}X_{2,j}\bigl)^{-1}M^2\Delta}\\
 & \quad \times\frac{\exp\Bigl(-\frac{1}{2\pi {\rm i}NM}\Li_{2}(X_{2,j})\Bigr)(1-X_{2,j})^{1-\frac{x_{2,j,\ell}}{NM}-\frac{1}{M}}}{\CYCDS_{N}\bigl(\e(x_{2,j,\ell}/N);\tq^{-1}\bigr) \CYCDS_{M}(q\e(x_{2,j,\ell}/M);q)}\\
 &\quad \times \frac{\exp\Bigl(-\frac{1}{2\pi {\rm i}NM}\Li_{2}(X_{1,j}^2)\Bigr)\bigl(1-X_{1,j}^2\bigr)^{1-\frac{2x_{1,j,k}}{NM}-\frac{1}{M}}}{\CYCDS_{N}\bigl(\e(2x_{1,j,k}/N);\tq^{-1}\bigr) \CYCDS_{M}(q\e(2x_{1,j,k}/M);q)}\\
 & \quad \times\frac{\exp\Bigl(-\frac{1}{2\pi {\rm i}NM}\Li_{2}(X_{1,j}^{-1}X_{2,j})\Bigr)\bigl(1-X_{1,j}^{-1}X_{2,j}\bigr)^{1-\frac{x_{2,j,\ell}-x_{1,j,k}}{NM} -\frac{1}{M}}}{\CYCDS_{N}\bigl(\e(x_{2,j,\ell}/N-x_{1,j,k}/N);\tq^{-1}\bigr)\CYCDS_{M}(q\e(x_{2,j,\ell}/M-x_{1,j,k}/M);q)} .
\end{align*}
Firstly, notice that
\begin{align*}
 &-\Li_{2}(X_{2,j})-\Li_{2}\bigl(X_{1,j}^2\bigr)-\Li_{2}\bigl(X_{1,j}^{-1}X_{2,j}\bigr)-2\pi {\rm i}x_{2,j,\ell}\log(1-X_{2,j})\\
 &-2(2\pi {\rm i})x_{1,j,k}\log\bigl(1-X_{1,j}^2\bigr)-2\pi {\rm i}(x_{2,j,\ell}-x_{1,j,k})\log\bigl(1-X_{1,j}^{-1}X_{2,j}\bigr)\\
 &+(2\pi {\rm i})^2x_{1,j,k}^2+(2\pi {\rm i})^2x_{1,j,k}x_{2,j,\ell}\\
 & =
 -\VC_{\rho_j}-4\pi^2\biggl(k^2+k\ell+a_jk+\frac{1}{24}\biggr),
\end{align*}
where $a_{j}=2,0,0,1,1,-1,1$ for $j=1,\dots,7$ or
\[
 2\pi {\rm i}a_{j}
  =
 -2\log\bigl(1-X_{1,j}^2\bigr)+\log\bigl(1-X_{1,j}^{-1}X_{2,j}\bigr)+2\log(X_1)+\log(X_2),
\]
where we note that for all $j$
\[
 0 =-\log(1-X_{2,j})-\log\bigl(1-X_{1,j}^{-1}X_{2,j}\bigr)+\log(X_{1,j}) .
\]
Therefore, we see that
\begin{align*}
 &\e\bigl(x_{1,j,k}^2/NM+x_{1,j,k}x_{2,j,\ell}/NM\bigr)
 \exp\biggl(-\frac{1}{2\pi {\rm i}NM}\Li_{2}(X_{2,j})\biggr)(1-X_{2,j})^{-\frac{x_{2,j,\ell}}{NM}}
 \\
 &\qquad\times\exp\biggl(-\frac{1}{2\pi {\rm i}NM}\Li_{2}\bigl(X_{1,j}^2\bigr)\biggr)\bigl(1-X_{1,j}^2\bigr)^{-\frac{2x_{1,j,k}}{NM}}
 \\
 &\qquad\times\exp\biggl(-\frac{1}{2\pi {\rm i}NM}\Li_{2}\bigl(X_{1,j}^{-1}X_{2,j}\bigr)\biggr)\bigl(1-X_{1,j}^{-1}X_{2,j}\bigr)^{-\frac{x_{2,j,\ell}-x_{1,j,k}}{NM}}\\
 & =
 \e\biggl(-\frac{\VC_{j}+\pi^2/6}{(2\pi {\rm i})^2NM}\biggr)
 \e\biggl(\frac{k^2+k\ell}{NM}\biggr)
 \e\biggl(\frac{ka_j+\ell\times0}{NM}\biggr) .
\end{align*}
With this, using the residue and Chinese remainder theorem and the cocycle of $\etaroots$, the final integral of equation~\eqref{eq:integral.4112.atrat} can then be shown to factorise into the bilinear sum of equation~\eqref{eq:lem.bilin.Q}.
\end{proof}

Including the trivial connection, we can apply the same methods to the state integral from equation~\eqref{eq:inhom.statint} written in a similar from to equation~\eqref{eq:stateint.4112.v2}. This has a similar outcome with some Appell--Lerch type versions of the functions $\WRT^{(j)}_{m}(q)$ of equation~\eqref{eq:per.fun} that appear in the first column of $\mathbf{Z}_{m}(q)^{-1}$,
\begin{align*}
 &K_n\WRT^{(j)}_{m}(q) =\frac{-{\rm i}}{M\etaroots(q)}\sqrt{\frac{\bigl(1-X_{1,j}^2\bigr)^2(1-X_{2,j})\bigl(1-X_{1,j}^{-1}X_{2,j}\bigr)}{\Delta_{\rho_j}}} \sum_{k,\ell\in\BZ/M\BZ}\frac{q^{(n-1)m}}{\bigl(1-q^{1-n+k}X_{1,j}^{\frac{1}{m}}\bigr)}\\
 &\qquad\times
 \frac{q^{k^2+k\ell-mk+\ell}X_{1,j}^{\frac{2k+\ell-m}{M}}X_{2,j}^{\frac{k+1}{M}}\prod_{i=0}^{M-1} \bigl(1-q^{i+1+\ell-k}X_{1,j}^{-1/M}X_{2,j}^{1/M}\bigr)^{-(i+1+\ell-k)/M-1/2}} {\prod_{i=0}^{M-1}\bigl(1-q^{i+1+\ell}X_{2,j}^{1/M}\bigr)^{(i+1+\ell)/M-1/2}\prod_{i=0}^{M-1}\bigl(1-q^{i+1+2k}X_{1,j}^{2/M}\bigr)^{(i+1+2k)/M-1/2}} ,
\end{align*}
so that
\[
 K_n\WRT^{(j)}_{m}(q)-K_n\WRT^{(j)}_{m-1}(q)
  =
 \WRT^{(j)}_{m}(q) .
\]
Then there are some additional residues at $z_1\in\tau\BZ$, which then factorise into $\WRT_{m}(q)$ and the function $1$ with a shift by $\tau^{1/2}$ compared to the other residues. This relies on the identity
\begin{gather*}
  \sum_{0\leq\ell\leq k\leq M/2}(-1)^{k}q^{-\frac{1}{2}k(k+1)+\ell(\ell+1)+mk}\frac{(q;q)_{2k+1}}{(q;q)_{\ell}(q;q)_{k-\ell}}\\
\qquad =
 -\sum_{M/2\leq k\leq \ell+k\leq M-1}(-1)^{\ell}q^{k(2k+1)+k\ell+\ell(\ell+1)/2-mk-m}\frac{(q;q)_{\ell+k}}{(q;q)_{2k-M}(q;q)_{\ell}} ,
\end{gather*}
which follows from the fact that $(q;q)_{n}=M\bigl(q^{-1};q^{-1}\bigr)_{M-1-n}^{-1}$.

\subsection[As $q$-series]{As $\boldsymbol{q}$-series}\label{sec:qserfac}

In this appendix, we will outline the proofs of Lemmas~\ref{lem:vanishingqseries},~\ref{lem:7x7stateint} and~\ref{lem:8x8stateint}.

\begin{proof}[Proofs of Lemma~\ref{lem:vanishingqseries}]
The proof of equation~\eqref{eq:z.qdif} follows from automatic $q$-holomorphic proofs described in~\cite{WilfZeil}. The vanishing of \smash{$Z^{(2)}_{m}(q)=0$} for $|q|>1$ follows from the fact that there is only one non-vanishing solution that is a power series in $q^{-m}$. Notice that \smash{$Z^{(1)}_m(q)$} is one such solution. Therefore, the coefficient of $q^{-m}$ determines \smash{$Z^{(2)}(q)$} in terms of \smash{$Z^{(1)}(q)$}. Therefore, for $|q|>1$ we want to compute
\[
 \sum_{j=0}^{\infty}
 \frac{q^{j}}{(q;q)_{j}^2}\biggl(-\frac{1}{2}+2E_{1}(q)-2\sum_{n=1}^{j}\frac{q^n}{1-q^n}\biggr) .
\]
Considering the family,
\[
 g_p(q)
  =
 \sum_{j=0}^{\infty}
 \frac{q^{j+jp}}{(q;q)_{j}^2}\biggl(-\frac{1}{2}-p+2E_{1}(q)-2\sum_{n=1}^{j}\frac{q^n}{1-q^n}\biggr) .
\]
This satisfies
\[
 2g_{p}(q)-(1-q^p)g_{p-1}(q)
  =
 g_{p+1}(q) .
\]
Therefore, for $p>0$ is determined by $p=0$. Moreover, we see that
\[
 \sum_{j=0}^{\infty}
 \frac{q^{j+jp}}{(q;q)_{j}^2}
  =
\begin{cases}
 q^{p^2/4}+\cdots & \text{for }p\to\infty\text{ even} ,\\
 2q^{(p^2-1)/4}+\cdots & \text{for }p\to\infty\text{ odd} .
 \end{cases}
 \]
satisfies the same equation. Therefore, a combination of these vanishes for all $p\geq0$. However,
\[
 g_{p}(q)
  =
 \sum_{j=0}^{\infty}
 \frac{q^{j+jp}}{(q;q)_{j}^2}\biggl(2j-p-2\sum_{n=j+1}^{n}\frac{q^{-n}}{1-q^{-n}}\biggr)
\]
implies
\[
 g_{p}(q)
  =
 \begin{cases}
 o\bigl(q^{p^2/4}\bigr) & \text{for }p\to\infty\text{ even} ,\\
 o\bigl(q^{(p^2-1)/4}\bigr) & \text{for }p\to\infty\text{ odd} .
 \end{cases}
\]
Therefore, we see that $g_{p}(q)=0$ for $p\geq0$.
\end{proof}

The proof of Lemma~\ref{lem:inhom.vanishingqseries} follows from similar methods.

\begin{proof}[Proofs of Lemma~\ref{lem:7x7stateint}]
From its formula in the upper half plane, we have
\begin{align*}
 \Phi(z+m\tau+n;\tau)
 & =
 \frac{(q\e(z);q)_{\infty}}{(\e(z/\tau);\tq)_{\infty}}\frac{1}{(q\e(z);q)_{m}}\frac{1}{\bigl(\tq^{-1}\e(z/\tau);\tq^{-1}\bigr)_{n}}\\
 & =
 \frac{(q\e(z);q)_{\infty}}{(\tq\e(z/\tau);\tq)_{\infty}}\frac{1}{1-\e(z/\tau)}\frac{1}{(q\e(z);q)_{m}}\frac{1}{\bigl(\tq^{-1}\e(z/\tau);\tq^{-1}\bigr)_{n}}.
\end{align*}
We will now factorise the state integral $\calS_{m,m'}$.
We will integrate over $z_2$ and then $z_1$ by pushing the contour to infinity and collecting residues. We have
\begin{align}
 &\int_{-{\rm i}\mathsf{b}(\BR-{\rm i}0)}\int_{-{\rm i}\mathsf{b}(\BR-{\rm i}0)}\Phi(z_2;\tau)\Phi(2z_1;\tau)\Phi(z_2-z_1;\tau)\nonumber\\
 &\quad\times
 \e\bigl(z_1^2/\tau+z_1z_2/\tau+z_1(-m-m'/\tau)+z_2(1+1/\tau)\bigr){\rm d}z_1{\rm d}z_2
  \nonumber\\
  =&
 \int_{-{\rm i}\mathsf{b}(\BR-{\rm i}0)}\! \sum_{j,j'=0}^{\infty}\!\mathop{\Res}\limits_{z_2=j\tau+j'}\!\!
 \frac{\e\bigl(2z_1^2/\tau\!+z_1z_2/\tau+z_1(1-m+(1-m')/\tau)+z_2(1+1/\tau)\bigr)} {\Phi(z_1+z_2;\tau)^{-1}\Phi(2z_1;\tau)^{-1}\Phi(z_2;\tau)^{-1}}{\rm d}z_1{\rm d}z_2
  \nonumber\\
 &+
 \int_{-{\rm i}\mathsf{b}(\BR-{\rm i}0)}\sum_{j,j'=0}^{\infty}\mathop{\Res}\limits_{z_2=j\tau+j'}
 \frac{\e\bigl(z_1^2/\tau+z_1z_2/\tau+z_1(-m-m'/\tau)+z_2(1+1/\tau)\bigr)}{\Phi(z_2;\tau)^{-1}\Phi(2z_1;\tau)^{-1}\Phi(z_2-z_1;\tau)^{-1}}{\rm d}z_1{\rm d}z_2
  \nonumber\\
  =
 &-\frac{\tau(q;q)_{\infty}}{(\tq;\tq)_{\infty}}\int_{-{\rm i}\mathsf{b}(\BR-{\rm i}0)}\sum_{j,j'=0}^{\infty} \frac{q^{j}\tq^{-j'}\e\bigl(2z_1^2/\tau+z_1((j+1-m)\tau+(j'+1-m'))/\tau\bigr)}{\Phi(z_1+j\tau+j';\tau)^{-1}\Phi(2z_1;\tau)^{-1}(q;q)_{j} \bigl(\tq^{-1};\tq^{-1}\bigr)_{j'}}{\rm d}z_1
  \nonumber\\
 &
 -\frac{\tau(q;q)_{\infty}}{(\tq;\tq)_{\infty}}\int_{-{\rm i}\mathsf{b}(\BR-{\rm i}0)}\sum_{j,j'=0}^{\infty} \frac{q^{j}\tq^{-j'}\e\bigl(z_1^2/\tau+z_1((j-m)\tau+(j'-m'))/\tau\bigr)}{\Phi(2z_1;\tau)^{-1}\Phi(j\tau+j'-z_1;\tau)^{-1}(q;q)_{j} \bigl(\tq^{-1};\tq^{-1}\bigr)_{j'}}{\rm d}z_1
\label{eq:fac.equ}\\
  =
 &-\frac{\tau(q;q)_{\infty}}{(\tq;\tq)_{\infty}}\sum_{k,j,k',j'=0}^{\infty}\mathop{\Res}\limits_{z_1=\frac{k}{2}\tau+\frac{k'}{2}}
 \frac{q^{j}\tq^{-j'}\e\bigl(2z_1^2/\tau+z_1((j+1-m)\tau+(j'+1-m'))/\tau\bigr)} {\Phi(z_1+j\tau+j';\tau)^{-1}\Phi(2z_1;\tau)^{-1}(q;q)_{j}\bigl(\tq^{-1};\tq^{-1}\bigr)_{j'}}{\rm d}z_1
 \nonumber\\
 &
 -\frac{\tau(q;q)_{\infty}}{(\tq;\tq)_{\infty}}\sum_{k,j,k',j'=0}^{\infty}\mathop{\Res}\limits_{z_1=\frac{k}{2}\tau+\frac{k'}{2}}
 \frac{q^{j}\tq^{-j'}\e\bigl(z_1^2/\tau+z_1((j-m)\tau+(j'-m'))/\tau\bigr)}{\Phi(2z_1;\tau)^{-1}\Phi(j\tau+j'-z_1;\tau)^{-1}(q;q)_{j} \bigl(\tq^{-1};\tq^{-1}\bigr)_{j'}}{\rm d}z_1.
 \nonumber
\end{align}

We can break up these sums into congruences. In particular, we get eight sums corresponding to $k,k'\equiv0,1\pmod{2}$ in the two sums in the last equality in equation~\eqref{eq:fac.equ}. The first of the sums with $k\equiv k'\equiv0\pmod{2}$ is given as follows:
\begin{align*}
 &\sum_{k,j,k',j'=0}^{\infty}\mathop{\Res}\limits_{z_1=k\tau+k'}\frac{\Phi(z_1+j\tau+j';\tau)\Phi(2z_1;\tau)}{(q;q)_{j}\bigl(\tq^{-1};\tq^{-1}\bigr)_{j'}}q^{j}\tq^{-j'}
 \\
 &\quad \times \e\bigl(2z_1^2/\tau+z_1((j+1-m)\tau+(j'+1-m))/\tau\bigr){\rm d}z_1\\
  =&
 \sum_{k,j,k',j'=0}^{\infty}\mathop{\Res}\limits_{z_1=0}\frac{\Phi(z_1+(k+j)\tau+(k'+j');\tau) \Phi(2z_1+2k\tau+2k';\tau)}{(q;q)_{j}\bigl(\tq^{-1};\tq^{-1}\bigr)_{j'}}\\
 &\quad\times\e(2(z_1+k\tau+k')^2/\tau+(z_1+k\tau+k')((j+1-m)\tau+(j'+1-m'))/\tau){\rm d}z_1\\
  =&
 \sum_{k,j,k',j'=0}^{\infty}\mathop{\Res}\limits_{z_1=0}
 \frac{(q\e(z_1);q)_{\infty}(q\e(2z_1);q)_{\infty}q^{2k^2+kj+j+k-mk}\e((4k+j+1-m)z_1)}{(q\e(z_1);q)_{k+j}(q\e(2z_1);q)_{2k}(q;q)_{j}}\\
 &\quad \times\frac{\e\bigl(8\pi {\rm i}E_2(q)z_1^2\bigr)(\e(z_1/\tau);\tq)_{\infty}^{-1}(\e(2z_1/\tau);\tq)_{\infty}^{-1}\tq^{-2k'^2-k'j'-j'-k'+m'k'}}{\bigl(\tq^{-1}\e(z_1/\tau);\tq^{-1}\bigr)_{k'+j'}
 \bigl(\tq^{-1}\e(2z_1/\tau);\tq^{-1}\bigr)_{2k'}}\\
 &\quad \times\frac{\e((4k'+j'+1-m')z_1/\tau)\e\bigl(8\pi {\rm i}E_2\bigl(\tq^{-1}\bigr)(z_1/\tau)^2\bigr)}{\bigl(\tq^{-1};\tq^{-1}\bigr)_{j'}}
 {\rm d}z_1,
\end{align*}
where we used the identity
\[
 8\pi {\rm i}E_2(q)-8\pi {\rm i}\tau^{-2}E_2(\tq) =\frac{2}{\tau} .
\]
This can now be expanded using Lemma~\ref{lem.epexpansions} to find
\begin{align*}
 &\sum_{k,j,k',j'=0}^{\infty}\mathop{\Res}\limits_{z_1=0}
 \frac{q^{2k^2+kj+j+k-mk}\e((4k+j+1-m)z_1)\e\bigl(8\pi iE_2(q)z_1^2\bigr)}{(q;q)_{k+j}(q;q)_{2k}(q;q)_{j}}\\
 & \times\sqrt{\frac{-(2\pi {\rm i}z_1)}{1-\e(z_1)}}\exp\left(-\sum_{\ell=1}^{\infty}\left(E_{\ell}(q)-\sum_{n=1}^{k+j}\Li_{1-\ell}(q^{n})\right)\frac{(2\pi {\rm i}z_1)^\ell}{\ell!}\right)\\
 &\quad\times\sqrt{\frac{-(4\pi {\rm i}z_1)}{1-\e(2z_1)}}\exp\left(-\sum_{\ell=1}^{\infty}\left(E_{\ell}(q)-\sum_{n=1}^{2k}\Li_{1-\ell}(q^{n})\right)\frac{(4\pi {\rm i} z_1)^\ell}{\ell!}\right)\\
 & \times\frac{\tq^{-2k'^2-k'j'-j'-k'+m'k'}\e((4k'+j'+1-m')z_1/\tau)\e\bigl(8\pi {\rm i} E_2\bigl(\tq^{-1}\bigr)(z_1/\tau)^2\bigr)}{\bigl(\tq^{-1};\tq^{-1}\bigr)_{k'+j'}\bigl(\tq^{-1};\tq^{-1}\bigr)_{2k'}\bigl(\tq^{-1};\tq^{-1}\bigr)_{j'}}\\
 & \times\frac{-1}{\sqrt{-(2\pi {\rm i}z_1/\tau)(1-\e(z_1/\tau))}}\exp\left(-\sum_{\ell=1}^{\infty}\left(E_{\ell}(\tq^{-1})-\sum_{n=1}^{k'+j'}\Li_{1-\ell}(\tq^{-n})\right)\frac{(2\pi {\rm i}z_1/\tau)^{\ell}}{\ell!}\right)\\
 & \times\frac{-1}{\sqrt{-(4\pi {\rm i}z_1/\tau)(1-\e(z_1/\tau))}}\exp\left(-\sum_{\ell=1}^{\infty}\left(E_{\ell}(\tq^{-1})-\sum_{n=1}^{2k'}\Li_{1-\ell}(\tq^{-n})\right)\frac{(4\pi {\rm i}z_1/\tau)^{\ell}}{\ell!}\right){\rm d}z_1 .
\end{align*}
Considering the $q$ terms expanded to second order, we have
\begin{align}
 &\sum_{k,j=0}^{\infty}
 \frac{q^{2k^2+kj+j+k}}{(q;q)_{k+j}(q;q)_{2k}(q;q)_{j}}\nonumber\\
 &\quad{}\times\exp\Bigg(\left(4k+j-m+\frac{1}{4}-3E_{1}(q)+\sum_{n=1}^{k+j}\Li_{0}(q^{n})+2\sum_{n=1}^{2k}\Li_{0}(q^{n})\right)(2\pi {\rm i}z_1)\nonumber\\
 &\qquad{}+\left(- \frac{5}{48}+\frac{3}{2}E_2(q)+\frac{1}{2}\sum_{n=1}^{k+j}\Li_{-1}(q^{n})+2\sum_{n=1}^{2k}\Li_{-1}(q^{n})\right)(2\pi {\rm i}z_1)^2+\cdots\Bigg).\label{eq:exp.tech.ap.1}
\end{align}
The second of the sums with $k\equiv k'\equiv0\pmod{2}$ is given as follows:
\begin{align*}
 &\sum_{k,j,k',j'=0}^{\infty}\mathop{\Res}\limits_{z_1=k\tau+k'}\frac{\Phi(2z_1;\tau)\Phi(j\tau+j'-z_1;\tau)}
 {(q;q)_{j}\bigl(\tq^{-1};\tq^{-1}\bigr)_{j'}}q^{j}\tq^{-j'}\\
 &\quad{}\times \e\bigl(z_1^2/\tau+z_1((j-m)\tau+(j'-m'))/\tau\bigr){\rm d}z_1\\
 & =
 \sum_{k,j,k',j'=0}^{\infty}\mathop{\Res}\limits_{z_1=0} \frac{\Phi(2z_1+2k\tau+2k';\tau)\Phi((j-k)\tau+(j'-k')-z_1;\tau)}{(q;q)_{j}\bigl(\tq^{-1};\tq^{-1}\bigr)_{j'}}q^{j}\tq^{-j'}\\
 &\quad{}\times\e\bigl((z_1+k\tau+k')^2/\tau+(z_1+k\tau+k')((j-m)\tau+(j'-m'))/\tau\bigr){\rm d}z_1\\
 & =
 \sum_{k,j,k',j'=0}^{\infty}\mathop{\Res}\limits_{z_1=0}\frac{(q\e(2z_1);q)_{\infty}(q\e(-z_1);q)_{\infty}q^{k^2+kj+j-mk}\e((2k+j+m)z_{1})\e\bigl(4\pi iE_2(q)z_1^2\bigr)}{(q\e(2z_1);q)_{2k}(q\e(-z_1);q)_{j-k}(q;q)_{j}}\\
 & \quad{}\times\frac{(\e(2z_1/\tau);\tq)_{\infty}^{-1}(\e(-z_1/\tau);\tq)_{\infty}^{-1}\tq^{-k'^2-k'j'-j'+m'k'}\e((2k'+j'-m')z_1/\tau)}{\bigl(\tq^{-1}\e(2z_1/\tau);\tq^{-1}\bigr)_{2k'}\bigl(\tq^{-1}\e(-z_1/\tau);\tq^{-1}\bigr)_{j'-k'} }\\
 &\quad\times\frac{\e\bigl(4\pi {\rm i}G_2(\tq^{-1})(z_1/\tau)^2\bigr)}{\bigl(\tq^{-1};\tq^{-1}\bigr)_{j'}} {\rm d}z_1.
\end{align*}
This can be expanded using Lemma~\ref{lem.epexpansions} and on a certain cone, we find
\begin{align*}
 &\sum_{k\leq j,k'\leq j'=0}^{\infty}\mathop{\Res}\limits_{z_1=0}\frac{q^{k^2+kj+j-mk}\e((2k+j-m)z_{1})\e\bigl(4\pi {\rm i} E_2(q)z_1^2\bigr)}{(q;q)_{2k}(q;q)_{j-k}(q;q)_{j}}\\
 & \times\sqrt{\frac{(2\pi {\rm i}z_1)}{1-\e(-z_1)}}\exp\left(-\sum_{\ell=1}^{\infty}\left(E_{\ell}(q)-\sum_{n=1}^{j-k}\Li_{1-\ell}(q^{n})\right)\frac{(-2\pi {\rm i}z_1)^\ell}{\ell!}\right)\\
 & \times\sqrt{\frac{-(4\pi {\rm i} z_1)}{1-\e(2z_1)}}\exp\left(-\sum_{\ell=1}^{\infty}\left(E_{\ell}(q)-\sum_{n=1}^{2k}\Li_{1-\ell}(q^{n})\right)\frac{(4\pi {\rm i}z_1)^\ell}{\ell!}\right)\\
 & \times\frac{\tq^{-k'^2-k'j'-j'+m'k'}\e((2k'+j'-m')z_1/\tau)\e\bigl(4\pi {\rm i} G_2\bigl(\tq^{-1}\bigr)(z_1/\tau)^2\bigr)}{\bigl(\tq^{-1};\tq^{-1}\bigr)_{2k'}\bigl(\tq^{-1};\tq^{-1}\bigr)_{j'-k'}\bigl(\tq^{-1};\tq^{-1}\bigr)_{j'}}\\
 & \times\frac{-1}{\sqrt{(2\pi {\rm i}z_1/\tau)(1-\e(-z_1/\tau))}}\exp \! \left(-\sum_{\ell=1}^{\infty}\left(E_{\ell}\bigl(\tq^{-1}\bigr)-\sum_{n=1}^{j'-k'}\Li_{1-\ell}\bigl(\tq^{-n}\bigr)\right) \!\frac{(-2\pi {\rm i}z_1/\tau)^{\ell}}{\ell!}\right)\\
 & \times\frac{-1}{\sqrt{-(4\pi {\rm i}z_1/\tau)(1-\e(2z_1/\tau))}}\exp\!\left(-\sum_{\ell=1}^{\infty}\left(E_{\ell}\bigl(\tq^{-1}\bigr)-\sum_{n=1}^{2k'}\Li_{1-\ell}\bigl(\tq^{-n}\bigr)\right)\! \frac{(4\pi {\rm i}z_1/\tau)^{\ell}}{\ell!}\right){\rm d}z_1 .
\end{align*}
Considering the $q$ terms and expanding to second order, we find
\begin{align}
 &\sum_{k\leq j=0}^{\infty}\frac{q^{k^2+kj+j-mk}}{(q;q)_{2k}(q;q)_{j-k}(q;q)_{j}}\nonumber\\
 &\quad\times\exp\Bigg(\left(2k+j-m-\frac{1}{4}-E_{1}(q)-\sum_{n=1}^{j-k}\Li_{0}(q^{n})+2\sum_{n=1}^{2k}\Li_{-1}(q^{n})\right)(2\pi {\rm i}z_{1})\nonumber\\
 &\quad\qquad+\left(-\frac{5}{48}-\frac{1}{2}E_2(q)+\frac{1}{2}\sum_{n=1}^{j-k}\Li_{-1}(q^{n})+2\sum_{n=1}^{2k}\Li_{-1}(q^{n})\right)(2\pi {\rm i}z_{1})^2+\cdots\Bigg)\nonumber\\
 & =
 \sum_{k,j=0}^{\infty}\frac{q^{2k^2+kj+j+k-mk}}{(q;q)_{2k}(q;q)_{j}(q;q)_{j+k}}\nonumber\\
 &\quad\times\exp\Bigg(\left(3k+j-m-\frac{1}{4}-E_{1}(q)-\sum_{n=1}^{j}\Li_{0}(q^{n})+2\sum_{n=1}^{2k}\Li_{-1}(q^{n})\right)(2\pi {\rm i}z_{1})\nonumber\\
 &\quad\qquad +\left(-\frac{5}{48}-\frac{1}{2}E_2(q)+\frac{1}{2}\sum_{n=1}^{j}\Li_{-1}(q^{n})+2\sum_{n=1}^{2k}\Li_{-1}(q^{n})\right)(2\pi {\rm i}z_{1})^2+\cdots\Bigg).\label{eq:exp.tech.ap.2}
\end{align}
On another cone, we find
\begin{align*}
 &\sum_{j<k,j'<k'=0}^{\infty}\mathop{\Res}\limits_{z_1=0}\frac{q^{k^2+kj+j-mk}\e((2k+j-m)z_{1})\e\bigl(4\pi {\rm i} E_2(q)z_1^2\bigr)\bigl(q^{-1};q^{-1}\bigr)_{k-j-1}}{(q;q)_{2k}(q;q)_{j}}\\
 &\quad\times-\sqrt{(2\pi {\rm i}z_1)(1-\e(-z_1))}\exp\left(\sum_{\ell=1}^{\infty}\left(E_{\ell}\bigl(q^{-1}\bigr)-\sum_{n=1}^{k-j-1}\Li_{1-\ell}(q^{-n})\right)\frac{(-2\pi {\rm i}z_1)^\ell}{\ell!}\right)\\
 & \times\sqrt{\frac{-(4\pi {\rm i}z_1)}{1-\e(2z_1)}}\exp\left(-\sum_{\ell=1}^{\infty}\left(E_{\ell}(q)-\sum_{n=1}^{2k}\Li_{1-\ell}(q^{n})\right)\frac{(4\pi {\rm i} z_1)^\ell}{\ell!}\right)\\
 & \times\frac{\tq^{-k'^2-k'j'-j'+m'k'}\e((2k'+j'-m')z_1/\tau)\e\bigl(4\pi {\rm i} G_2\bigl(\tq^{-1}\bigr)(z_1/\tau)^2\bigr)(\tq;\tq)_{k'-j'-1}}{\bigl(\tq^{-1};\tq^{-1}\bigr)_{2k'}\bigl(\tq^{-1};\tq^{-1}\bigr)_{j'}}\\
 & \times\sqrt{\frac{1-\e(z_1/\tau)}{-(2\pi {\rm i}z_1/\tau)}}\exp\left(\sum_{\ell=1}^{\infty}\left(E_{\ell}(\tq)-\sum_{n=1}^{k'-j'-1}\Li_{1-\ell}(\tq^{n})\right)\frac{(2\pi {\rm i}z_1/\tau)^\ell}{\ell!}\right)\\
 & \times\frac{-1}{\sqrt{-(4\pi {\rm i}z_1/\tau)(1-\e(2z_1/\tau))}}\exp\left(-\sum_{\ell=1}^{\infty}\left(E_{\ell}\bigl(\tq^{-1}\bigr)-\sum_{n=1}^{2k'}\Li_{1-\ell}(\tq^{-n})\right)\frac{(4\pi {\rm i} z_1/\tau)^{\ell}}{\ell!}\right){\rm d}z_1,
\end{align*}
and so taking the $q$ terms and expanding to first order, we find that
\begin{align}
 &-(2\pi {\rm i}z_1)\sum_{j<k=0}^{\infty}\mathop{\Res}\limits_{z_1=0}\frac{q^{k^2+kj+j-mk}\bigl(q^{-1};q^{-1}\bigr)_{k-j-1}}{(q;q)_{2k}(q;q)_{j}}\label{eq:exp.tech.ap.3}\\
 &\times\exp\left(\left(2k+j-m-3/4-E_1(q)+2\sum_{n=1}^{2k}\Li_{0}(q^{n})+\sum_{n=1}^{k-j-1}\Li_{0}(q^{-n})\right)(2\pi {\rm i}z_{1})+\cdots\right).\nonumber
\end{align}
Therefore, the contribution of the sums with $k\equiv k'\equiv0\pmod{2}$ to the residue is given by
\[
 Z_{m}^{(1)}(q)Z_{m'}^{(2)}\bigl(\tq^{-1}\bigr)+\tau Z_{m}^{(2)}(q)Z_{m'}^{(1)}\bigl(\tq^{-1}\bigr) .
\]
When we restrict to congruences when $k$ or $k'$ are odd we find simple poles and the computation is much simpler. There it is not hard to see the bilinear combinations of the series $Z^{(4)}_m(q),\dots,Z^{(9)}_m(q)$.
\end{proof}

\begin{proof}[Proofs of Lemma~\ref{lem:8x8stateint}]
To factorise the state integral of equation~\eqref{eq:inhom.statint}, we apply the same kind of computations to find Appell--Lerch type sums in $\tq^{-1}$ to find Lemma~\ref{lem:8x8stateint}. This lemma is written in terms of the following explicit $q$-series.

For $n>1$, we find
\begin{align}
& K^{(1)}_{m,n}(q)
 =
 (q;q)_{\infty}^2\sum_{k,j=0}^{\infty}\frac{q^{k(2k+1)+jk+j-mk+(n-1)m}}{(q;q)_{j}(q;q)_{2k}(q;q)_{k+j}\bigl(1-q^{1-n+k}\bigr)},
 \nonumber\\
& K^{(2)}_{m,n}(q)
 =
 (q;q)_{\infty}^2\sum_{k,j=0}^{\infty}
 \frac{q^{k(2k+1)+jk+j-mk+(n-1)m}}{(q;q)_{j}(q;q)_{2k}(q;q)_{k+j}\bigl(1-q^{1-n+k}\bigr)}\nonumber\\
 &\hphantom{K^{(2)}_{m,n}(q)=}{}\times
 \left(-k-\frac{1}{2}+2E_{1}(q) -\sum_{\ell=1}^{j}\frac{q^\ell}{1-q^\ell}-\sum_{\ell=1}^{k+j}\frac{q^\ell}{1-q^\ell}\right)\nonumber\\
 &\hphantom{K^{(2)}_{m,n}(q)=}{} +(q;q)_{\infty}^2\sum_{k=0}^{\infty}\sum_{j=-k}^{-1}
 \frac{q^{k(2k+1)+jk+j-mk+(n-1)m}\bigl(q^{-1};q^{-1}\bigr)_{-j-1}}{(q;q)_{2k}(q;q)_{k+j}\bigl(1-q^{1-n+k}\bigr)},\nonumber\\
 & K^{(4)}_{m,n}(q)
=
 (q;q)_{\infty}\bigl(q^{3/2};q\bigr)_{\infty}\sum_{k,j=0}^{\infty}\frac{q^{(2k+1)(2k+2)/2+(k+1/2)j-m(k+1/2)+j+(n-1)m}} {(q;q)_{j}(q;q)_{2k+1}\bigl(q^{3/2};q\bigr)_{k+j}\bigl(1-q^{3/2-n+k}\bigr)},\nonumber\\
 & K^{(5)}_{m,n}(q)
=
 (q;q)_{\infty}(-q;q)_{\infty}\sum_{k,j=0}^{\infty}(-1)^{j+m}\frac{q^{k(2k+1)+jk+j-mk+(n-1)m}}{(q;q)_{j}(q;q)_{2k}(-q;q)_{k+j}\bigl(1+q^{1-n+k}\bigr)},\nonumber\\
& K^{(6)}_{m,n}(q)
 =
 (q;q)_{\infty}\bigl(-q^{3/2};q\bigr)_{\infty}\sum_{k,j=0}^{\infty}(-1)^{j+m}\frac{q^{(2k+1)(2k+2)/2+(k+1/2)j-m(k+1/2)+j+(n-1)m}}{(q;q)_{j} (q;q)_{2k+1}\bigl(-q^{3/2};q\bigr)_{k+j}\bigl(1+q^{3/2-n+k}\bigr)},\nonumber\\
& K^{(7)}_{m,n}(q)
 =
 (q;q)_{\infty}\bigl(q^{1/2};q\bigr)_{\infty}\nonumber\\
&\hphantom{K^{(7)}_{m,n}(q)=}{}\times
\sum_{k,j=0}^{\infty}\frac{q^{(2k+1)(2k+2)/2+(j-k-1/2)(k+1/2)+(j-k-1/2)-m(k+1/2)+(n-1)m}} {\bigl(q^{1/2};q\bigr)_{j-k}(q;q)_{2k+1}(q;q)_{j}\bigl(1-q^{3/2-n+k}\bigr)},\nonumber\\
& K^{(8)}_{m,n}(q)
 =
 (q;q)_{\infty}(-q;q)_{\infty}\sum_{k,j=0}^{\infty}(-1)^{j+m}\frac{q^{k(2k+1)+(j-k)k+(j-k)-mk+(n-1)m}} {(-q;q)_{j-k}(q;q)_{2k}(q;q)_{j}\bigl(1+q^{1-n+k}\bigr)},\nonumber\\
& K^{(9)}_{m,n}(q)
 =
 (q;q)_{\infty}\bigl(-q^{1/2};q\bigr)_{\infty} \label{eq:theqser.ler}\\
 &\phantom{K^{(9)}_{m,n}(q)=}{}\times
 \sum_{k,j=0}^{\infty}(-1)^{j+m} \frac{q^{(2k+1)(2k+2)/2+(j-k-1/2)(k+1/2)+(j-k-1/2)-m(k+1/2)+(n-1)m}}{\bigl(-q^{1/2};q\bigr)_{j-k}(q;q)_{2k+1}(q;q)_{j}\bigl(1+q^{3/2-n+k}\bigr)},\nonumber
\end{align}
and for $n=1$, we have
\begin{gather*}
 K^{(1)}_{m,1}(q)
 =
 (q;q)_{\infty}^2\sum_{k=1}^{\infty}\sum_{j=0}^{\infty}\frac{q^{k(2k+1)+jk+j-mk}}{(q;q)_{j}(q;q)_{2k}(q;q)_{k+j}(1-q^{k})}\\
 \hphantom{K^{(1)}_{m,1}(q)=}{}
 +(q;q)_{\infty}^2\sum_{j=0}^{\infty}\frac{q^{j}}{(q;q)_{j}^2}\left(m-\frac{1}{2}+2E_{1}(q)-2\sum_{n=1}^{j}\frac{q^n}{1-q^n}\right),
 \\
 K^{(2)}_{m,1}(q)
 =
 (q;q)_{\infty}^2\sum_{k=1}^{\infty}\sum_{j=0}^{\infty}
 \frac{q^{k(2k+1)+jk+j-mk}}{(q;q)_{j}(q;q)_{2k}(q;q)_{k+j}(1-q^{k})}\\
 \hphantom{K^{(2)}_{m,1}(q)=}{} \times
 \left(-k-\frac{1}{2}+2E_{1}(q)-\sum_{\ell=1}^{j}\frac{q^\ell}{1-q^\ell}-\sum_{\ell=1}^{k+j}\frac{q^\ell}{1-q^\ell}\right)\\
 \hphantom{K^{(2)}_{m,1}(q)=}{}+(q;q)_{\infty}^2\sum_{k=1}^{\infty}\sum_{j=-k}^{-1}
 \frac{q^{k(2k+1)+jk+j-mk}\bigl(q^{-1};q^{-1}\bigr)_{-j-1}}{(q;q)_{2k}(q;q)_{k+j}(1-q^{k})}\\
 \hphantom{K^{(2)}_{m,1}(q)=}{}+(q;q)_{\infty}^2\sum_{j=0}^{\infty}
 \frac{q^{j}}{(q;q)_{j}^2}\biggl(\!\bigl(j-m-2E_1(q)\bigr)
 \!\biggl(\frac{1}{2}-2E_1(q)+2\sum_{\ell=1}^{j}\frac{q^\ell}{1-q^\ell}\biggr)\!+2E_2(q)\!\biggr).
\end{gather*}
Using these functions and elementary $q$-holonomic methods, one can prove that for $j=1,2,4,5,\allowbreak 6,7,8,9$, we have
\begin{align}
 \frac{1}{2q(1-q)}L^{(j)}_{m}(q)
 ={}&q^{3}K^{(j)}_{-m-6,3}(q)
 +q^{2}K^{(j)}_{-m-5,3}(q) + \left(\frac{1}{2} + q\right)K^{(j)}_{-m-5,2}(q)\nonumber\\
 & +\bigl(-q^{-1} - 1\bigr)K^{(j)}_{-m-4,3}(q) + K^{(j)}_{-m-4,2}(q) + \frac{1}{2}q^{-1}K^{(j)}_{-m-4,1}(q)\nonumber\\
 & +\bigl(-q^{-2} - q^{-1}\bigr)K^{(j)}_{-m-3,3}(q) + \frac{1}{2}q^{-2}K^{(j)}_{-m-3,2}(q)\nonumber\\
 & +q^{-4}K^{(j)}_{-m-2,3}(q) + q^{-3}K^{(j)}_{-m-2,2}(q)
 +q^{-5}K^{(j)}_{-m-1,3}(q) .\label{eq:LinK.gen}
\end{align}
Then to finish the proof consider
\begin{align*}
 &\frac{2-2\tq}{q^{m}\tq^{2-m'}}\int_{\BR-{\rm i} 0}\int_{\BR-{\rm i}0}\Phi_{\mathsf{b}}(x_2+c_{\mathsf{b}})\Phi_{\mathsf{b}}(2x_1+c_{\mathsf{b}})\Phi_{\mathsf{b}}(x_2-x_1+c_{\mathsf{b}})\\
 &\qquad{}\times\e\bigl(-x_1^2-x_1x_2+ix_1(m\mathsf{b}-m'\mathsf{b}^{-1})-{\rm i} x_2\bigl(\mathsf{b}+\mathsf{b}^{-1}\bigr)\bigr)\frac{1}{1-\e\bigl(-{\rm i}x_1\mathsf{b}^{-1}\bigr)\tq^n} {\rm d}x_1{\rm d}x_2 .
\end{align*}
Has a similar residue formula as that in equation~\eqref{eq:fac.equ}. Then using the expansions in equation~\eqref{eq:exp.tech.ap.1}, equation~\eqref{eq:exp.tech.ap.2}, and equation~\eqref{eq:exp.tech.ap.3} we see that the corresponding congruence sums with $k\equiv k'\equiv0\pmod{2}$ leads to a contribution
\[
 Z_{m}^{(1)}(q)K_{m',n}^{(2)}\bigl(\tq^{-1}\bigr)+\tau Z_{m}^{(2)}(q)K_{m',n}^{(1)}\bigl(\tq^{-1}\bigr)+\tau^2 Z_{m}^{(3)}(q)K_{m',n}^{(0)}\bigl(\tq^{-1}\bigr) .
\]
Then the odd congruences give the other contributions
\[
\begin{aligned}
&Z_{m}^{(4)}(q)K_{m',n}^{(5)}\bigl(\tq^{-1}\bigr)+Z_{m}^{(5)}(q)K_{m',n}^{(4)}\bigl(\tq^{-1}\bigr)+Z_{m}^{(6)}(q)K_{m',n}^{(6)}\bigl(\tq^{-1}\bigr)\\
&\qquad{}+Z_{m}^{(7)}(q)K_{m',n}^{(8)}\bigl(\tq^{-1}\bigr)+Z_{m}^{(8)}(q)K_{m',n}^{(7)}\bigl(\tq^{-1}\bigr)+Z_{m}^{(9)}(q)K_{m',n}^{(9)}\bigl(\tq^{-1}\bigr) .
\end{aligned}
\]
Combining these into the full integral in equation~\eqref{eq:inhom.statint} giving $\calS^{(0)}_{m,m'}(\tau)$ and using equation~\eqref{eq:LinK.gen} completes the proof.
\end{proof}

\section{Numerical data}\label{app:numerical}

\subsection{Coefficients of asymptotic series}

Let
\[
\delta_{\rho}^{3k}D_{k}A_{k}^{\rho}=v_k\cdot\bigl(1,\xi_{\rho},\dots,\xi_{\rho}^{6}\bigr),
\]
 where $A_{k}^{\rho}$ are the coefficients of the formal series $\Phi^{(\rho)}(h)$ and $D_{k}$ is the universal denominator~\cite[Section~9]{GZ:RQMOD} or $\mathsf{A144618}$ of~\cite{oeis}. Compare this with equation~\eqref{eq:first.two.coeffs}. Then we have the first ten coefficients given by
\begin{align*}
& v_{1} =\text{\small $\begin{pmatrix}1497746\\ 1345119\\ -3675733\\ 2082815\\ -839488\\ -283405\\ 383432\end{pmatrix}$},\quad
v_{2} =\text{\small $\begin{pmatrix}3014838521575\\ 2732414541176\\ -7414786842283\\ 4197826806919\\ -1690529009777\\ -574198051621\\ 771765277669\end{pmatrix}$},\quad
  v_{3} =\text{\small $\begin{pmatrix}32075969439962558848\\ 29072051320986984192\\ -78890163112729764262\\ 44663022669399364904\\ -17986393087279219908\\ -6109213574944197629\\ 8211237042735127719\end{pmatrix}$},\\
& v_{4}  =\text{\small $\begin{pmatrix}98903732906139198096381046\\ 89641412486672113132053947\\ -243251651355458316156314530\\ 137714936589966968925726953\\ -55459634015776148007047322\\ -18837283785124500422973646\\ 25318706488062016968836637\end{pmatrix}$},\quad
  v_{5} =\text{\small $\begin{pmatrix}373656849964771673102447848336968\\ 338663940408450448690705586196244\\ -919001166686457533093349035199303\\ 520285008027629237579984234248816\\ -209525682773484396469829882941570\\ -71166981528115309886518232968370\\ 95653701284221268807576227329105\end{pmatrix}$},\\
&v_{6} =\text{\small $\begin{pmatrix}21989671118201887803945372073321048495066\\ 19930341623289247502670115082673523380863\\ -54083133803870016705986580318428675191313\\ 30618724681985189739403786413017035235444\\ -12330567084572849612900391145280625338946\\ -4188170291786657171757720239835038734160\\ 5629211488385838304799798635163554877439\end{pmatrix}$},\\
&v_{7} =\text{\small $\begin{pmatrix}-1054030545825412712812408078160552195229863984\\ -955320738849619688293413314788684756250452751\\ 2592365967496492262855439169469379443787684774\\ -1467646829075563688463535229534257371862182237\\ 591040870262700339422927393396401704097747915\\ 200751498052648379128268740205002933658138908\\ -269824902145113475694954495714866370536391072\end{pmatrix}$},\\
&v_{8} =\text{\small $\begin{pmatrix}785068551753725112218001620111973623952147096632968703\\ 711546996317316237545794569697964227777285521066597855\\ -1930859597739726343641718529915337408555715366354050216\\ 1093140398209113767916417765637919618616403605757376033\\ -440222156636833263384654075794264680557120599871459754\\ -149524782239719590802932762607629951810005663265867143\\ 200972396856152746359957812805528788099931903889495638\end{pmatrix}$},\\
&v_{9} =\text{\small $\begin{pmatrix}-546829073043613191824846643992763184725397522350496553360578870\\ -495618610061486366105609510787883683088888496176511424880507566\\ 1344914608604265087316615193816366867727891325200855048212591954\\ -761412426117321614544155613318478678354942966429698151498485708\\ 306630896511181763733968906570771016524016131993573668550481490\\ 104149501195206951752312520707872472447564279840830737141285971\\ -139984653868389669281650143397062122245612674367744717538309452\end{pmatrix}$},\\
&v_{10} = \text{\small $\left( \begin{matrix}55773889227532732828184524484609228876938777009880981456945186409638709\\ 50550672631241763416652769995927583211210706476150220837684914984478733\\ -137174707963638567039716829430329426944220354934166406966672359632621335\\ 77660341053862422297075354814431049552802650479442395163628006479192885\\ -31274850769300096685615679903427192732422153751016963837955145275755771\\ -10622739406361023735394438986275519058827807031940355968987688901516933\\ 14277749599073307932548287383306128369126364441142578078536796870888011\end{matrix} \right)$}.
\end{align*}
The one hundredth coefficient is given by $v_{100}$ with integral coefficients
\begin{align*}
&(v_{100})_1
=
\text{\scriptsize  $\begin{array}{@{}l@{}}
-114556048165446096903495795730947849708437279660658295476857439604114304994175860750170\\
2226219399306807433289133146923592607781713175774346292294872664755350832771967051497175\\
1334313782984397480690877472331622677387425915948755657272716966854004507060650209013241\\
6569706356532192405608071940707397415719345482190248035899758781995872272394920146262330\\
1698548060290220636218128535516038952334985685437490687755438347915422766743031977132921\\
3292609804728940391733411498064054755102595914134649830076930702285202226549672024572442\\
2741333959931616255085383320293741396945674084219989034366072333443935961567195438650179\\
5926711644098970714058781335872578746020301408937843017961016290954258605034585621333384\\
7322147861563730543658684156317706213094956485020905848674804554490149990732271461422778\\
2281287274267428053381108128078094921230370426655504834857948597417043431844274571507595\\
077515850057994343342
\end{array}$},\\[1mm]
&(v_{100})_2=
\text{\scriptsize  $\begin{array}{@{}l@{}}
-103827890952986718569282719960333743848174227847715178535429357211487930754122543256928\\
6364771358135387050255871031142739581927289814480155391971798017173295004155625332050227\\
4552419760012000843565268420484614969812367355718817073566254303660362148101472318443023\\
7150076080342818946875553065491550558002582449943712909886590635148221826112420724177028\\
0730359573816938240186730714475732140488169600119347794067930446705057263186951018406293\\
5294030444393780222067546014674431986945472501721336924732896437921205013780082295797140\\
1420996854491250153954071321269813883064162338546674074075232238975129219946201764336834\\
6043667303157009459477632965305685142276803043428619047097428173800728222713337389236249\\
3915728398979584622063083634277443950368540162602226327517742165963822258371902591576445\\
3326882617264326027274951309060728350395424921219472167658345935987680321999223347050252\\
920811379128326403892\end{array}$},\\[1mm]
&(v_{100})_3=
\text{\scriptsize  $\begin{array}{@{}l@{}}
2817481920340295055281346730690472644380317933071702198598981495233282399086191435667449\\
6318545087786559340588181231963365406120047446514342724006236372067633836753208996888304\\
7784544580349968001472287099983364797488008832833025087828779817376659479276529487964064\\
9799525886383878320853736065860615741140404261055267296172055845675881830691974456608354\\
9748482145914988717471723321568151130409397730543067106720446514892664882995860806845723\\
4843124130079766196008201885611511607802526887254002686408169780007566008120907107167410\\
1792661774536529007063128316059331914983920449666700983429133272838450621725348999813875\\
4427982048981672947342450227940868000582847370423943946310831419318625881988514349327639\\
3750477411195956944990997115743753860437766714701707838664070755329703657474544852901710\\
1459110753136357107434769089560412327126660393670228149990349172522970297062571687376935\\
84071673080459868758\end{array}$},\\[1mm]
&(v_{100})_4=
\text{\scriptsize  $\begin{array}{@{}l@{}}
-159509438798818852740403948954708034054610961555102441028658317201518481964882236020489\\
4526799952567533755570836430393792790381924179609668417088998737430330219980796803695379\\
3312110356011537472775765180374738843769484794042634477790332491319457344242422789257153\\
0202729885176121344727985617671903426889715715344945581792699129500328034897328362210132\\
7517218122563497996549327492790306206127518554631425135318382559197078897066162615676358\\
1770094267848826178566510687812386156251320318287628107064688818880045911108107159895752\\
8030213030069261610567623973995227099703139785973985230453804758568537620851861601286735\\
9763828390894194681264188506211920648451232529446588028801054008025311995957596039694674\\
7094507748532691760075874182877531587172664832988050640411901221776807894698254871948954\\
2051609414231763114152221244838930246570222277129894592099968500846202161763328423193297\\
355797112560595448063\end{array}$},\\[1mm]
&(v_{100})_5=
\text{\scriptsize  $\begin{array}{@{}l@{}}
6423656949005579635602380044834745805954008050422449757542181237662231935350381351794282\\
6345530520200575928837990407786626073345472662850049121938319972295907828869161004930804\\
7993144584499473846590267810548893831624583994913266780345038302207654608206078395361455\\
7040407260356741811074829257916925751026364921789888314066490354946561780131889416228616\\
0129263647153160879755569142812948767636047073193357707155704484393729503906526211820999\\
9030240478071441375193637129884062811997809370516862524875910273528679367653302910891415\\
1009779756107707217696125995167983346647078858750106017345906644337215675333990670042745\\
9925628867181510009160177921849999665405522221961363954722660010756647216901550860112516\\
1431002630893268513862917464594248805266993199441672822494617884430439806978212793583385\\
8904767686962286325345777540345688304253609632969892351149049178284066587117846175566832\\
1869460062345737528\end{array}$},
\\[1mm]
&(v_{100})_6=
\text{\scriptsize  $\begin{array}{@{}l@{}}
2181843626001527857043303377779018779842300046160536046727256908658290916944915026680454\\
7635664935136342238494294282216288689263298203300473168339347861015679674659575853857484\\
3272758027896194993453060494080374508914094672090391820927332158128265301293326997247408\\
4995805087693039029285239020092406352948611944824187109825858265021259418523302156261650\\
8092858553305488071316604868135326927393070682672073489736085898663010188719273481024248\\
9400372295729818032649504717891526934095477264744797523231335241235254598394417984278017\\
7750127603963651262570092473539686740530733191906639274388446202472691228579371886844066\\
1886682276538044627554398689381929689042197266102879028526074642438047560942764074898339\\
1291014038520316914301735361963268483007842819634202309355619045861602386968733207276110\\
4482184907706156431872570848820846282918582897113159425324906185999289823652080974172674\\
7058264455132676090\end{array}$},\\[1mm]
&(v_{100})_7=
\text{\scriptsize  $\begin{array}{@{}l@{}}
-293255965007763301055959616729773341327807076738201811172617652334185666107000427163406\\
1347194485867193959491892898983163471167273115485473516378390204101418762987067995590870\\
3021890755912110839030136531006257006903513073402976123712152658902037753367293325813950\\
4816431669392563492098383670641296976783763858683346978194946508598226243012839351583086\\
3419700678522902019923577150756869056361690035488669558001977155588085780209293795421552\\
7431786286302346537925124448672877893689099626289350461034118531699511804908446632109387\\
2067515765901114529461240564985381816167395965900356210130711747940097796984287932290626\\
9172128562958007775380904832673308287880950992069799640010522494309497291173909385099974\\
5403459339567654527745136673145209530674651373846980239803787758365037797734419230567386\\
8776417408486199172676339218746795467065463496428945361248282032768001859399532612619115\\
60556187926242200365
\end{array}$}.
\end{align*}

\subsection{Generating series of Stokes constants}\label{app:stokes}

The conjectural generating series of Stokes constants for the series $\Phi^{(\rho_{0})}$ of $4_1(-1,2)$ are given
\begin{align*}
\mathsf{S}_{+}^{(\rho_{1},\rho_{0})}(q)
 ={}&-q^3 - 3q^4 - 4q^5 - 4q^6 - 2q^7 + 3q^8 + 9q^9 + 18q^{10} + 27q^{11} + 36q^{12} + 41q^{13}\\
& + 38q^{14} + 29q^{15} + 3q^{16} - 35q^{17} - 92q^{18} - 163q^{19} - 252q^{20} + \cdots , \\
\mathsf{S}_{+}^{(\rho_{2},\rho_{0})}(q)
 ={}&-1 + q^3 + 3q^4 + q^5 - 5q^7 - 9q^8 - 15q^9 - 20q^{10} - 26q^{11} - 24q^{12} - 19q^{13}\\
&  - 6q^{14} + 20q^{15} + 53q^{16} + 101q^{17} + 155q^{18} + 216q^{19} + 274q^{20} + \cdots , \\
\mathsf{S}_{+}^{(\rho_{3},\rho_{0})}(q)
 ={}&-q + q^2 - q^4 - q^5 - 2q^6 - q^7 - q^8 + 4q^9 + 6q^{10} + 12q^{11} + 15q^{12} + 21q^{13}\\
&  + 18q^{14} + 19q^{15} + 8q^{16} - 5q^{17} - 31q^{18} - 59q^{19} - 99q^{20} + \cdots , \\
\mathsf{S}_{+}^{(\rho_{4},\rho_{0})}(q)
 ={}&-q - q^2 + q^3 + q^4 + 4q^5 + 3q^6 + 3q^7 - 8q^9 - 14q^{10} - 25q^{11} - 35q^{12}\\
& - 45q^{13} - 45q^{14} - 46q^{15} - 26q^{16} + 50q^{18} + 113q^{19} + 200q^{20} + \cdots , \\
\mathsf{S}_{+}^{(\rho_{5},\rho_{0})}(q)
 ={}&-2q^2 - 2q^3 - 2q^4 + q^6 + 6q^7 + 10q^8 + 15q^9 + 19q^{10} + 20q^{11} + 15q^{12}\\
&  + 3q^{13}- 13q^{14} - 47q^{15} - 81q^{16} - 129q^{17} - 176q^{18} - 227q^{19} - 261q^{20} + \cdots , \\
\mathsf{S}_{+}^{(\rho_{6},\rho_{0})}(q)
 ={}&2q^2 + q^3 - 5q^5 - 6q^6 - 12q^7 - 13q^8 - 12q^9 - 6q^{10} + 3q^{11} + 25q^{12} + 51q^{13}\\
& + 78q^{14} + 117q^{15} + 144q^{16} + 171q^{17} + 173q^{18} + 157q^{19} + 96q^{20} + \cdots , \\
\mathsf{S}_{+}^{(\rho_{7},\rho_{0})}(q)
 ={}&q + 2q^2 + q^3 + q^4 - 4q^5 - 7q^6 - 13q^7 - 15q^8 - 13q^9 - 9q^{10} + 5q^{11} + 26q^{12}\\
& + 56q^{13} + 86q^{14} + 132q^{15} + 161q^{16} + 193q^{17} + 196q^{18} + 182q^{19} \\
& + 114q^{20} + \cdots ,
\end{align*}
and
\begin{align*}
\mathsf{S}_{-}^{(\rho_{1},\rho_{0})}(q)
 ={}&-1 + q + 3q^2 + 2q^3 + q^4 - q^5 - 2q^6 - 10q^7 - 15q^8 - 21q^9 - 28q^{10} - 28q^{11} \\
&  - 20q^{12}- 5q^{13} + 22q^{14} + 65q^{15} + 116q^{16} + 174q^{17} + 239q^{18} + 304q^{19}\\
&  + 353q^{20} + \cdots , \\
\mathsf{S}_{-}^{(\rho_{2},\rho_{0})}(q)
 ={} &q - q^2 - 3q^3 - q^4 - 3q^5 + q^6 + q^7 + 6q^8 + 14q^9 + 19q^{10} + 27q^{11} + 31q^{12} \\
& + 28q^{13} + 18q^{14} - 37q^{16} - 77q^{17} - 137q^{18} - 200q^{19} - 272q^{20} + \cdots , \\
\mathsf{S}_{-}^{(\rho_{3},\rho_{0})}(q)
={}&
-1 + 2q + q^2 + q^3 - q^5 - 2q^6 - 6q^7 - 8q^8 - 13q^9 - 14q^{10} - 14q^{11} - 5q^{12} \\
&  - q^{13} + 20q^{14} + 39q^{15} + 69q^{16} + 96q^{17} + 133q^{18} + 157q^{19} + 184q^{20} + \cdots , \\
\mathsf{S}_{-}^{(\rho_{4},\rho_{0})}(q)
={}&
1 - 2q - q^2 - q^3 - q^4 + 2q^5 + 3q^6 + 10q^7 + 12q^8 + 17q^9 + 18q^{10} + 16q^{11} \\
&  + 3q^{12}- 11q^{13} - 41q^{14} - 76q^{15} - 117q^{16} - 159q^{17} - 204q^{18} - 232q^{19} \\
& - 249q^{20} + \cdots , \\
\mathsf{S}_{-}^{(\rho_{5},\rho_{0})}(q)
={}&
q - q^2 - 2q^3 - 3q^4 - 3q^5 - q^6 - q^7 + 6q^8 + 11q^9 + 21q^{10} + 31q^{11} + 40q^{12} \\
&  + 40q^{13} + 39q^{14} + 24q^{15} + 7q^{16} + 52q^{17} + 112q^{18} + 188q^{19} + 279q^{20} + \cdots , \\
\mathsf{S}_{-}^{(\rho_{6},\rho_{0})}(q)
={}&
1 - 3q - q^2 + 2q^3 + 4q^4 + 8q^5 + 6q^6 + 13q^7 + 8q^8 + 3q^9 - 10q^{10} - 29q^{11} \\
&  - 60q^{12}- 84q^{13} - 117q^{14} - 142q^{15} - 147q^{16} - 135q^{17} - 93q^{18} - 12q^{19}\\
& + 115q^{20} + \cdots , \\
\mathsf{S}_{-}^{(\rho_{7},\rho_{0})}(q)
={}&
1 - 3q - q^2 + 3q^3 + 2q^4 + 7q^5 + 6q^6 + 12q^7 + 7q^8 + 3q^9 - 6q^{10} - 26q^{11} \\
&  - 52q^{12}- 74q^{13} - 103q^{14} - 129q^{15} - 133q^{16} - 125q^{17} - 90q^{18} - 22q^{19} \\
& + 92q^{20} + \cdots .
\end{align*}
There are $7\times 8\times 2=112$ of these series and we include a PARI/GP~\cite{PARI} code to compute them.
\begin{lstlisting}
/* infinite Pochhammer symbol */
qpochinfty(a,q,N,ss=1,WMW=ss)=local(s,t,qn);t=1;s=t;qn=1;for(k=1,N,qn=qn*q^ss;t=-qn/q^(1/2+ss/2)*t/(1-qn)*a*q^(ss/2-WMW/2);s=s+t);s^ss;

/* Eisenstein series */
/* ss=/=0 gives the vector of Eisenstein series from 1 to n */
{Gf(q,n,N,ss)=
 local(qk,t,s,ex,bernt);
 ex=exp(x+O(x^(n+1)));qk=1;t=1;s=1;
 for(k=1,N,qk=qk*q;t=-t*qk/(1-qk)*ex;s=s+t);
 t=log(s/polcoeff(s,0,x));
 if(ss,bernt=bernvec(floor(n/2));
 vector(n,j,if(j==1,-1/4,0)+if(Mod(j,2)==Mod(0,2),-(-1)^j*bernt[j/2+1]/j/2,0)-j!*polcoeff(t,j,x)),
 -(-1)^n*bernfrac(n)/n/2-n!*polcoeff(t,n,x)
 )};

/* canonical basis */
{SI8x8(q)=
 [1/q^4, (q + 2)/q^2, (-q^3 + 2*q^2 - q - 1)/q^2, -q^3 - q^2 - 2*q - 1, q^4 + 2*q, q^6 + q^5 + q^3, -q^6 - q^4, -q^8;
 1/q^4, (q + 2)/q^2, (q^2 - q - 1)/q^2, -q^2 - 2*q - 2, q^3 + 2*q, q^5 + q^4 + q^3, -q^5 - q^3, -q^7;
 1/q^3, (q^2 + q + 1)/q^2, -1/q, -q^3 - 2*q^2 - q - 1, q^4 + 2*q^2, q^6 + q^5 + q^4, -q^6 - q^4, -q^8;
 0, (-q + 1)/q^3, (-q + 1)/q, (q^2 - 1)/q, -q + 1, -q^3 + q^2, 0, 0;
 1/q^4, (2*q + 1)/q^2, (q^3 - q - 1)/q^2, -q^3 - q^2 - 3*q, 2*q^2 + q, q^6 + q^5 + q^4 + q^3 - q^2, -q^5 - q^4, -q^8;
 1/q^4, (2*q + 1)/q^2, (q^3 + q^2 - 2*q - 1)/q^2, -q^3 - q^2 - 3*q, q^2 + 2*q, q^6 + q^5 + q^4 + q^3 - q^2, -q^5 - q^4, -q^8;
 0, 0, 0, -q + 1, 0, 0, 0, 0;
 1/q^4, (2*q + 1)/q^2, (q^3 - q - 1)/q^2, -q^3 - q^2 - 3*q, q^2 + 2*q, q^6 + q^5 + q^4 + q^3 - q^2, -q^5 - q^4, -q^8]/(1-q)};

/* zhat series computed inductively */
{Zhatf1(q,m,N,ss=1)=
 local(t0,t1,s,temp0,qk,qj);
 t0=(-1+q)*(ss+1);t1=1;qj=1;tj=1;
 for(j=1,N,qj=qj*q;tj=q*tj/(1-qj)^2;t1=t1+tj);
 t1=qpochinfty(q,q,N,ss,1)^2*t1;
 s=t1/q^m;qk=q;
 for(k=1,N,qk=qk/q;temp0=t1;
 t1=-(qk^2*t0+(q*qk + qk^3 - qk^2 + qk)*(1+qk)*(1-q*qk^2)*temp0)
 /(q^2*(1+qk)*(1-q*qk^2)*(1-qk^2/q^2)*(1-qk^2/q));
 t0=temp0;s=s+(qk/q^2)^m*t1
 );
 s};

{Zhatf2(q,m,N,ss=1)=
 local(t0,t1,tj,tjn,tjc,s,temp0,qk,qj);
 t0=0;qj=1;tj=1;tjc=-1/2+2*ss*Gf(q^ss,1,N);tjn=0;t1=tj*(tjc-2*tjn);
 for(j=1,N,qj=qj*q;tj=q*tj/(1-qj)^2;
 tjn=tjn+qj/(1-qj);t1=t1+tj*(tjc-2*tjn)
 );
 t1=qpochinfty(q,q,N,ss,1)^2*t1;
 s=t1/q^m;qk=q;
 for(k=1,N,qk=qk/q;temp0=t1;
 t1=-(qk^2*t0+(q*qk + qk^3 - qk^2 + qk)*(1+qk)*(1-q*qk^2)*temp0)
 /(q^2*(1+qk)*(1-q*qk^2)*(1-qk^2/q^2)*(1-qk^2/q));
 t0=temp0;s=s+(qk/q^2)^m*t1
 );
 s};

{Zhatf3(q,m,N,ss=1)=
 local(t0,t1,tj,tjn,tjc1,tjc2,tG1,tG2,s,temp0,qk,qj);
 t0=-2*(1-q);qj=1;tj=1;tG1=ss*Gf(q^ss,1,N);tG2=ss*Gf(q^ss,2,N);
 tjc1=(1/2-2*tG1);tjc2=(-2*m-4*tG1);tjn=0;
 t1=tj*((tjc1+2*tjn)*tjc2/2+2*tG2);
 for(j=1,N,qj=qj*q;tj=q*tj/(1-qj)^2;tjn=tjn+qj/(1-qj);
 tjc2=tjc2+2;t1=t1+tj*((tjc1+2*tjn)*tjc2/2+2*tG2)
 );
 t1=qpochinfty(q,q,N,ss,1)^2*t1;s=t1/q^m;qk=q;
 for(k=1,N,qk=qk/q;temp0=t1;
 t1=-(qk^2*t0+(q*qk + qk^3 - qk^2 + qk)*(1+qk)*(1-q*qk^2)*temp0)
 /(q^2*(1+qk)*(1-q*qk^2)*(1-qk^2/q^2)*(1-qk^2/q));
 t0=temp0;s=s+(qk/q^2)^m*t1
 );
 s};

{Zhatf4(q2,m,N,ss=1)=
 local(t0,t1,t0j,t1j,s,temp0,qk,qkj,qpc);
 qpc=qpochinfty(q2^2,q2^2,N,ss,1)*qpochinfty(q2^3,q2^2,N,ss);
 t0j=q2^(2)/qpoch(q2^2,q2^2,1);t1j=q2^(12)/qpoch(q2^3,q2^2,1)/qpoch(q2^2,q2^2,3);t0=t0j;t1=t1j;qj=1;
 for(j=1,N,qj=qj*q2^2;t0j=t0j*q2^3/(1-qj)/(1-qj*q2);
 t1j=t1j*q2^5/(1-qj)/(1-qj*q2^3);t0=t0+t0j;t1=t1+t1j
 );
 t0=qpc*t0;t1=qpc*t1;s=t0/q2^(3*m)+t1/q2^(5*m);qk=1;
 for(k=1,N,qk=qk/q2^2;temp0=t1;
 t1=-((qk/q2)^2*t0+(q2^2*(qk/q2) + (qk/q2)^3 - (qk/q2)^2 + (qk/q2))*(1+(qk/q2))*(1-q2^2*(qk/q2)^2)*temp0)
 /(q2^4*(1+(qk/q2))*(1-q2^2*(qk/q2)^2)*(1-(qk/q2)^2/q2^4)*(1-(qk/q2)^2/q2^2));
 t0=temp0;s=s+(qk/q2^5)^m*t1
 );
 s};

{Zhatf5(q,m,N,ss=1)=
 local(t0,t1,t0j,t1j,s,temp0,qk,qj,qpc);
 qpc=qpochinfty(q,q,N,ss,1)*qpochinfty(-q,q,N,ss);
 t0j=(-1)^m;t0=t0j;t1j=(-1)^m*q^(3)/qpoch(-q,q,1)/qpoch(q,q,2);
 t1=t1j;qj=1;
 for(j=1,N,qj=qj*q;t0j=-t0j*q/(1-qj^2);t0=t0+t0j;
 t1j=-t1j*q^2/(1-qj)/(1+qj*q);t1=t1+t1j
 );
 t0=qpc*t0;t1=qpc*t1;s=t0/q^m+t1/q^(2*m);qk=1;
 for(k=1,N,qk=qk/q;temp0=t1;
 t1=-((-qk)^2*t0+(q*(-qk) + (-qk)^3 - (-qk)^2 + (-qk))*(1+(-qk))*(1-q*(-qk)^2)*temp0)
 /(q^2*(1+(-qk))*(1-q*(-qk)^2)*(1-(-qk)^2/q^2)*(1-(-qk)^2/q));
 t0=temp0;s=s+(qk/q^2)^m*t1
 );
 s};

{Zhatf6(q2,m,N,ss=1)=
 local(t0,t1,t0j,t1j,s,temp0,qk,qj,qpc);
 qpc=qpochinfty(q2^2,q2^2,N,ss,1)*qpochinfty(-q2^3,q2^2,N,ss);
 t0j=(-1)^(m)*q2^2/qpoch(q2^2,q2^2,1);t0=t0j;
 t1j=(-1)^(m)*q2^12/qpoch(-q2^3,q2^2,1)/qpoch(q2^2,q2^2,3);
 t1=t1j;qj=1;
 for(j=1,N,qj=qj*q2^2;t0j=-t0j*q2^3/(1-qj)/(1+qj*q2);
 t0=t0+t0j;t1j=-t1j*q2^5/(1-qj)/(1+qj*q2^3);t1=t1+t1j
 );
 t0=qpc*t0;t1=qpc*t1;s=t0/q2^(3*m)+t1/q2^(5*m);qk=1;
 for(k=1,N,qk=qk/q2^2;temp0=t1;
 t1=-((-qk/q2)^2*t0+(q2^2*(-qk/q2) + (-qk/q2)^3 - (-qk/q2)^2 + (-qk/q2))*(1+(-qk/q2))*(1-q2^2*(-qk/q2)^2)*temp0)
 /(q2^4*(1+(-qk/q2))*(1-q2^2*(-qk/q2)^2)*(1-(-qk/q2)^2/q2^4)*(1-(-qk/q2)^2/q2^2));
 t0=temp0;s=s+(qk/q2^5)^m*t1
 );
 s};

{Zhatf7(q4,m,N,ss=1)=
 local(t0,t1,t0j,t1j,s,temp0,qk,qj,q22,qpc);
 q22=q4^2;qpc=qpochinfty(q22^2,q22^2,N,ss,1)*qpochinfty(q22,q22^2,N,ss);
 t0j=q4/qpoch(q22^2,q22^2,1);t0=t0j;
 t1j=q4^9/qpoch(q22,q22^2,-1)/qpoch(q22^2,q22^2,3);t1=t1j;qj=1;
 for(j=1,N,qj=qj*q22^2;t0j=t0j*q22^3/(1-qj)/(1-qj/q22);
 t0=t0+t0j;t1j=t1j*q22^5/(1-qj)/(1-qj/q22^3);t1=t1+t1j
 );
 t0=t0*qpc;t1=t1*qpc;s=t0/q22^(3*m)+t1/q22^(5*m);qk=1;
 for(k=1,N,qk=qk/q22^2;temp0=t1;
 t1=-((qk/q22)^2*t0+(q22^2*(qk/q22) + (qk/q22)^3 - (qk/q22)^2 + (qk/q22))*(1+(qk/q22))*(1-q22^2*(qk/q22)^2)*temp0)
 /(q22^4*(1+(qk/q22))*(1-q22^2*(qk/q22)^2)*(1-(qk/q22)^2/q22^4)*(1-(qk/q22)^2/q22^2));
 t0=temp0;s=s+(qk/q22^5)^m*t1
 );
 s};

{Zhatf8(q,m,N,ss=1)=
 local(t0,t1,t0j,t1j,s,temp0,qk,qj,qpc);
 qpc=qpochinfty(q,q,N,ss,1)*qpochinfty(-q,q,N,ss);
 t0j=(-1)^(m);t0=t0j;t1j=(-1)^(m)*q/qpoch(-q,q,-1)/qpoch(q,q,2);t1=t1j;qj=1;
 for(j=1,N,qj=qj*q;t0j=-t0j*q/(1-qj^2);t0=t0+t0j;
 t1j=-t1j*q^2/(1-qj)/(1+qj/q);t1=t1+t1j
 );
 t0=t0*qpc;t1=t1*qpc;s=t0/q^m+t1/q^(2*m);qk=1;
 for(k=1,N,qk=qk/q;temp0=t1;
 t1=-((-qk)^2*t0+(q*(-qk) + (-qk)^3 - (-qk)^2 + (-qk))*(1+(-qk))*(1-q*(-qk)^2)*temp0)
 /(q^2*(1+(-qk))*(1-q*(-qk)^2)*(1-(-qk)^2/q^2)*(1-(-qk)^2/q));
 t0=temp0;s=s+(qk/q^2)^m*t1
 );
 s};

{Zhatf9(q4,m,N,ss=1)=
 local(t0,t1,t0j,t1j,s,temp0,qk,q22,qpc);
 q22=q4^2;qpc=qpochinfty(q22^2,q22^2,N,ss,1)*qpochinfty(-q22,q22^2,N,ss);
 t0j=(-1)^(m)*q4/qpoch(q22^2,q22^2,1);t0=t0j;t1j=(-1)^(m)*q4^9/qpoch(-q22,q22^2,-1)/qpoch(q22^2,q22^2,3);t1=t1j;qj=1;
 for(j=1,N,qj=qj*q22^2;t0j=-t0j*q22^3/(1-qj)/(1+qj/q22);
 t0=t0+t0j;t1j=-t1j*q22^5/(1-qj)/(1+qj/q22^3);t1=t1+t1j
 );
 t0=t0*qpc;t1=t1*qpc;s=t0/q22^(3*m)+t1/q22^(5*m);qk=1;
 for(k=1,N,qk=qk/q22^2;temp0=t1;
 t1=-((-qk/q22)^2*t0+(q22^2*(-qk/q22) + (-qk/q22)^3 - (-qk/q22)^2 + (-qk/q22))*(1+(-qk/q22))*(1-q22^2*(-qk/q22)^2)*temp0)
 /(q22^4*(1+(-qk/q22))*(1-q22^2*(-qk/q22)^2)*(1-(-qk/q22)^2/q22^4)*(1-(-qk/q22)^2/q22^2));t0=temp0;s=s+(qk/q22^5)^m*t1
 );
 s};

/* 8x8 matrix of q-series */
{ZHAT8x8(q4,m,N,s=1)=
if(s==1,
[
Zhatf1(q4^4,m,N,s),Zhatf2(q4^4,m,N,s),Zhatf4(q4^2,m,N,s),Zhatf5(q4^4,m,N,s),Zhatf6(q4^2,m,N,s),Zhatf7(q4,m,N,s),Zhatf8(q4^4,m,N,s),Zhatf9(q4,m,N,s);
Zhatf1(q4^4,m+1,N,s),Zhatf2(q4^4,m+1,N,s),Zhatf4(q4^2,m+1,N,s),Zhatf5(q4^4,m+1,N,s),Zhatf6(q4^2,m+1,N,s),Zhatf7(q4,m+1,N,s),Zhatf8(q4^4,m+1,N,s),Zhatf9(q4,m+1,N,s);
Zhatf1(q4^4,m+2,N,s),Zhatf2(q4^4,m+2,N,s),Zhatf4(q4^2,m+2,N,s),Zhatf5(q4^4,m+2,N,s),Zhatf6(q4^2,m+2,N,s),Zhatf7(q4,m+2,N,s),Zhatf8(q4^4,m+2,N,s),Zhatf9(q4,m+2,N,s);
Zhatf1(q4^4,m+3,N,s),Zhatf2(q4^4,m+3,N,s),Zhatf4(q4^2,m+3,N,s),Zhatf5(q4^4,m+3,N,s),Zhatf6(q4^2,m+3,N,s),Zhatf7(q4,m+3,N,s),Zhatf8(q4^4,m+3,N,s),Zhatf9(q4,m+3,N,s);
Zhatf1(q4^4,m+4,N,s),Zhatf2(q4^4,m+4,N,s),Zhatf4(q4^2,m+4,N,s),Zhatf5(q4^4,m+4,N,s),Zhatf6(q4^2,m+4,N,s),Zhatf7(q4,m+4,N,s),Zhatf8(q4^4,m+4,N,s),Zhatf9(q4,m+4,N,s);
Zhatf1(q4^4,m+5,N,s),Zhatf2(q4^4,m+5,N,s),Zhatf4(q4^2,m+5,N,s),Zhatf5(q4^4,m+5,N,s),Zhatf6(q4^2,m+5,N,s),Zhatf7(q4,m+5,N,s),Zhatf8(q4^4,m+5,N,s),Zhatf9(q4,m+5,N,s);
Zhatf1(q4^4,m+6,N,s),Zhatf2(q4^4,m+6,N,s),Zhatf4(q4^2,m+6,N,s),Zhatf5(q4^4,m+6,N,s),Zhatf6(q4^2,m+6,N,s),Zhatf7(q4,m+6,N,s),Zhatf8(q4^4,m+6,N,s),Zhatf9(q4,m+6,N,s);
Zhatf1(q4^4,m+7,N,s),Zhatf2(q4^4,m+7,N,s),Zhatf4(q4^2,m+7,N,s),Zhatf5(q4^4,m+7,N,s),Zhatf6(q4^2,m+7,N,s),Zhatf7(q4,m+7,N,s),Zhatf8(q4^4,m+7,N,s),Zhatf9(q4,m+7,N,s)
],
[
Zhatf3(q4^4,m,N,s),Zhatf1(q4^4,m,N,s),Zhatf4(q4^2,m,N,s),Zhatf5(q4^4,m,N,s),Zhatf6(q4^2,m,N,s),Zhatf7(q4,m,N,s),Zhatf8(q4^4,m,N,s),Zhatf9(q4,m,N,s);
Zhatf3(q4^4,m+1,N,s),Zhatf1(q4^4,m+1,N,s),Zhatf4(q4^2,m+1,N,s),Zhatf5(q4^4,m+1,N,s),Zhatf6(q4^2,m+1,N,s),Zhatf7(q4,m+1,N,s),Zhatf8(q4^4,m+1,N,s),Zhatf9(q4,m+1,N,s);
Zhatf3(q4^4,m+2,N,s),Zhatf1(q4^4,m+2,N,s),Zhatf4(q4^2,m+2,N,s),Zhatf5(q4^4,m+2,N,s),Zhatf6(q4^2,m+2,N,s),Zhatf7(q4,m+2,N,s),Zhatf8(q4^4,m+2,N,s),Zhatf9(q4,m+2,N,s);
Zhatf3(q4^4,m+3,N,s),Zhatf1(q4^4,m+3,N,s),Zhatf4(q4^2,m+3,N,s),Zhatf5(q4^4,m+3,N,s),Zhatf6(q4^2,m+3,N,s),Zhatf7(q4,m+3,N,s),Zhatf8(q4^4,m+3,N,s),Zhatf9(q4,m+3,N,s);
Zhatf3(q4^4,m+4,N,s),Zhatf1(q4^4,m+4,N,s),Zhatf4(q4^2,m+4,N,s),Zhatf5(q4^4,m+4,N,s),Zhatf6(q4^2,m+4,N,s),Zhatf7(q4,m+4,N,s),Zhatf8(q4^4,m+4,N,s),Zhatf9(q4,m+4,N,s);
Zhatf3(q4^4,m+5,N,s),Zhatf1(q4^4,m+5,N,s),Zhatf4(q4^2,m+5,N,s),Zhatf5(q4^4,m+5,N,s),Zhatf6(q4^2,m+5,N,s),Zhatf7(q4,m+5,N,s),Zhatf8(q4^4,m+5,N,s),Zhatf9(q4,m+5,N,s);
Zhatf3(q4^4,m+6,N,s),Zhatf1(q4^4,m+6,N,s),Zhatf4(q4^2,m+6,N,s),Zhatf5(q4^4,m+6,N,s),Zhatf6(q4^2,m+6,N,s),Zhatf7(q4,m+6,N,s),Zhatf8(q4^4,m+6,N,s),Zhatf9(q4,m+6,N,s);
Zhatf3(q4^4,m+7,N,s),Zhatf1(q4^4,m+7,N,s),Zhatf4(q4^2,m+7,N,s),Zhatf5(q4^4,m+7,N,s),Zhatf6(q4^2,m+7,N,s),Zhatf7(q4,m+7,N,s),Zhatf8(q4^4,m+7,N,s),Zhatf9(q4,m+7,N,s)])};

/* Stokes matrix in region I */
SMI8x8(q4,N)=[-1,0,0,0,0,0,0,0; 0,-1,0,0,0,0,0,0; 0,0,-1,0,0,0,0,0; 0,0,0,-1,0,0,0,0; 0,0,0,0,-1,0,0,0; 0,0,0,0,0,-1,0,0; 0,0,0,0,0,0,1,0; 0,0,0,0,0,0,0,1]*SI8x8(q4^4)*ZHAT8x8(q4,-3,N)*matdiagonal([-1,-1,1,1,-1,-1,-1,1])

/* Stokes matrix in region II */
SMII8x8(q4,N)=[1,0,0,0,0,0,0,0; 0,1,0,0,0,0,0,0; 0,0,1,0,0,0,0,0; 0,0,0,1,0,0,0,0; 0,0,0,0,1,0,0,0; 0,0,0,0,0,1,0,0; 0,0,0,0,0,0,0,1; 0,0,0,0,0,0,1,0]*SI8x8(q4^4)*ZHAT8x8(q4,-3,N)*matdiagonal([1,-1,-1,-1,1,1,1,-1])

SMII8x8(q4+O(q4^150),150)*SMI8x8(q4+O(q4^150),150)^(-1)-matrix(8,8,i,j,if(i==j,1,0)+O(q4^10))
/*
[O(q4^10) O(q4^10) O(q4^10) O(q4^10) O(q4^10) O(q4^10) O(q4^10) O(q4^10)]

[-1 + q4^4 + 3*q4^8 + O(q4^10) -q4^4 - 2*q4^8 + O(q4^10) 1 + q4^4 + O(q4^10) -q4^8 + O(q4^10) q4^4 + 2*q4^8 + O(q4^10) -1 - q4^4 + q4^8 + O(q4^10) -1 + 3*q4^8 + O(q4^10) -1 + 3*q4^8 + O(q4^10)]

[q4^4 - q4^8 + O(q4^10) q4^8 + O(q4^10) -q4^4 - q4^8 + O(q4^10) O(q4^10) -q4^8 + O(q4^10) q4^4 + q4^8 + O(q4^10) q4^4 + O(q4^10) q4^4 + O(q4^10)]

[-1 + 2*q4^4 + q4^8 + O(q4^10) -q4^4 - q4^8 + O(q4^10) 1 + O(q4^10) -q4^8 + O(q4^10) q4^4 + q4^8 + O(q4^10) -1 + q4^8 + O(q4^10) -1 + q4^4 + 2*q4^8 + O(q4^10) -1 + q4^4 + 2*q4^8 + O(q4^10)]

[1 - 2*q4^4 - q4^8 + O(q4^10) q4^4 + q4^8 + O(q4^10) -1 + O(q4^10) q4^8 + O(q4^10) -q4^4 - q4^8 + O(q4^10) 1 - q4^8 + O(q4^10) 1 - q4^4 - 2*q4^8 + O(q4^10) 1 - q4^4 - 2*q4^8 + O(q4^10)]

[-q4^4 + q4^8 + O(q4^10) -q4^8 + O(q4^10) q4^4 + q4^8 + O(q4^10) O(q4^10) q4^8 + O(q4^10) -q4^4 - q4^8 + O(q4^10) -q4^4 + O(q4^10) -q4^4 + O(q4^10)]

[1 - 3*q4^4 - q4^8 + O(q4^10) q4^4 + O(q4^10) -1 + q4^4 + 2*q4^8 + O(q4^10) q4^8 + O(q4^10) -q4^4 + O(q4^10) 1 - q4^4 - 3*q4^8 + O(q4^10) -2*q4^4 - 3*q4^8 + O(q4^10) -2*q4^4 - 3*q4^8 + O(q4^10)]

[1 - 3*q4^4 - q4^8 + O(q4^10) q4^4 + O(q4^10) -1 + q4^4 + 2*q4^8 + O(q4^10) q4^8 + O(q4^10) -q4^4 + O(q4^10) 1 - q4^4 - 3*q4^8 + O(q4^10) -2*q4^4 - 3*q4^8 + O(q4^10) -2*q4^4 - 3*q4^8 + O(q4^10)]
*/

/* Stokes matrix in region IV */
SMIV8x8(q4,N)=[-1,0,0,0,0,0,0,0; 0,1,0,0,0,0,0,0; 0,0,1,0,0,0,0,0; 0,0,0,1,0,0,0,0; 0,0,0,0,1,0,0,0; 0,0,0,0,0,1,0,0; 0,0,0,0,0,0,-1,0; 0,0,0,0,0,0,0,-1]*SI8x8(1/q4^4)*ZHAT8x8(1/q4,-3,N,-1)*matdiagonal([-1,-1,1,1,-1,-1,-1,1])

/* Stokes matrix in region III */
SMIII8x8(q4,N)=[1,0,0,0,0,0,0,0; 0,1,0,0,0,0,0,0; 0,0,1,0,0,0,0,0; 0,0,0,1,0,0,0,0; 0,0,0,0,1,0,0,0; 0,0,0,0,0,1,0,0; 0,0,0,0,0,0,0,1; 0,0,0,0,0,0,1,0]*SI8x8(1/q4^4)*ZHAT8x8(1/q4,-3,N,-1)*matdiagonal([1,1,1,1,-1,-1,-1,1])

SMIV8x8(q4+O(q4^100),100)*SMIII8x8(q4+O(q4^100),100)^(-1)-matrix(8,8,i,j,if(i==j,1,0)+O(q4^10))
/*
[O(q4^10) O(q4^10) O(q4^10) O(q4^10) O(q4^10) O(q4^10) O(q4^10) O(q4^10)]

[O(q4^10) -q4^4 - 2*q4^8 + O(q4^10) q4^8 + O(q4^10) -q4^4 - q4^8 + O(q4^10) q4^4 + q4^8 + O(q4^10) -q4^8 + O(q4^10) -q4^4 + O(q4^10) -q4^4 + O(q4^10)]

[-1 + O(q4^10) 1 + q4^4 + O(q4^10) -q4^4 - q4^8 + O(q4^10) 1 + O(q4^10) -1 + O(q4^10) q4^4 + q4^8 + O(q4^10) 1 - q4^4 - 2*q4^8 + O(q4^10) 1 - q4^4 - 2*q4^8 + O(q4^10)]

[-q4^4 + q4^8 + O(q4^10) -q4^8 + O(q4^10) O(q4^10) -q4^8 + O(q4^10) q4^8 + O(q4^10) O(q4^10) -q4^8 + O(q4^10) -q4^8 + O(q4^10)]

[-q4^4 - q4^8 + O(q4^10) q4^4 + 2*q4^8 + O(q4^10) -q4^8 + O(q4^10) q4^4 + q4^8 + O(q4^10) -q4^4 - q4^8 + O(q4^10) q4^8 + O(q4^10) q4^4 + O(q4^10) q4^4 + O(q4^10)]

[-2*q4^8 + O(q4^10) -1 - q4^4 + q4^8 + O(q4^10) q4^4 + q4^8 + O(q4^10) -1 + q4^8 + O(q4^10) 1 - q4^8 + O(q4^10) -q4^4 - q4^8 + O(q4^10) -1 + q4^4 + 3*q4^8 + O(q4^10) -1 + q4^4 + 3*q4^8 + O(q4^10)]

[2*q4^8 + O(q4^10) 1 - 3*q4^8 + O(q4^10) -q4^4 + O(q4^10) 1 - q4^4 - 2*q4^8 + O(q4^10) -1 + q4^4 + 2*q4^8 + O(q4^10) q4^4 + O(q4^10) -2*q4^4 - 3*q4^8 + O(q4^10) -2*q4^4 - 3*q4^8 + O(q4^10)]

[q4^4 + 2*q4^8 + O(q4^10) 1 - 3*q4^8 + O(q4^10) -q4^4 + O(q4^10) 1 - q4^4 - 2*q4^8 + O(q4^10) -1 + q4^4 + 2*q4^8 + O(q4^10) q4^4 + O(q4^10) -2*q4^4 - 3*q4^8 + O(q4^10) -2*q4^4 - 3*q4^8 + O(q4^10)]
*/
\end{lstlisting}

\subsection*{Acknowledgements}

This work was carried out throughout my doctoral studies and therefore owes a great debt to my supervisors Stavros Garoufalidis and Don Zagier. Their constant sharing of ideas and perspectives has been invaluable to the completion of this work. I would also like to thank Dongmin Gang, Mauricio Romo, and Masahito Yamazaki for sharing some of their code early on in this project. Again I would like to thank Dongmin Gang for sharing a computation of the 3d index of $4_1(-1,2)$ proposed in~\cite{Gang3d} and the question of whether this could be found in the Stokes constants. I thank Jie Gu, Marcos Mari\~no and Matthias Storzer for various conversations throughout this project and Sergei Gukov for suggestions improving the manuscript. Finally, I would like to thank the referees for many helpful suggestions. The work of the author has been supported by the Max-Planck-Gesellschaft with the Max Planck institute for mathematics in Bonn and the Southern University of Science and Technology’s International Center for Mathematics in Shenzhen.

\LastPageEnding

\end{document}